\documentclass[final,3p]{elsarticle}

\usepackage[T1]{fontenc}
\usepackage{lmodern}
\usepackage{amssymb}
\usepackage{amsmath}
\usepackage{amsthm}
\usepackage{xcolor}
\usepackage{url}
\usepackage{booktabs}
\usepackage{hyperref}
\usepackage{multirow}
\usepackage{caption}
\usepackage{subcaption}
\usepackage{natbib}
\usepackage{algorithm}
\usepackage{algpseudocode}
\usepackage{float}
\usepackage{graphicx}
\graphicspath{{figures/}}

\usepackage[section]{placeins}

\newtheorem{theorem}{Theorem}
\newtheorem{lemma}{Lemma}
\newtheorem{corollary}{Corollary}
\newtheorem{proposition}{Proposition}
\newtheorem{definition}{Definition}
\newtheorem{remark}{Remark}
\newtheorem{assumption}{Assumption}

\makeatletter
\renewcommand{\theALG@line}{\arabic{ALG@line}}
\makeatother

\newcommand{\jump}[1]{[\![\,#1\,]\!]}
\newcommand{\R}{\mathbb{R}}
\newcommand{\dd}{\mathrm{d}}
\newcommand{\Rm}{\mathcal{R}^{\mathrm{m}}}
\newcommand{\Rp}{\mathcal{R}^{\mathrm{p}}}
\newcommand{\Sm}{\mathcal{S}}
\newcommand{\Lc}{\mathcal{L}}
\newcommand{\naSub}{n/a{\hypersetup{pdfborder={0 0 0}}\hyperlink{note:na-sub}{\textcolor{blue}{$^{\diamond}$}}}}
\newcommand{\naIface}{n/a{\hypersetup{pdfborder={0 0 0}}\hyperlink{note:na-iface}{\textcolor{blue}{$^{\ast}$}}}}
\newcommand{\naSubNote}{\protect\hypertarget{note:na-sub}{\textcolor{blue}{$^{\diamond}$}}}
\newcommand{\naIfaceNote}{\protect\hypertarget{note:na-iface}{\textcolor{blue}{$^{\ast}$}}}

\journal{Elsevier}
\date{}

\begin{document}

\begin{frontmatter}

\title{Well-Balanced WE-PINNs for the Shallow Water Equations over Smooth and Discontinuous Bathymetries}

\author[a,b]{Ismail Oubarka\corref{cor1}}
\ead{ismail.oubarka@um6p.ma, oubarka.ismail5@gmail.com}
\cortext[cor1]{Corresponding author}

\author[a]{Imad Kissami}
\author[b]{Mohamed Boubekeur}
\author[a,c]{Imad Elmahi}
\author[b]{Fayssal Benkhaldoun}

\affiliation[a]{organization={Mohammed VI Polytechnic University, College of Computing},
            addressline={Lot 660},
            city={Ben Guerir},
            postcode={43150},
            country={Morocco}}

\affiliation[b]{organization={LAGA, CNRS, UMR 7539, Université Sorbonne Paris Nord},
            city={Villetaneuse},
            postcode={93430},
            country={France}}

\affiliation[c]{organization={ENSAO, LMCS, Mohammed I University},
               addressline={B.P. 669},
               city={Oujda},
               postcode={60000},
               country={Morocco}}

\begin{abstract}
Physics-informed neural networks (PINNs) solve partial differential equations
without a mesh, but they do not in general preserve the steady states that
well-balanced schemes reproduce exactly. We present a well-balanced extension of
the weak and entropy PINN (WE-PINN) framework to the shallow water equations
over smooth and discontinuous bathymetry. The residual is a weak space--time
control-volume identity enforced on randomly placed control volumes that may
overlap, so no mesh is needed. In free-surface variables, the geometric source
is integrated by parts onto the bed, so that the network is never
differentiated, and a cell-local bed-slope correction makes the quadrature
reproduce the exact bed increment. Bed discontinuities, located from the
bathymetry alone, separate sub-domains carrying independent networks, coupled
through conservation of discharge and total head. We prove that the fully
discrete residual vanishes identically on the lake at rest for any
control-volume family, quadrature order and nominal bed slope, so that the
equilibrium is an exact global minimiser of the training objective; that the
continuous residual characterises weak solutions; that the slope correction is
consistent; and that no approximation with a continuous free surface can satisfy
the head condition at a step with non-zero discharge, so that the decomposition
is necessary. Numerically, the residual vanishes to round-off on four
lake-at-rest configurations over smooth, discontinuous, compound and tabulated
beds, and the trained network reaches the equilibrium to round-off under an
extended training budget. On seven Riemann problems at a bottom step, covering
sub-, super- and critical transitions, the network attains relative $L^1$
free-surface errors between $6.9\times10^{-4}$ and $3.5\times10^{-3}$, and
reaches critical states without prescribing them.
\end{abstract}


\end{frontmatter}

\section{Introduction}\label{sec:introduction}

The shallow water equations with variable bottom topography are the canonical
hyperbolic \emph{balance} law of computational hydraulics. Unlike the
homogeneous conservation laws for which most shock-capturing methods were
designed, they admit non-trivial steady states in which the flux gradient and
the geometric source cancel exactly while neither vanishes separately. The most
important of these is the lake at rest, a horizontal free surface over
motionless water. A discretisation that fails to reproduce this cancellation
generates spurious waves whose amplitude scales with the truncation error of the
bed treatment and which, over realistic bathymetries, can exceed the physical
signal. Methods that reproduce it exactly are called well balanced, or said to
satisfy the exact conservation property ($C$-property)
\cite{BermudezVazquez1994}, and the requirement has shaped the numerical
literature on this system since its early analyses
\cite{GreenbergLeRoux1996,LeVeque1998}.

In the finite-volume and discontinuous Galerkin setting, two mechanisms achieve
this property. The first modifies the reconstruction so that the discrete flux
and source balance by construction, as in the hydrostatic reconstruction
\cite{Audusse2004} and in central-upwind schemes
\cite{KurganovLevy2002,KurganovPetrova2007}. The second changes variables so
that the equilibrium becomes trivial: the surface-gradient method
\cite{Zhou2001} and the pre-balanced formulation \cite{Rogers2003} reconstruct
the free-surface elevation rather than the depth, after which the lake at rest
is a constant state. High-order extensions of both routes exist
\cite{XingShu2005,Noelle2006}, together with path-conservative formulations
\cite{Pares2006,Castro2007,Gallardo2007}, and well-balanced schemes have since
been extended to \emph{moving} equilibria, characterised by constancy of the
discharge and of the total head \cite{NoelleXingShu2007,CaleffiValiani2017}.
Surveys are given in \cite{XingShu2014review,Bouchut2004}. Two features common
to these methods do not transfer to a neural formulation: the equilibrium is
preserved by an explicit update rule acting on cell averages, and the mechanism
that preserves it is tied to the mesh through the reconstruction operator. The
second feature has been addressed in \cite{BihloMacLachlan2018}, which gives a
necessary and sufficient criterion for exact preservation of the lake at rest in
mesh-based and meshless schemes alike, by pairing a mimetic design of the
discrete derivative operators with a compatible averaging of the depth in the
source. That criterion is formulated for the strong form, requires a discrete
bed derivative and is developed for smooth topography.

The problem is more delicate when the bed itself is discontinuous. The source
$-gh\,\partial_x b$ is then a measure supported on the bed jump, multiplied by a
solution that is generally discontinuous at the same location, and the resulting
nonconservative product is not determined by the equations; its definition
requires additional structure, such as a family of paths \cite{DalMaso1995}.
Conservation of discharge and of total head at the step was adopted in
\cite{AlcrudoBenkhaldoun2001}, where the wet-bed solutions of the Riemann
problem are classified into five configurations and closed-form solutions are
derived for a complete family of test problems. A momentum-based closure, which
yields different intermediate states, leads to the exact solutions of
\cite{Bernetti2008}; the conditions under which the momentum balance can be
closed at a discontinuous bed are examined in \cite{ValianiCaleffi2017}, and the
problem is studied within the theory of nonconservative systems in
\cite{LeFlochThanh2007}. These analyses identify regimes in which the solution is
not unique, and the system is resonant where a characteristic speed approaches
the vanishing speed of the stationary wave
\cite{IsaacsonTemple1995,GoatinLeFloch2004,Andrianov2005}. Any method for this
problem must therefore select a closure explicitly, and it inherits whatever
multiplicity that closure carries. Wetting and drying constitutes a further body
of work \cite{Brufau2002,Bokhove2005,Gallardo2007,Bunya2009,XingZhangShu2010};
the present construction assumes a strictly wet bed and does not address it.

Physics-informed neural networks (PINNs)~\cite{Raissi2019} replace the marching
update by the minimisation of a residual functional over the parameters of a
network defined on the whole space--time domain. For hyperbolic problems the
standard formulation is structurally inadequate: the strong-form residual
presupposes pointwise derivatives that do not exist across a shock and diverges
as a smooth approximation sharpens towards a discontinuous entropy solution,
whereas suitable weak residual norms remain bounded \cite{de2024wpinns}.
Moreover, the resulting optimisation problem does not favour the entropy
solution \cite{CIPINN2024}, and standard formulations smear discontinuities and
misplace shocks on canonical benchmarks
\cite{fuks2020limitations,krishnapriyan2021characterizing}. Three families of
remedies have emerged. Stabilised strong-form residuals (training points
clustered near discontinuities \cite{Mao2020}, localised artificial viscosity
\cite{coutinho2023physics}, gradient annihilation \cite{ferrer2024gradient},
lift-and-embed formulations \cite{sun2026lift}, explicit enforcement of jump
conditions \cite{liu2024discontinuity}, interface flux continuity in a domain
decomposition \cite{jagtap2020conservative}) retain pointwise differential
information, and with it a sensitivity to hyperparameters and, in some cases, to
prior knowledge of where the discontinuities lie. Finite-volume-informed
objectives \cite{Patel2022,cassia2025godunov,patsatzis2025gorinns,SFVnet2026}
embed a discrete update rule, an approximate Riemann solver or a mesh-dependent
stabilisation into the loss, and thereby reintroduce the discretisation choices
that the neural representation was meant to avoid; in a well-balanced Roe
finite-volume loss for the shallow water equations \cite{LiuFVPINN2026}, the
balance is inherited from the embedded scheme, at the price of an unstructured
mesh and of a loss landscape on which purely physics-based training fails unless
sparse measurements are supplied. Weak and dual formulations
\cite{de2024wpinns,Chaumet2024,wang2025wf,CIPINN2024} avoid differentiating the
approximation but introduce auxiliary test-function networks, potentials or
saddle-point problems.

Neural solvers for the shallow water system with topography are recent, and the
discrete equilibrium structure has received uneven attention. Label-free
networks have been shown to be competitive with a finite-volume model on the
accuracy--cost trade-off \cite{Qi2024}, in a comparison that does not address
that structure. The two-dimensional equations with topographic and rainfall
sources have been solved in strong form over smooth beds \cite{Tian2025}, with
well-balancedness treated as an empirical training difficulty. To address
well-balancedness directly, \cite{Brecht2025} builds it into the network ansatz
through a small-weight initialisation, which keeps the solution near the
equilibrium rather than exactly at it; the network uses the depth and velocity
rather than the free surface, and all topographies are smooth. For discontinuous beds, \cite{Dazzi2024} adds the bed
elevation as a network output and obtains its slope by automatic differentiation,
which makes the free surface continuous; by Theorem~\ref{thm:no_continuous}
below, such a free surface cannot satisfy the head condition at a bed step with
non-zero discharge.

Across this body of work, the equilibrium is approached by training, imposed
approximately through the ansatz, or inherited from an embedded mesh-based
scheme, but in no case is it satisfied exactly by the residual itself; and
wherever the bed-slope source is treated explicitly, the bed is differentiated,
an operation without meaning where it jumps. The gap is structural rather than
one of measured accuracy.

The WE-PINN framework~\cite{Oubarka2026WEPINNs}, on which we build, belongs to
the weak family and uses no auxiliary networks. It enforces the space--time
divergence identity over sampled control volumes, evaluates the resulting
boundary integrals by Gauss--Legendre quadrature, imposes the initial condition
exactly through a temporal ansatz, supplements the conservation residual by
integral Kru\v{z}kov entropy inequalities and a total-variation control, and
trains with Adam followed by L-BFGS. It is, however, built on an identity for a
homogeneous law and on a globally continuous approximation, and both hypotheses
fail here. We therefore separate what is inherited, what is adapted and what is
new. Inherited are the control-volume residual, the face quadrature, the exact
initial condition and the training schedule, together with the use of
overlapping control volumes. Here these are placed at random and form no
partition of the space--time domain, so that the method is mesh-free: no mesh is
generated, and no property of the method depends on one. Adapted from shallow
water theory are the free-surface change of variables, which plays the role it
plays in surface-gradient and pre-balanced schemes \cite{Zhou2001,Rogers2003},
and the energy closure at a bed step \cite{AlcrudoBenkhaldoun2001}, of whose two
conditions one is forced by mass conservation and the other is a modelling
choice.

New are three constructions and their analysis. The first is the evaluation of
the geometric source over a space--time control volume by means of the
elementary Leibniz identity
$b\,\partial_x\eta=\partial_x(b\,\eta)-\eta\,\partial_xb$, an integration by
parts that transfers the derivative from the solution to the bed; the residual
therefore never differentiates the network and remains defined for
$\eta\in L^\infty$. The second is a cell-local shift of the bed slope, defined so
that the quadrature reproduces the exact bed increment across the cell, which
makes the lake-at-rest cancellation exact for any control-volume family, any
quadrature order and any nominal bed slope. The principle behind it, that the
source must be discretised consistently with the flux so that the two balance
exactly at equilibrium, is the classical principle of well-balanced schemes
\cite{BermudezVazquez1994,LeVeque1998}; what is new is its realisation inside a
quadrature rule on arbitrary space--time cells. In contrast with the criterion
of \cite{BihloMacLachlan2018}, which constrains discrete derivative operators in
the strong form and presupposes a smooth bed and a bed derivative, the shift acts
inside a quadrature rule in the weak form, requires no bed derivative, holds for
an arbitrary and even inconsistent nominal slope, tolerates a discontinuous $b$,
and admits into the residual no information about the bathymetry beyond the
values of $b$ at the two lateral faces of the cell. Well-balancedness is thereby
decoupled from the mesh and from the accuracy, and indeed the existence, of a bed
derivative. The third construction is a decomposition of the domain into
sub-domains of smooth bed, located from the bathymetry alone, which are
parameterised independently and coupled through the stationary step transition,
so that the free surface can jump exactly where the exact solution requires.
Unlike the domain decompositions of cPINN \cite{jagtap2020conservative} and
XPINN \cite{jagtap2020xpinn}, whose interfaces are chosen freely and coupled by
continuity of the flux or of the solution, the interfaces here are dictated by
the bed discontinuities and coupled by the hydraulic step closure, across which
the free surface jumps.

On the theoretical side, we prove that the common zeros of the continuous
residuals are the weak solutions of the balance law; that the fully discrete
residual vanishes identically on the lake at rest, so that the equilibrium is an
exactly representable global minimiser of the training objective; that the
bed-slope correction is of the order of the quadrature error it cancels; that no
approximation with a continuous free surface can satisfy the head condition at a
bed step at non-zero discharge, with an explicit lower bound on the defect; and
that the formulation reduces exactly to the original one when the bed is
constant. Classical results used along the way (the free-surface form of the
equations, the conjugate-depth structure at a step, integration by parts for
functions of bounded variation, Gauss--Legendre error bounds, floating-point
error analysis and the momentum balance across a step) are stated with their
sources, and their proofs, where given, are included only to keep the paper
self-contained. The limitations of the analysis, notably the absence of a
convergence estimate for the trained network and of a mechanism selecting
among non-unique step transitions, are discussed in Section~\ref{sec:open}.

The numerical experiments test the analysis directly. On four lake-at-rest
configurations, over smooth, discontinuous, compound and tabulated beds, the
exact equilibrium is fed through the residual operator, which annihilates it
bitwise or to a round-off floor identified by its dependence on the working
precision; the trained network approaches the equilibrium to a level set by the
optimisation budget alone. Seven Riemann problems at a bottom step with
closed-form solutions, covering sub-, super- and critical transitions and both
signs of the free-surface jump, show that the trained network recovers the
stationary transition and the surrounding wave pattern, and a dam break over a
smooth bump confirms that without a bed jump the decomposition and interface
terms vanish and only the source treatment remains.

The paper is organised as follows. Section~\ref{sec:governing} sets out the
system, its equilibria, the closure at a bed discontinuity and the standing
assumptions. Section~\ref{sec:method} develops the formulation,
Section~\ref{sec:theory} its analysis, and Section~\ref{sec:implementation} the
implementation details required for reproduction. Section~\ref{sec:results}
reports the numerical experiments and Section~\ref{sec:conclusion} concludes.
The appendices give the tabulated bathymetry of one test case, the exact Riemann
solver at a bottom step and the well-balanced finite-volume reference.

\section{Governing equations and mathematical model}
\label{sec:governing}

\subsection{The shallow water equations in free-surface variables}
\label{sec:eta_form}

Let $\Omega_x=(x_{\mathrm L},x_{\mathrm R})$ with $L:=|\Omega_x|$, let $T>0$
and set $\Omega=\Omega_x\times(0,T)$. Over a fixed bed
$b:\overline{\Omega_x}\to\R$ the one-dimensional shallow water equations read
\begin{equation}
\partial_t h+\partial_x q=0,
\qquad
\partial_t q+\partial_x\!\Big(\frac{q^{2}}{h}+\tfrac12 g h^{2}\Big)
=-\,g\,h\,\partial_x b ,
\label{eq:swe_conservative}
\end{equation}
with depth $h>0$, discharge $q=hu$, velocity $u$ and gravitational acceleration
$g$; $c=\sqrt{gh}$ is the celerity and $F=u/c$ the Froude number. Throughout,
the free-surface elevation $\eta:=h+b$ replaces the depth as the first unknown.

The following reformulation is classical \cite{Zhou2001,Rogers2003}; we record
it with its short proof because the construction relies on its exact form.

\begin{lemma}[Free-surface form~\cite{Zhou2001,Rogers2003}]
\label{lem:eta_form}
Let $b\in C^{1}(\Omega_x)$ and let $(h,q)$ be a classical solution of
\eqref{eq:swe_conservative} with $h>0$. Then $(\eta,q)$, $\eta=h+b$, is a
classical solution of
\begin{equation}
\partial_t\eta+\partial_x q=0,
\qquad
\partial_t q+\partial_x\Phi=g\,b\,\partial_x\eta,
\qquad
\Phi:=\frac{q^{2}}{h}+\tfrac12 g\,\eta^{2},\quad h=\eta-b ,
\label{eq:swe_eta}
\end{equation}
and conversely.
\end{lemma}

\begin{proof}
Since $b$ is independent of time, $\partial_t\eta=\partial_th$ and the mass
equations coincide. Substituting $h=\eta-b$ in the hydrostatic term gives
\[
\partial_x\big(\tfrac12 g h^{2}\big)
=\partial_x\big(\tfrac12 g\eta^{2}\big)
-g\,b\,\partial_x\eta-g\,\eta\,\partial_x b+g\,b\,\partial_x b ,
\qquad
-g\,h\,\partial_x b=-g\,\eta\,\partial_x b+g\,b\,\partial_x b .
\]
The terms $-g\eta\,\partial_xb$ and $gb\,\partial_xb$ appear on both sides of
the momentum equation and cancel, leaving the second equation of
\eqref{eq:swe_eta}. Each step is an algebraic identity between $C^1$ functions
and can be read in either direction.
\end{proof}

This change of variables underlies the surface-gradient and pre-balanced
finite-volume formulations. The momentum flux $\Phi$
and the right-hand side $g\,b\,\partial_x\eta$ of \eqref{eq:swe_eta} are the
two quantities on which the construction of Section~\ref{sec:method} operates.

\subsection{Steady states and the lake at rest}
\label{sec:steady}

As is well known \cite{Chow1959,XingShu2014review}, setting
$\partial_t\equiv0$ in \eqref{eq:swe_eta} characterises smooth steady states by
constancy of the discharge $q$ and of the total head
\begin{equation}
H:=\eta+\frac{u^{2}}{2g} .
\label{eq:head}
\end{equation}
Among them, the lake at rest
\begin{equation}
q\equiv0,\qquad \eta\equiv\eta^\star>\max_{\overline{\Omega_x}}b
\label{eq:lake_at_rest}
\end{equation}
is the state whose exact preservation defines the $C$-property. The two
formulations realise \eqref{eq:lake_at_rest} differently. In the conservative
variables, $\partial_x(\tfrac12gh^{2})=-gh\,\partial_xb\neq0$ wherever
$\partial_xb\neq0$, so the flux gradient and the source are individually
non-zero and cancel only in sum; a discretisation must reproduce that
cancellation. In the free-surface variables,
\begin{equation}
\Phi\equiv\tfrac12 g(\eta^\star)^{2},
\qquad
\partial_x\Phi\equiv0,
\qquad
g\,b\,\partial_x\eta\equiv0 ,
\label{eq:lar_pointwise}
\end{equation}
so both sides of the momentum equation vanish separately and pointwise, for
every bed, with no cancellation between large quantities; this observation is
the basis of the surface-gradient and pre-balanced schemes
\cite{Zhou2001,Rogers2003}. Preserving
\eqref{eq:lake_at_rest} is therefore a matter of not destroying a structure the
continuous equations already possess. The reformulation has one cost: the
right-hand side of \eqref{eq:swe_eta} is not a flux divergence, so the
space--time divergence theorem does not convert it into a boundary integral.
That conversion is precisely the hypothesis on which the inherited WE-PINN
residual rests, and Section~\ref{sec:source} restores the missing structure.

\subsection{The stationary step transition}
\label{sec:sst}

Let $b$ be discontinuous at $x_0\in\Omega_x$. For a function $v$ with one-sided
limits at $x_0$ we write $v_\pm:=v(x_0^\pm)$ and $\jump{v}:=v_+-v_-$, and assume
$\jump{b}\neq0$. The exact solution then carries a stationary wave pinned at
$x_0$. The two conditions satisfied by this wave have different logical status.

The mass equation of \eqref{eq:swe_eta} is in divergence form and contains no
product with $b$; it therefore holds across $x_0$ in $\mathcal D'(\Omega)$, and
its Rankine--Hugoniot condition \cite{LeVeque2002} for a discontinuity of speed
$s$ is $-s\jump{\eta}+\jump{q}=0$. The wave is pinned to the bed, so $s=0$ and
\begin{equation}
\jump{q}=0 .
\label{eq:sst_q}
\end{equation}
Condition \eqref{eq:sst_q} is a consequence of the equations. The momentum
equation supplies no counterpart: its right-hand side is the nonconservative
product $g\,b\,\partial_x\eta$, in which a function discontinuous at $x_0$
multiplies a measure charging $x_0$, and such a product is not determined by
\eqref{eq:swe_eta} \cite{DalMaso1995}. A closure must be prescribed. Following
\cite{AlcrudoBenkhaldoun2001} we impose conservation of the total head,
\begin{equation}
\jump{H}=\jump{\;\eta+\frac{u^{2}}{2g}\;}=0 ,
\label{eq:sst_H}
\end{equation}
which, since $\eta$ carries the bed elevation, is the classical
energy-conserving (Bernoulli) condition of open-channel hydraulics
\cite{Chow1959,Henderson1966},
$h_-+b_-+u_-^{2}/(2g)=h_++b_++u_+^{2}/(2g)$. We call
\eqref{eq:sst_q}--\eqref{eq:sst_H} the \emph{stationary step transition}. It
does not conserve momentum; the defect, recalled in
Proposition~\ref{prop:momentum_defect}, equals the horizontal force
exerted by the step on the flow, which the one-dimensional equations do not determine
\cite{Bernetti2008,ValianiCaleffi2017}. The
momentum-based closure of \cite{Bernetti2008} is equally admissible and
produces different intermediate states; in the formulation below the closure
enters through a single residual, so replacing \eqref{eq:sst_H} modifies one
term of the objective and nothing else. In resonant regimes the closure carries
its own non-uniqueness \cite{IsaacsonTemple1995,GoatinLeFloch2004}, which the
method inherits.

Eliminating $u_+=q/h_+$ between \eqref{eq:sst_q} and \eqref{eq:sst_H} reduces
the transition to a scalar equation for the conjugate depth. The following
lemma is the classical specific-energy analysis of open-channel hydraulics
\cite{Chow1959,Henderson1966}, used in the same form in
\cite{AlcrudoBenkhaldoun2001}; we include a short proof for completeness.

\begin{lemma}[Conjugate depths at a step~\cite{Chow1959,AlcrudoBenkhaldoun2001}]
\label{lem:conjugate}
Let $q\neq0$, let $H_-$ be the head of the left trace and set $A:=H_--b_+$. A
depth $h_+>0$ satisfies \eqref{eq:sst_q}--\eqref{eq:sst_H} if and only if
\begin{equation}
f(h_+)=0,
\qquad
f(h):=h^{3}-A\,h^{2}+\frac{q^{2}}{2g} .
\label{eq:sst_cubic}
\end{equation}
Assume $A>0$. Then $f$ has two positive roots, counted with multiplicity, if
and only if
\begin{equation}
\frac{q^{2}}{2gA^{3}}\ \le\ \frac{4}{27} ,
\label{eq:sst_admissible}
\end{equation}
and none otherwise. When \eqref{eq:sst_admissible} is strict, the roots satisfy
$0<h^{\mathrm{sup}}<\tfrac23A<h^{\mathrm{sub}}<A$, and the Froude number of
the corresponding state is
\begin{equation}
F^{2}=\frac{q^{2}}{g\,h^{3}}=\frac{2(A-h)}{h} ,
\label{eq:froude_root}
\end{equation}
so that $h^{\mathrm{sup}}$ is supercritical and $h^{\mathrm{sub}}$
subcritical. Equality in \eqref{eq:sst_admissible} holds exactly when the roots
coalesce at the critical depth $\tfrac23A=(q^{2}/g)^{1/3}$.
\end{lemma}

\begin{proof}
With $u_+=q/h_+$, \eqref{eq:sst_H} reads $h_+-A+q^{2}/(2gh_+^{2})=0$, and
multiplication by $h_+^{2}>0$ gives \eqref{eq:sst_cubic}. Since
$f'(h)=h(3h-2A)$, $f$ has a local maximum $f(0)=q^{2}/(2g)>0$ and a local
minimum $f(\tfrac23A)=q^{2}/(2g)-\tfrac{4}{27}A^{3}$; the product of the three
roots is $-q^{2}/(2g)<0$, so exactly one root is negative, and $f$ has two
positive roots if and only if $f(\tfrac23A)\le0$, which is
\eqref{eq:sst_admissible}. If the inequality is strict, $f(0)>0$,
$f(\tfrac23A)<0$ and $f(A)=q^{2}/(2g)>0$ locate the roots. At a root,
$q^{2}/(2g)=h^{2}(A-h)$, whence \eqref{eq:froude_root}, which exceeds one if
and only if $h<\tfrac23A$. At equality, $\tfrac23A$ is a double root and
$q^{2}=2g\cdot\tfrac49A^{2}\cdot\tfrac13A=gh^{3}$, i.e.\ $F=1$.
\end{proof}

Combining Lemma~\ref{lem:conjugate} with the wave curves of the homogeneous
system yields the five wet-bed configurations of the Riemann problem at a step
classified in \cite{AlcrudoBenkhaldoun2001}. For the construction that follows,
the essential consequence of \eqref{eq:sst_q}--\eqref{eq:sst_H} is that only
$q$ and $H$ are continuous across $x_0$, whereas
$\jump{\eta}=\jump{h}+\jump{b}\neq0$ in general.
Theorem~\ref{thm:no_continuous} converts this observation into a quantitative
obstruction for approximation spaces with a continuous free surface.

\subsection{Standing assumptions and notation}
\label{sec:assumptions}

The bed is allowed finitely many jump discontinuities
$x_{\mathrm L}=e_0<e_1<\dots<e_M<e_{M+1}=x_{\mathrm R}$ and is smooth between
them. We set
\[
\Omega_k:=(e_k,e_{k+1}),\qquad {\Omega_k^{T}}:=\Omega_k\times(0,T),\qquad k=0,\dots,M .
\]
The jump abscissae are not supplied to the method; Section~\ref{sec:detection}
recovers them from $b$.

\begin{assumption}[Bathymetry]
\label{ass:bed}
$b|_{\Omega_k}$ extends to a function of $C^{1}(\overline{\Omega_k})$ for every
$k$, and $b$ is evaluated exactly at every node used below.
\end{assumption}

\begin{assumption}[Admissible states]
\label{ass:states}
$(\eta,q)\in L^\infty({\Omega_k^{T}})^{2}$ and $h:=\eta-b\ge h_{\min}>0$ almost everywhere
on ${\Omega_k^{T}}$, for every $k$.
\end{assumption}

The positive constant $h_{\min}$ is the depth floor of
Section~\ref{sec:interface}. Under Assumptions~\ref{ass:bed}--\ref{ass:states},
$\Phi\in L^\infty({\Omega_k^{T}})$. The implementation evaluates the momentum flux through
the regularised depth $h_\delta:=\max\{\eta-b,h_{\min}\}$, which equals $h$
under Assumption~\ref{ass:states}, so the analysis concerns the residual that is
actually assembled. Dry states ($h\to0$) are excluded.

The admissible control volumes on $\Omega_k$ are the axis-aligned rectangles
\begin{equation}
\mathcal K_k:=\big\{K=[x_1,x_2]\times[t_1,t_2]\ :\
[x_1,x_2]\subset\overline{\Omega_k},\ [t_1,t_2]\subset[0,T],\
x_1<x_2,\ t_1<t_2\big\},
\label{eq:cell_set}
\end{equation}
parameterised by $(x_1,x_2,t_1,t_2)\in\R^{4}$; a property holds \emph{for
almost every} $K\in\mathcal K_k$ when it holds for almost every parameter with
respect to four-dimensional Lebesgue measure. For $K\in\mathcal K_k$ we write
$\Delta x=x_2-x_1$, $\Delta t=t_2-t_1$ and $|K|=\Delta x\,\Delta t$. These are
the extents of the single control volume $K$, not the spacing of a mesh: the
control volumes of Section~\ref{sec:cells} are placed at random, overlap, and
form no structured grid, so that $\Delta x$ and $\Delta t$ are attributes of
each cell and no neighbouring cell, node spacing or grid index is implied. The
Gauss--Legendre nodes and weights on $[-1,1]$ are $\{\xi_p,w_p\}_{p=1}^{Q}$,
normalised by
\begin{equation}
\sum_{p=1}^{Q}w_p=2 ,
\label{eq:weights_sum}
\end{equation}
and their images on $K$ are
$x_p=\tfrac12(x_1+x_2)+\tfrac{\Delta x}{2}\xi_p$ and
$t_p=\tfrac12(t_1+t_2)+\tfrac{\Delta t}{2}\xi_p$.

\section{Method}
\label{sec:method}

\subsection{The inherited weak space--time residual}
\label{sec:wepinns}

For a homogeneous system $\partial_t U+\partial_x F(U)=0$, the space--time flux
$J(U):=(F(U),U)$ turns the law into $\nabla_{x,t}\cdot J(U)=0$, whose integral
form over an axis-aligned rectangle
$D=[x_1,x_2]\times[t_1,t_2]\subset\overline\Omega$ is the classical integral
conservation law \cite{LeVeque2002,godlewski2013numerical}
\begin{equation}
\int_{x_1}^{x_2}\!\big[U(x,t_2)-U(x,t_1)\big]\dd x
+\int_{t_1}^{t_2}\!\big[F(U(x_2,t))-F(U(x_1,t))\big]\dd t=0 .
\label{eq:we_rect}
\end{equation}
Three properties of \eqref{eq:we_rect} are used unchanged. It involves only the
traces of $U$ and $F(U)$ on the four faces of $D$, so it is meaningful for
$U\in L^\infty$ and requires no derivative of the approximation; it remains
bounded when $U$ is discontinuous and encodes the Rankine--Hugoniot conditions
implicitly; and it holds on every bounded Lipschitz region, so a family of test
volumes need not form a partition of $\Omega$. WE-PINNs~\cite{Oubarka2026WEPINNs}
penalise the volume-weighted mean square of the left-hand side of
\eqref{eq:we_rect} over a sampled family of control volumes, evaluate the face
integrals by Gauss--Legendre quadrature, and impose the initial condition exactly
through a temporal ansatz $U_\theta=(1-\tau)U_0+\tau\mathcal N_\theta$ with
$\tau(0)=0$.

Three hypotheses of that construction fail for \eqref{eq:swe_eta}. The law is
not homogeneous, so \eqref{eq:we_rect} is not available; the ansatz is globally
continuous in $x$, whereas Section~\ref{sec:sst} requires $\jump{\eta}\neq0$ at
a bed step; and the entropy machinery has no established counterpart for a
stationary wave at a bed discontinuity. The first two failures are addressed in
Sections~\ref{sec:source}--\ref{sec:correction} and
Section~\ref{sec:detection}, respectively. The entropy and total-variation
terms of the original framework are not used in the formulation analysed here.

The input of the method is $\Omega_x$, $T$, the bed $b$ and the initial conditions
$\eta_0,u_0$, each supplied as a function that can be evaluated pointwise. No
information about the solution is used: the number and position of the bed
discontinuities, the domain decomposition, the normalisation scales and the
boundary treatment are all derived from these inputs.

\subsection{Locating the bed discontinuities and decomposing the domain}
\label{sec:detection}

The discontinuities of $b$ are a property of the data, not of the flow, so
locating them requires no knowledge of the solution; the operation is the
analogue of a finite-volume code evaluating $b$ at cell interfaces. With
$\mathrm{osc}(b):=\max_{\overline{\Omega_x}}b-\min_{\overline{\Omega_x}}b$, set
\begin{equation}
\mathrm{tol}:=\max\big\{10^{-4}\max\{\mathrm{osc}(b),1\},\,10^{-12}\big\}.
\label{eq:tol}
\end{equation}
The detector scans $b$ on a uniform grid, flags every interval across which the
increment of $b$ exceeds $\mathrm{tol}$, bisects it to a candidate abscissa
$x_0$, and confirms the candidate by comparing the symmetric increments
$j_k:=b(x_0+d_k)-b(x_0-d_k)$ at the probe distances $d_1=10^{-7}L$ and
$d_2=10^{-2}d_1$, accepting when
\begin{equation}
|j_2|\ge\tfrac12|j_1|
\qquad\text{and}\qquad
|j_2|\ge\mathrm{tol}.
\label{eq:detect_test}
\end{equation}
At a genuine jump $j_k\to\jump{b}\neq0$, so $j_2/j_1\to1$ and both inequalities
hold. If instead $b$ is Lipschitz with constant $L_b$ near $x_0$, then
$|j_2|\le2L_bd_2=2L_b\times10^{-9}L$ and the second inequality fails whenever
$2L_bd_2<\mathrm{tol}$, that is, unless the slope exceeds
$\mathrm{tol}/(2d_2)$. The confirmation test is thus exact in the limit
$d_1\to0$; the
bracketing stage that precedes it can still miss a genuine jump immediately
preceded by a steep smooth feature, a situation in which the well-balanced
property nevertheless survives (Remark~\ref{rem:wb_weak}).

The accepted abscissae define the sub-domains $\Omega_k$ of
Section~\ref{sec:assumptions}. Each $\Omega_k$ carries its own network, so the
global approximation is smooth inside every sub-domain and admits a jump at
every interior node $e_m$. A single total budget $N$ of control volumes ---
a number of space--time cells, not a product of counts in $x$ and in $t$ --- is
shared among the sub-domains in proportion to their widths,
\begin{equation}
N_k=\max\Big\{8,\ \Big[\,N\,\frac{|\Omega_k|}{L}\,\Big]\Big\} ,
\label{eq:cells_share}
\end{equation}
where $[\,\cdot\,]$ denotes rounding to the nearest integer. Each sub-domain is
flagged constant-bed or variable-bed by sampling $b$ at interior points; the
flag decides whether the volume term of the source is assembled. All bed and
initial-condition values attached to $\Omega_k$ are evaluated at abscissae pulled a
distance $10^{-9}|\Omega_k|$ inside $\Omega_k$, which makes every trace
one-sided at a bed jump; the shift is applied to the bed and the initial conditions, which genuinely
jump, and never to the network input.

\subsection{Space--time control volumes}
\label{sec:cells}

\begin{figure}[H]
  \centering
  \includegraphics[width=0.95\textwidth]{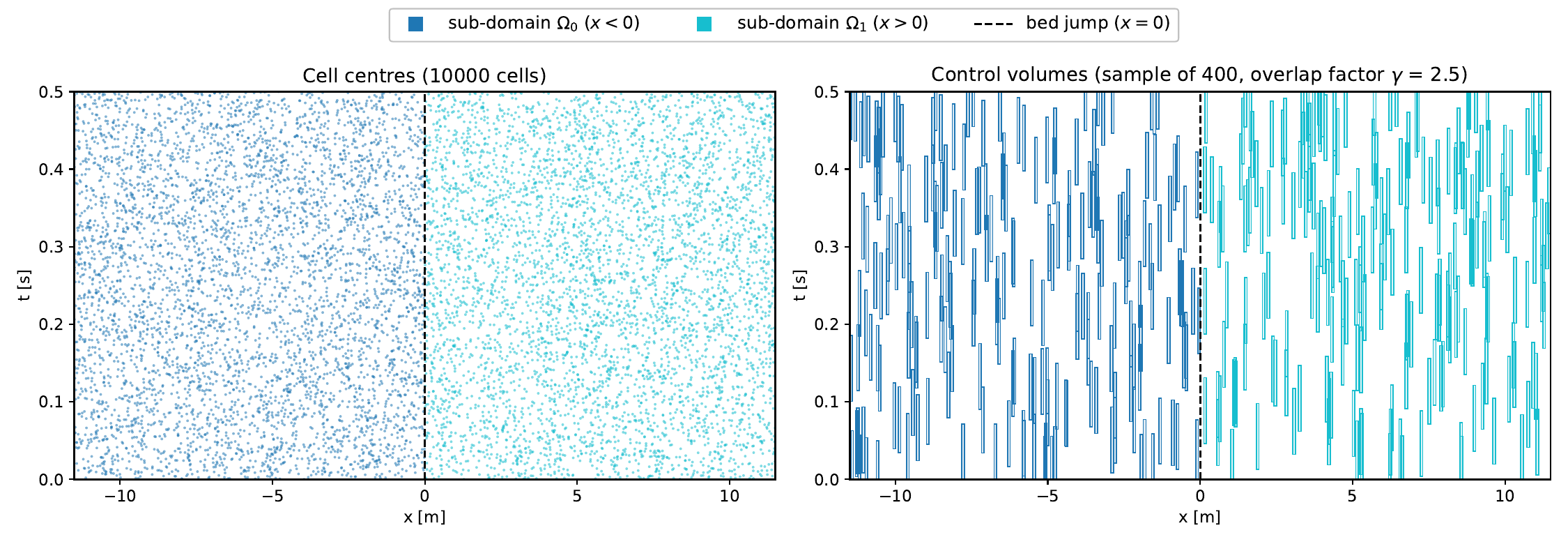}
  \caption{An unstructured, overlapping space--time control-volume family for a
  domain with a single bed jump: cell centres (left) and a sample of the control
  volumes (right), coloured by sub-domain. The dashed line marks the bed jump;
  no cell straddles it. The family is generated once from the geometry alone,
  by placing the cell centres independently and uniformly at random.}
  \label{fig:cells}
\end{figure}

On each sub-domain the control volumes are axis-aligned space--time rectangles
of equal size centred on randomly placed points. With $L_x:=|\Omega_k|$ and
$L_t:=T$, the centres are
\begin{equation}
x_i^{c}=e_k+L_x\,U_i,
\qquad
t_i^{c}=L_t\,V_i,
\qquad i=0,\dots,N_k-1,
\label{eq:cell_centres}
\end{equation}
where $U_0,\dots,U_{N_k-1},V_0,\dots,V_{N_k-1}$ are independent random
variables, uniformly distributed on $[0,1]$ and drawn once with a fixed seed.
The half-widths are
\begin{equation}
h_x=\min\Big\{\frac{L_x}{n_x^{\mathrm{eq}}}\frac{\gamma}{2},\frac{L_x}{2}\Big\},
\qquad
h_t=\min\Big\{\frac{L_t}{n_t^{\mathrm{eq}}}\frac{\gamma}{2\alpha_{\mathrm c}},\frac{L_t}{2}\Big\},
\qquad
n_x^{\mathrm{eq}}=\sqrt{\tfrac{N_kL_x}{L_t}},\quad
n_t^{\mathrm{eq}}=\sqrt{\tfrac{N_kL_t}{L_x}},
\label{eq:cell_sizes}
\end{equation}
with overlap factor $\gamma>0$ and aspect ratio $\alpha_{\mathrm c}>0$, and the
cell $K_i=[x_1^i,x_2^i]\times[t_1^i,t_2^i]$ is obtained by clipping
$x_i^{c}\mp h_x$ to $\overline{\Omega_k}$ and $t_i^{c}\mp h_t$ to $[0,T]$.
The counts $n_x^{\mathrm{eq}}$ and $n_t^{\mathrm{eq}}$ are those of an
\emph{equivalent} lattice with $N_k$ cells and serve only to fix the cell size;
no such lattice is built. Every unclipped cell therefore has the same extents
$\Delta x=2h_x$, equal to $\gamma\sqrt{L_xL_t/N_k}$ when the minimum in
\eqref{eq:cell_sizes} is not active, and $\Delta t=2h_t$, clipped cells near
$\partial\Omega_k$ or $t\in\{0,T\}$ being smaller. In what follows $\Delta x$
always denotes this cell width and never a grid spacing, and refining means
increasing $N_k$, which shrinks the cells as $N_k^{-1/2}$ while their centres
remain random. For
$\gamma>1$ the cells overlap and, after clipping, may leave parts of ${\Omega_k^{T}}$
uncovered. Both are admissible, since \eqref{eq:we_rect} holds on every bounded
Lipschitz region: the family is a Monte Carlo collection of test regions, not a
mesh (Fig.~\ref{fig:cells}). The random placement is compatible with the
analysis of Section~\ref{sec:theory}: the well-balanced property is proved cell
by cell, and consistency is a statement about all rectangles. Coverage must
nevertheless be monitored, because the weak identity is enforced only where a
cell reaches. The family is constructed once from the geometry alone and is
never adapted to the solution, which permits the pre-computation of
Section~\ref{sec:implementation}. Clipping guarantees
\begin{equation}
K_i\in\mathcal K_k,
\qquad\text{hence}\qquad
b\in C^{1}(K_i)
\label{eq:cell_inclusion}
\end{equation}
under Assumption~\ref{ass:bed}, which is the hypothesis under which the
integration by parts of Section~\ref{sec:source} is legitimate.

\begin{remark}[Cell shape]
\label{rem:cell_shape}
Nothing in the formulation requires rectangles: any bounded Lipschitz region
may serve as a test region. Rectangles are used because the boundary splits
into four axis-aligned faces carrying one-dimensional Gauss--Legendre rules,
because clipping delivers \eqref{eq:cell_inclusion} immediately, and because the
bed-slope correction takes an elementary form, the two horizontal faces sharing
the abscissae $x_p$ and the two lateral faces the instants $t_p$. For another
shape, the boundary quadrature must reproduce $\oint_{\partial K}n\,\dd s=0$
exactly for Theorem~\ref{thm:wb} to survive, and the shift
\eqref{eq:bp_correction} must be redefined accordingly.
\end{remark}

\subsection{A source term that never differentiates the network}
\label{sec:source}

Let $K=[x_1,x_2]\times[t_1,t_2]\in\mathcal K_k$ and let $(\eta,q)$ satisfy
Assumption~\ref{ass:states}. Integrating \eqref{eq:swe_eta} over $K$ and
applying the divergence theorem to the conservative part yields the weak
residuals
\begin{align}
\Rm(K)&=\int_{x_1}^{x_2}\!\big[\eta(x,t_2)-\eta(x,t_1)\big]\dd x
+\int_{t_1}^{t_2}\!\big[q(x_2,t)-q(x_1,t)\big]\dd t ,
\label{eq:res_mass}\\
\Rp(K)&=\int_{x_1}^{x_2}\!\big[q(x,t_2)-q(x,t_1)\big]\dd x
+\int_{t_1}^{t_2}\!\big[\Phi(x_2,t)-\Phi(x_1,t)\big]\dd t-\Sm(K) ,
\label{eq:res_mom}
\end{align}
in which $\Sm(K)$ represents $\int_Kg\,b\,\partial_x\eta$. As it stands, this
integral requires a derivative of the solution. We remove that requirement by
integrating by parts in $x$ (the elementary Leibniz identity), which is
legitimate because $b\in C^{1}(K)$ by
\eqref{eq:cell_inclusion}, and which transfers the derivative from the solution
onto the bed:
\begin{equation}
\Sm(K):=
g\!\int_{t_1}^{t_2}\!\big[b(x_2)\eta(x_2,t)-b(x_1)\eta(x_1,t)\big]\dd t
-g\!\int_{t_1}^{t_2}\!\!\int_{x_1}^{x_2}\!\eta(x,t)\,b'(x)\,\dd x\,\dd t .
\label{eq:source_ibp}
\end{equation}
The right-hand side of \eqref{eq:source_ibp} is well defined for
$\eta\in L^\infty(K)$, and Lemma~\ref{lem:source_pairing} identifies it with the
pairing of $b$ against the distributional derivative $D_x\eta$, which is
meaningful because $b$ is continuous on the cell even where $\eta$ is not. We
take \eqref{eq:source_ibp} as the definition of the source contribution.
Consequently the residual pair \eqref{eq:res_mass}--\eqref{eq:res_mom} never
differentiates the approximation and remains a statement about point values of
$\eta$ and $q$, whether or not the solution is discontinuous inside $K$; by
Proposition~\ref{prop:consistency}, its common zeros over $\mathcal K_k$ are
exactly the weak solutions of \eqref{eq:swe_eta} on ${\Omega_k^{T}}$. When $b\equiv b_0$
is constant on $\Omega_k$, the volume term vanishes identically and is not
assembled, while the face term reduces to
$g\,b_0\int_{t_1}^{t_2}[\eta(x_2,t)-\eta(x_1,t)]\dd t$.

\subsection{The quadrature-consistent bed-slope correction}
\label{sec:correction}

This is the element that makes the method well balanced. It applies to a
quadrature rule on arbitrary space--time cells the classical requirement that
the discrete source balance the discrete flux exactly at equilibrium
\cite{BermudezVazquez1994,LeVeque1998,Noelle2006}. The face integrals of
\eqref{eq:source_ibp} are evaluated with the $Q$-point rule already used for
the fluxes, and the volume term with the tensor-product $Q\times Q$ rule:
\begin{equation}
\Sm_Q(K)=
\frac{g\Delta t}{2}\sum_{p}w_p\big[b(x_2)\eta(x_2,t_p)-b(x_1)\eta(x_1,t_p)\big]
-\frac{g\Delta x\Delta t}{4}\sum_{p,r}w_pw_r\,\eta(x_p,t_r)\,\tilde b'(x_p) .
\label{eq:source_quad}
\end{equation}
At the lake at rest, $\eta\equiv\eta^\star$, the first term collapses to
$g\eta^\star\Delta t\,(b(x_2)-b(x_1))$, the \emph{exact} bed increment across
the cell, while the second collapses to
$g\eta^\star\Delta t\,\frac{\Delta x}{2}\sum_pw_p\,\tilde b'(x_p)$, a
\emph{quadrature} of the same increment. If $\tilde b'$ were the nominal bed
slope, the two would differ by the quadrature error and the $C$-property would
be lost. The remedy is to make the quadrature exact for this one integrand by
construction. Let $b'(x_p)$ denote any available approximation of the bed
derivative at the nodes, let
$G_Q(K)=\frac{\Delta x}{2}\sum_pw_p\,b'(x_p)$ be the corresponding
approximation of $\int_{x_1}^{x_2}b'$, and use the shifted slope
\begin{equation}
\tilde b'(x_p):=b'(x_p)+\frac{\big(b(x_2)-b(x_1)\big)-G_Q(K)}{\Delta x} ,
\label{eq:bp_correction}
\end{equation}
a shift by a constant on each cell, which by \eqref{eq:weights_sum} satisfies
\begin{equation}
\frac{\Delta x}{2}\sum_{p}w_p\,\tilde b'(x_p)=b(x_2)-b(x_1) .
\label{eq:correction_property}
\end{equation}
The two terms of \eqref{eq:source_quad} then cancel identically at the
equilibrium. Identity \eqref{eq:correction_property} is the single algebraic
fact behind Theorem~\ref{thm:wb}, and three of its features deserve emphasis.
First, the shift is defined \emph{from} the computed value $G_Q(K)$, so the
cancellation is exact whatever $b'$ is, accurate or not, and indeed whether or
not $b$ is differentiable on the cell. Second, since $\Delta x$ and the nodes
are cell-dependent, \eqref{eq:bp_correction} is defined for overlapping,
non-uniform and unstructured families exactly as for a mesh. Third, the only
information about the bathymetry entering the equilibrium cancellation is the
pair of face values $b(x_1)$, $b(x_2)$. Consistency is unaffected: the shift is
of the order of the quadrature error it cancels when $b$ is smooth, and equals
$J/\Delta x+O(1)$ when $b$ carries a jump of height $J$ inside the cell
(Propositions~\ref{prop:shift_order}--\ref{prop:source_consistency}).

\subsection{Coupling at a bed discontinuity, boundary conditions and positivity}
\label{sec:interface}

At each detected jump $e_m$ the traces of the two adjacent networks are coupled
through the stationary step transition. Writing $b_\mp=b(e_m^\mp)$,
$\eta_\mp,q_\mp$ for the corresponding traces and
$u_\mp=q_\mp/\max\{\eta_\mp-b_\mp,h_{\min}\}$, the interface residuals are
\begin{equation}
r_q(t)=\tau(t)\big(q_--q_+\big),
\qquad
r_H(t)=\tau(t)\Big[\Big(\eta_-+\frac{u_-^{2}}{2g}\Big)
-\Big(\eta_++\frac{u_+^{2}}{2g}\Big)\Big] ,
\label{eq:iface_res}
\end{equation}
evaluated at $N_\Gamma Q$ Gauss points in time, $N_\Gamma$ being a prescribed
number of time intervals on the interface, independent of the control-volume
family and used for the outer boundary conditions as well. The factor $\tau$ is
the temporal blending of \eqref{eq:ansatz}, with $\tau(0)=0$ and $\tau>0$ for
$t>0$. It deactivates the conditions at $t=0$, where they are incompatible with
a Riemann initial condition with $q_l\neq q_r$; \eqref{eq:sst_q}--\eqref{eq:sst_H} are in
any case statements about $t>0$. Both residuals vanish identically at the lake
at rest, so the coupling does not interfere with the $C$-property.

At each outer boundary one scalar condition is imposed, matching the number of
incoming characteristics for subcritical flow. The transmissive treatment
is the standard characteristic condition \cite{toro2009riemann}: it freezes
the incoming Riemann invariant at the value carried by the initial condition,
\begin{equation}
\big(u+2\sqrt{gh}\big)(x_{\mathrm L},t)=W_{\mathrm L},
\qquad
\big(u-2\sqrt{gh}\big)(x_{\mathrm R},t)=W_{\mathrm R} ,
\label{eq:bc_transmissive}
\end{equation}
with $W_{\mathrm L}=u_0(x_{\mathrm L})+2\sqrt{gh_0(x_{\mathrm L})}$ and
$W_{\mathrm R}=u_0(x_{\mathrm R})-2\sqrt{gh_0(x_{\mathrm R})}$. This is exact
for a Riemann problem truncated to a finite window as long as no wave has
reached the boundary, and remains exact after a rarefaction has reached it,
since $u+2c$ is constant through a $1$-rarefaction; it degrades only if a shock
crosses the boundary. The wall treatment imposes $q=0$ at both ends, and a
Dirichlet treatment pinning both components to the initial condition is also
available. Positivity is encouraged by the one-sided quadratic barrier
\begin{equation}
\Lc_{\mathrm P}
=\Big\langle\big(\max\{0,\;h_{\mathrm{safe}}-h_\theta\}\big)^{2}\Big\rangle ,
\label{eq:barrier}
\end{equation}
averaged over the cell faces and the interface nodes, and the hard floor
$h_{\min}<h_{\mathrm{safe}}$ is applied inside $q^{2}/h$ and $q/h$ so that the
objective remains finite. On wet beds satisfying Assumption~\ref{ass:states}
with $h\ge h_{\mathrm{safe}}$, both mechanisms are inactive.

\subsection{Approximation and training objective}
\label{sec:loss}

On $\Omega_k$ the solution is represented by
\begin{equation}
\eta_{\theta_k}=\big(1-\tau(t)\big)\eta_0+\tau(t)\,\mathcal N^{(1)}_{\theta_k},
\qquad
q_{\theta_k}=\big(1-\tau(t)\big)q_0+\tau(t)\,\mathcal N^{(2)}_{\theta_k},
\qquad
\tau(t)=\tanh(\alpha t),\ \ \alpha=\frac{c_\tau}{T},
\label{eq:ansatz}
\end{equation}
where $\mathcal N_{\theta_k}=(\mathcal N^{(1)}_{\theta_k},\mathcal N^{(2)}_{\theta_k})$
is a feedforward network acting on the affinely rescaled coordinates of
$\overline{{\Omega_k^{T}}}$ and $q_0=(\eta_0-b)u_0$. Since $\tau(0)=0$, the initial
condition holds exactly and pointwise, including for discontinuous Riemann
initial conditions, so no initial-condition term appears in the objective; since
$\eta_0,q_0$ are read one-sidedly, the ansatz is discontinuous at a bed jump
exactly where the initial conditions are. The constant $c_\tau>0$ concentrates the blending
over the initial fraction $2/c_\tau$ of the window, so that the network is free
to develop the solution early; it is an input and is not tied to the control
volumes, which carry no layer structure in time.

The two cell residuals have different physical dimensions and may differ by
orders of magnitude, so they are normalised by scales built from the initial
conditions alone. With $\eta_{\mathrm{sc}}=\max_{\overline{\Omega_x}}|\eta_0|$ and
$u_{\mathrm{sc}}=\sqrt{g\max_{\overline{\Omega_x}}h_0}$,
\begin{equation}
S_{\mathrm m}=\frac{\eta_{\mathrm{sc}}u_{\mathrm{sc}}}{L},\quad
S_{\mathrm p}=\frac{g\eta_{\mathrm{sc}}^{2}}{L},\quad
S_q=\eta_{\mathrm{sc}}u_{\mathrm{sc}},\quad
S_H=\eta_{\mathrm{sc}},\quad
S_W=u_{\mathrm{sc}} .
\label{eq:scales}
\end{equation}
Dividing a cell residual by $|K|S$ compares it with the size of its constituent
terms, uniformly over cells of different sizes. Denoting by
$\{K_i^k\}_{i=1}^{N_k}$ the cells of $\Omega_k$ and by $\Rm_Q,\Rp_Q$ the
quadrature residuals of Definition~\ref{def:disc_res}, the objective is
\begin{equation}
\Lc(\theta)=\Lc_{\mathrm m}+\Lc_{\mathrm p}
+\lambda_{\mathrm I}\Lc_{\mathrm I}
+\lambda_{\mathrm B}\Lc_{\mathrm B}
+\lambda_{\mathrm P}\Lc_{\mathrm P} ,
\label{eq:total_loss}
\end{equation}
with
\begin{align}
\Lc_{\mathrm m}&=\frac{1}{M+1}\sum_{k=0}^{M}\frac{1}{N_k}\sum_{i=1}^{N_k}
\bigg(\frac{\Rm_Q(K_i^k)}{|K_i^k|\,S_{\mathrm m}}\bigg)^{2},
\qquad
\Lc_{\mathrm p}=\frac{1}{M+1}\sum_{k=0}^{M}\frac{1}{N_k}\sum_{i=1}^{N_k}
\bigg(\frac{\Rp_Q(K_i^k)}{|K_i^k|\,S_{\mathrm p}}\bigg)^{2},
\label{eq:loss_cells}\\
\Lc_{\mathrm I}&=\big\langle(r_q/S_q)^{2}\big\rangle
+\big\langle(r_H/S_H)^{2}\big\rangle ,
\qquad
\Lc_{\mathrm B}=\sum_{\text{ends}}\big\langle r_{\mathrm B}^{2}\big\rangle ,
\label{eq:loss_ifc}
\end{align}
where $\langle\cdot\rangle$ is the arithmetic mean over the corresponding node
set and $r_{\mathrm B}$ is the boundary residual of
Section~\ref{sec:interface}, normalised by $S_W$ for the transmissive and by
$S_q$ for the wall treatment. The well-balanced property, the interface closure
and the boundary treatment are thus separate terms of one objective and can be
modified independently. An entropy residual built on the classical energy--entropy pair
\cite{Bouchut2004}
$E=\tfrac12qu+\tfrac12gh^{2}+ghb$,
$G=(\tfrac12qu+gh^{2}+ghb)u$ is implemented but disabled throughout.

\section{Theoretical analysis}\label{sec:theory}

This section establishes the properties on which the method rests.
Section~\ref{sec:consistency} identifies the source functional with the
geometric source and characterises the zeros of the continuous residuals.
Section~\ref{sec:wb} proves that the fully discrete residual annihilates the
lake at rest under hypotheses weaker than those used to construct it, and
deduces that the equilibrium is an exact global minimiser of the objective.
Section~\ref{sec:shift} quantifies the perturbation introduced by the bed-slope
correction, Section~\ref{sec:necessity} proves that the sub-domain
decomposition is necessary, and Section~\ref{sec:reduction} records the
reduction to the homogeneous formulation and the conservation properties.
Section~\ref{sec:open} discusses the limitations of the analysis and directions
for future work.

\subsection{The source pairing and consistency of the weak residual}
\label{sec:consistency}

For $L^\infty$ data the traces in \eqref{eq:res_mass}--\eqref{eq:source_ibp}
exist, as limits of averages, for almost every $K\in\mathcal K_k$; for the
neural approximation \eqref{eq:ansatz}, which is continuous on
$\overline{{\Omega_k^{T}}}$ whenever $(\eta_0,q_0)$ are continuous on
$\overline{\Omega_k}$, they are then classical and every statement below holds
for every $K\in\mathcal K_k$.

The next lemma is a direct consequence of the integration-by-parts formula for
functions of bounded variation \cite[Ch.~3]{AmbrosioFuscoPallara2000}; it
identifies the functional \eqref{eq:source_ibp} with the geometric source
rather than with a substitute for it.

\begin{lemma}[Source pairing]
\label{lem:source_pairing}
Let Assumption~\ref{ass:bed} hold, and let $\eta(\cdot,t)\in BV(\Omega_k)$
for almost every $t$, with $\sup_t\mathrm{TV}(\eta(\cdot,t))<\infty$. Let
$K\in\mathcal K_k$ be such that, for almost every $t\in[t_1,t_2]$, $x_1$ and
$x_2$ are points of approximate continuity of $\eta(\cdot,t)$. Then
\begin{equation}
\Sm(K)=g\int_{t_1}^{t_2}\!\int_{[x_1,x_2]}b(x)\,
\dd\big(D_x\eta(\cdot,t)\big)(x)\,\dd t ,
\label{eq:source_pairing}
\end{equation}
where $D_x\eta(\cdot,t)$ is the distributional derivative of $\eta(\cdot,t)$,
a finite Radon measure. In particular $\Sm(K)=\int_Kg\,b\,\partial_x\eta$
whenever the latter is classically defined.
\end{lemma}

\begin{proof}
For fixed admissible $t$, the integration-by-parts formula for a $BV$ function
against $b\in C^{1}([x_1,x_2])$ \cite{AmbrosioFuscoPallara2000} gives
$\int_{[x_1,x_2]}b\,\dd(D_x\eta(\cdot,t))
=b(x_2)\eta(x_2,t)-b(x_1)\eta(x_1,t)-\int_{x_1}^{x_2}\eta\,b'\,\dd x$,
the boundary terms being unambiguous at points of approximate continuity. The
uniform total-variation bound makes the left-hand side bounded and measurable
in $t$, and $\eta\,b'\in L^\infty(K)$ justifies Fubini's theorem for the last
term; multiplying by $g$ and integrating over $[t_1,t_2]$ gives
\eqref{eq:source_pairing}.
\end{proof}

For homogeneous conservation laws, the equivalence between the integral form
\eqref{eq:we_rect} on rectangles and the distributional form is classical
\cite{Dafermos2016,godlewski2013numerical}. The next result extends it to
\eqref{eq:swe_eta} with the source in the form \eqref{eq:source_ibp}, and
licenses the use of \eqref{eq:res_mass}--\eqref{eq:res_mom} as a training
residual. The proof follows the standard mollification argument; only the
treatment of the source terms $C$ and $D$ below is specific to the present
setting.

\begin{proposition}[Consistency]
\label{prop:consistency}
Let Assumptions~\ref{ass:bed}--\ref{ass:states} hold. The following are
equivalent:
\begin{enumerate}
\item[(i)] $\Rm(K)=\Rp(K)=0$ for almost every $K\in\mathcal K_k$;
\item[(ii)] for every $\varphi\in C_c^{1}({\Omega_k^{T}})$,
\begin{align}
&\int_{{\Omega_k^{T}}}\big(\eta\,\partial_t\varphi+q\,\partial_x\varphi\big)\,\dd x\,\dd t=0 ,
\label{eq:distrib_mass}\\
&\int_{{\Omega_k^{T}}}\Big(q\,\partial_t\varphi+\Phi\,\partial_x\varphi
-g\,\eta\,\big(b\,\partial_x\varphi+b'\,\varphi\big)\Big)\,\dd x\,\dd t=0 .
\label{eq:distrib_mom}
\end{align}
\end{enumerate}
Identities \eqref{eq:distrib_mass}--\eqref{eq:distrib_mom} state that
$(\eta,q)$ solves \eqref{eq:swe_eta} in $\mathcal D'({\Omega_k^{T}})$, the source being
understood through the Leibniz identity
$g\,b\,\partial_x\eta=\partial_x(g\,b\,\eta)-g\,\eta\,b'$, valid in
$\mathcal D'({\Omega_k^{T}})$ because $b\in C^{1}$. If moreover
$(\eta,q)\in C(\overline{{\Omega_k^{T}}})^{2}$, ``almost every'' in (i) may be replaced by
``every''.
\end{proposition}

\begin{proof}
Set $A:=q$, $B:=\Phi$, $C:=g\,b\,\eta$ and $D:=g\,\eta\,b'$, all in
$L^\infty({\Omega_k^{T}})$. By \eqref{eq:res_mom}--\eqref{eq:source_ibp},
\begin{equation}
\Rp(K)=\int_{x_1}^{x_2}\!\big[A(x,t_2)-A(x,t_1)\big]\dd x
+\int_{t_1}^{t_2}\!\big[(B-C)(x_2,t)-(B-C)(x_1,t)\big]\dd t
+\int_K D .
\label{eq:Rp_ABCD}
\end{equation}
Let $\rho_\varepsilon$ be a standard mollifier supported in the ball of radius
$\varepsilon$, write $f_\varepsilon:=f*\rho_\varepsilon$, let
$\omega\Subset {\Omega_k^{T}}$ with $\mathrm{dist}(\omega,\partial {\Omega_k^{T}})>\varepsilon$, and
denote by $\Rm_\varepsilon(K)$, $\Rp_\varepsilon(K)$ the residuals formed with
the mollified functions for $K\subset\omega$.

\emph{(i)$\Rightarrow$(ii).} Every term of \eqref{eq:res_mass} and
\eqref{eq:Rp_ABCD} is an integral over a face of $K$ or over $K$, so by
Fubini's theorem
$\Rm_\varepsilon(K)=\int\rho_\varepsilon(z)\Rm(K-z)\,\dd z$, and likewise for
$\Rp$. Translations preserve Lebesgue measure on the parameter space, so by (i)
and Fubini's theorem on the product set
$\{(K,z):\Rm(K-z)\neq0\}$, for almost every $K$ one has $\Rm(K-z)=0$ for almost
every $z$; hence $\Rm_\varepsilon(K)=\Rp_\varepsilon(K)=0$ for almost every,
and by continuity in the parameters for every, $K\subset\omega$. The mollified
functions are smooth, so the fundamental theorem of calculus turns these
identities into
$\int_K(\partial_t\eta_\varepsilon+\partial_xq_\varepsilon)=0$ and
$\int_K(\partial_tA_\varepsilon+\partial_x(B_\varepsilon-C_\varepsilon)+D_\varepsilon)=0$
for every rectangle $K\subset\omega$; continuous integrands with vanishing
integrals over all rectangles vanish identically, so
\begin{equation}
\partial_t\eta_\varepsilon+\partial_x q_\varepsilon=0,
\qquad
\partial_t A_\varepsilon+\partial_x(B_\varepsilon-C_\varepsilon)+D_\varepsilon=0
\qquad\text{on }\omega .
\label{eq:pointwise_moll}
\end{equation}
Multiplying by $\varphi\in C_c^1({\Omega_k^{T}})$ with $\mathrm{supp}\,\varphi\subset\omega$,
integrating by parts and letting $\varepsilon\to0$, using
$f_\varepsilon\to f$ in $L^1_{\mathrm{loc}}$, gives
\eqref{eq:distrib_mass}--\eqref{eq:distrib_mom}.

\emph{(ii)$\Rightarrow$(i).} Convolution commutes with constant-coefficient
derivatives, so (ii) implies \eqref{eq:pointwise_moll} on $\omega$, and
integration over rectangles gives $\Rm_\varepsilon(K)=\Rp_\varepsilon(K)=0$ for
every $K\subset\omega$. For $f\in L^\infty({\Omega_k^{T}})$, Fubini's theorem and
$f_\varepsilon\to f$ in $L^1(\omega)$ give, along a subsequence
$\varepsilon_j\to0$, convergence of the horizontal face integrals for almost
every $(t_1,t_2)$ and of the vertical ones for almost every $(x_1,x_2)$, while
$\int_KD_{\varepsilon_j}\to\int_KD$ for every $K$. Applying this to
$f\in\{\eta,q,A,B-C\}$ yields $\Rm(K)=\Rp(K)=0$ for almost every
$K\subset\omega$, and exhausting ${\Omega_k^{T}}$ by such $\omega$ gives (i). If
$(\eta,q)$, hence $\Phi$, are continuous on $\overline{{\Omega_k^{T}}}$, both residuals are
continuous in the parameters, and vanishing almost everywhere implies vanishing
everywhere.
\end{proof}

Proposition~\ref{prop:consistency} concerns the continuous residual on the
whole family $\mathcal K_k$, whereas \eqref{eq:total_loss} involves finitely
many control volumes and a quadrature rule; it is a consistency statement, and
a zero of $\Lc$ is not by itself a weak solution. Two errors separate the two
levels. The first is a sampling error, arising from the finiteness of the
family and from incomplete coverage of ${\Omega_k^{T}}$. The second is a quadrature
error, governed by the regularity of the integrand on each face. On a lateral
face the integrand $t\mapsto\Phi(x_j,t)$ is discontinuous whenever a wave of
speed $s$ crosses the cell, that is whenever
\begin{equation}
\frac{s\,\Delta t}{\Delta x}\ \gtrsim\ 1 ,
\label{eq:cross_ratio}
\end{equation}
and, as is well known \cite{DavisRabinowitz1984}, the Gauss--Legendre rule
then converges at a rate governed by the discontinuity rather than by $Q$. The
minimiser of the discrete objective is
the network whose quadrature-visible residual vanishes, which need not resolve
the discontinuity; this mechanism predicts a smearing of fast shocks that
correlates with $s_{\max}h_t/h_x$ and is absent when the lateral integrands are
Lipschitz, as for rarefactions.

\subsection{Exact preservation of the lake at rest}
\label{sec:wb}

\begin{definition}[Discrete residuals]
\label{def:disc_res}
Let $K\in\mathcal K_k$, $Q\ge1$, and let $\beta=(\beta_1,\dots,\beta_Q)\in\R^{Q}$
be a nominal bed slope at the nodes $x_p$. Set
\begin{equation}
G_Q(K;\beta):=\frac{\Delta x}{2}\sum_{p=1}^{Q}w_p\beta_p ,
\qquad
\tilde\beta_p:=\beta_p+\frac{\big(b(x_2)-b(x_1)\big)-G_Q(K;\beta)}{\Delta x} ,
\label{eq:def_shift}
\end{equation}
and define
\begin{align}
\Rm_Q(K)&:=\frac{\Delta x}{2}\sum_{p}w_p\big[\eta(x_p,t_2)-\eta(x_p,t_1)\big]
+\frac{\Delta t}{2}\sum_{p}w_p\big[q(x_2,t_p)-q(x_1,t_p)\big] ,
\label{eq:def_RmQ}\\
\Sm_Q(K)&:=\frac{g\Delta t}{2}\sum_{p}w_p\big[b(x_2)\eta(x_2,t_p)-b(x_1)\eta(x_1,t_p)\big]
-\frac{g\Delta x\Delta t}{4}\sum_{p,r}w_pw_r\,\eta(x_p,t_r)\,\tilde\beta_p ,
\label{eq:def_SQ}\\
\Rp_Q(K)&:=\frac{\Delta x}{2}\sum_{p}w_p\big[q(x_p,t_2)-q(x_p,t_1)\big]
+\frac{\Delta t}{2}\sum_{p}w_p\big[\Phi(x_2,t_p)-\Phi(x_1,t_p)\big]-\Sm_Q(K) .
\label{eq:def_RpQ}
\end{align}
If the sub-domain carrying $K$ is flagged constant-bed, the double sum in
\eqref{eq:def_SQ} is omitted.
\end{definition}

Definition~\ref{def:disc_res} is \eqref{eq:source_quad}--\eqref{eq:bp_correction}
with the nominal slope left unspecified: $\beta$ is not required to approximate
$b'$.

\begin{lemma}[Exactness of the corrected quadrature]
\label{lem:correction}
For every $K\in\mathcal K_k$, every $Q\ge1$ and every $\beta\in\R^{Q}$,
\begin{equation}
\frac{\Delta x}{2}\sum_{p=1}^{Q}w_p\,\tilde\beta_p=b(x_2)-b(x_1) .
\label{eq:correction_exact}
\end{equation}
\end{lemma}

\begin{proof}
By \eqref{eq:def_shift} and \eqref{eq:weights_sum},
$\frac{\Delta x}{2}\sum_pw_p\tilde\beta_p
=G_Q(K;\beta)+\big[(b(x_2)-b(x_1))-G_Q(K;\beta)\big]\frac{1}{2}\sum_pw_p
=b(x_2)-b(x_1)$.
\end{proof}

\begin{theorem}[Exact discrete $C$-property]
\label{thm:wb}
Let $K$ be a space--time rectangle, $Q\ge1$, and $\beta\in\R^{Q}$ arbitrary.
Let $\eta^\star\in\R$ and assume that on $\overline K$
\[
\eta\equiv\eta^\star,\qquad q\equiv0,\qquad
\eta^\star-\max_{\overline K}b\ \ge\ h_{\min} .
\]
Then, in exact arithmetic, $\Rm_Q(K)=0$ and $\Rp_Q(K)=0$. The conclusion holds
whether or not $b$ is continuous on $K$, whether or not $\beta$ approximates
$b'$, and for every control-volume family; the only quantities of the
bathymetry that enter are $b(x_1)$ and $b(x_2)$.
\end{theorem}

\begin{proof}
\emph{Mass residual.} Since $\eta$ is constant and $q\equiv0$ on $\overline K$,
both sums of \eqref{eq:def_RmQ} vanish term by term, so $\Rm_Q(K)=0$.

\emph{Momentum residual, conservative part.} The first sum of
\eqref{eq:def_RpQ} vanishes because $q\equiv0$. The hypothesis
$\eta^\star-\max_{\overline K}b\ge h_{\min}$ gives $h_\delta=h\ge h_{\min}>0$,
so $\Phi=q^{2}/h+\tfrac12g\eta^{2}=\tfrac12g(\eta^\star)^{2}$ at every node, a
value independent of $x$ and $t$, and the second sum vanishes term by term.

\emph{Source.} Substituting $\eta\equiv\eta^\star$ in \eqref{eq:def_SQ} and
using \eqref{eq:weights_sum} twice,
\[
\Sm_Q(K)
=g\,\eta^\star\,\Delta t\,\big(b(x_2)-b(x_1)\big)
-g\,\eta^\star\,\Delta t\,\frac{\Delta x}{2}\sum_pw_p\tilde\beta_p ,
\]
and the two terms coincide by Lemma~\ref{lem:correction}. If the sub-domain is
flagged constant-bed, the double sum is omitted and $b(x_1)=b(x_2)$, so the
remaining face term vanishes. Hence $\Sm_Q(K)=0$ and $\Rp_Q(K)=0$.
\end{proof}

\begin{remark}[Scope of the hypotheses]
\label{rem:wb_weak}
\emph{(i)} No property of the control-volume family enters the proof, which is
cell-local and uses only \eqref{eq:weights_sum}. The unstructured, overlapping,
non-partitioning family of Section~\ref{sec:cells} therefore carries a property
that, in the finite-volume setting, is carried by the mesh through the
reconstruction operator. \emph{(ii)} No property of $\beta$ enters: the accuracy
of the nominal bed slope is irrelevant to well-balancedness and matters only for
consistency (Section~\ref{sec:shift}); the two requirements are decoupled
because \eqref{eq:def_shift} defines the shift from the computed value
$G_Q(K;\beta)$. \emph{(iii)} Assumption~\ref{ass:bed} is not invoked. If $b$
jumps strictly inside $K$, as happens when the detector misses a jump, then
$b(x_2)-b(x_1)$ still records the total bed increment across the cell, $\beta$
carries no information about it, and the shift absorbs the difference; the
conclusion is unchanged. \emph{(iv)} The hypothesis
$\eta^\star-\max_{\overline K}b\ge h_{\min}$ is not needed for the cancellation
itself. With the regularised depth $h_\delta=\max\{\eta-b,h_{\min}\}$ used in the
implementation, $q\equiv0$ gives $q^{2}/h_\delta=0$ at every node whatever the
value of $h_\delta$, so the flux is $\tfrac12g(\eta^\star)^{2}$ at every node and the
proof goes through unchanged. The hypothesis is retained because it makes the
state physically admissible and gives $h_\delta=h$, so that the discrete residual
represents the balance law \eqref{eq:swe_eta}.
\end{remark}

Theorem~\ref{thm:wb} concerns a single control volume. Assembling it over the
family gives a statement about the objective.

\begin{lemma}[Exact representability of the equilibrium]
\label{lem:representable}
Let $\eta_0\equiv\eta^\star$ and $u_0\equiv0$. There exists $\theta^\star$ such
that $\eta_{\theta^\star}\equiv\eta^\star$ and $q_{\theta^\star}\equiv0$ on
$\overline\Omega$.
\end{lemma}

\begin{proof}
The argument is elementary. The last layer of each sub-domain network is affine,
$\mathcal N_{\theta_k}=W_L\,\sigma(\cdots)+c_L$. Choosing $W_L=0$ and
$c_L=(\eta^\star,0)^{\!\top}$ in every sub-domain gives
$\mathcal N_{\theta^\star}\equiv(\eta^\star,0)$ irrespective of the hidden
parameters, and \eqref{eq:ansatz} returns
$(1-\tau)\eta^\star+\tau\eta^\star=\eta^\star$ and $0$ for every $(x,t)$.
\end{proof}

\begin{corollary}[The equilibrium is an exact global minimiser]
\label{cor:minimiser}
Let $\eta_0\equiv\eta^\star$ and $u_0\equiv0$ with
$\eta^\star-\max_{\overline{\Omega_x}}b\ge h_{\mathrm{safe}}\ge h_{\min}$, and
let $\theta^\star$ be as in Lemma~\ref{lem:representable}. Then, for both the
transmissive and the wall boundary treatment,
\[
\Lc(\theta^\star)=0=\min_\theta\Lc(\theta) .
\]
\end{corollary}

\begin{proof}
$\Lc_{\mathrm m}$ and $\Lc_{\mathrm p}$ vanish by Theorem~\ref{thm:wb} applied
to each $K_i^k$. At every interface $q_-=q_+=0$ and
$\eta_-=\eta_+=\eta^\star$, hence $u_\mp=0$ and $r_q=r_H=0$ by
\eqref{eq:iface_res}, whatever $\tau$. The barrier is inactive since
$h=\eta^\star-b\ge h_{\mathrm{safe}}$. For the wall treatment the boundary
residual is $q/S_q=0$; for the transmissive treatment
$W_{\mathrm L}=2\sqrt{g(\eta^\star-b(x_{\mathrm L}))}=(u+2\sqrt{gh})|_{x_{\mathrm L}}$
for all $t$, and symmetrically at $x_{\mathrm R}$. Hence $\Lc(\theta^\star)=0$,
and $\Lc\ge0$ as a finite sum of squares and one-sided squares.
\end{proof}

Corollary~\ref{cor:minimiser} is the variational counterpart of the
$C$-property of a finite-volume scheme, and the analogy must be read precisely.
In a finite-volume scheme the equilibrium is a fixed point of the update map
and is preserved for all time by the algorithm itself. Here it is a global
minimiser of the objective lying in the range of the parameterisation, but the
training algorithm is not guaranteed to reach it. The deviation of a trained
network from the lake at rest is therefore an \emph{optimisation} error, not a
discretisation error: it is controlled by no discretisation parameter, and
further optimisation reduces it without the construction ever eliminating it.
Uniqueness is not claimed: the zero set of $\Lc$ contains at least the family
obtained by varying the hidden parameters in Lemma~\ref{lem:representable}.

In finite precision the statement degrades in a way that is itself a
diagnostic. The following proposition is an application of the standard
rounding-error analysis of sums and inner products \cite{Higham2002} to the
residual of Definition~\ref{def:disc_res}.

\begin{proposition}[Floating-point behaviour]
\label{prop:roundoff}
Adopt the standard model of floating-point arithmetic~\cite{Higham2002} with
unit round-off $\mathbf u$ and $\gamma_n:=n\mathbf u/(1-n\mathbf u)$. Evaluate
$\Rm_Q(K)$ and $\Rp_Q(K)$ at the state of Theorem~\ref{thm:wb}, using the same
stored values of $b(x_1)$ and $b(x_2)$ in the face term of \eqref{eq:def_SQ}
and in the shift \eqref{eq:def_shift}, and let $n_Q=O(Q^{2})$ be the number of
operations required by \eqref{eq:def_SQ}. Then
\begin{equation}
\mathrm{fl}\big(\Rm_Q(K)\big)=0,
\qquad
\big|\mathrm{fl}\big(\Rp_Q(K)\big)\big|
\ \le\ \gamma_{c\,n_Q}\,g\,|\eta^\star|\,\Delta t\,
\Big(|b(x_1)|+|b(x_2)|+\Delta x\max_{p}|\beta_p|\Big)
+O(\mathbf u^{2}) ,
\label{eq:roundoff_bound}
\end{equation}
where $c$ is a small absolute constant, and $\mathrm{fl}(\Rp_Q(K))=0$ exactly
on a constant-bed sub-domain.
\end{proposition}

\begin{proof}
Every term of \eqref{eq:def_RmQ} and of the two sums in \eqref{eq:def_RpQ} is a
difference of identical stored values: $q^2/h_\delta$ evaluates to $0$ exactly
whatever $h_\delta$, so $\Phi$ evaluates to the single stored value
$\mathrm{fl}(\tfrac12g(\eta^\star)^{2})$ at every node. Hence
$\mathrm{fl}(\Rm_Q(K))=0$ and $\mathrm{fl}(\Rp_Q(K))=-\mathrm{fl}(\Sm_Q(K))$.
Write $\Sm_Q(K)=T_1-T_2$ for the face and volume terms; in exact arithmetic
$T_1=T_2=g\eta^\star\Delta t\,(b(x_2)-b(x_1))$. Because $T_1$ and $T_2$ are
differences of nearly equal products, no relative bound on $\mathrm{fl}(T_j)$
is available; the rounding error is instead controlled by the absolute values
of the summands. Each $T_j$ is a sum of at most $n_Q$ rounded products
$x^{(j)}_i$, so that
$|\mathrm{fl}(T_j)-T_j|\le\gamma_{n_Q}\sum_i|x^{(j)}_i|$
\cite[Ch.~3]{Higham2002}. By \eqref{eq:weights_sum},
$\sum_i|x^{(1)}_i|\le g|\eta^\star|\Delta t\,(|b(x_1)|+|b(x_2)|)$ for the face
term and $\sum_i|x^{(2)}_i|\le g|\eta^\star|\Delta t\,\Delta x\max_p|\tilde\beta_p|$
for the volume term. Forming the shift \eqref{eq:def_shift} perturbs each
$\tilde\beta_p$ by at most
$\gamma_{cQ}\big(2\max_p|\beta_p|+(|b(x_1)|+|b(x_2)|)/\Delta x\big)$, and the
stored weights satisfy $\sum_p\mathrm{fl}(w_p)=2(1+\theta)$ with
$|\theta|\le Q\mathbf u$, which perturbs the identity of
Lemma~\ref{lem:correction} by at most
$|\theta|\,(|b(x_1)|+|b(x_2)|+\Delta x\max_p|\beta_p|)$. Since $T_1=T_2$ in exact
arithmetic, the final subtraction adds only a factor $1+\varepsilon$,
$|\varepsilon|\le\mathbf u$, and summing these contributions with
$\Delta x\max_p|\tilde\beta_p|\le2\Delta x\max_p|\beta_p|+|b(x_1)|+|b(x_2)|$
gives \eqref{eq:roundoff_bound}. On a constant-bed
sub-domain the double sum is not assembled and $b(x_1)$, $b(x_2)$ are the same
stored number, so the face term is a difference of identical partial sums and
$\mathrm{fl}(\Sm_Q(K))=0$.
\end{proof}

Dividing \eqref{eq:roundoff_bound} by $|K|$ shows that the residual density at
the equilibrium is bounded by
$\gamma_{c\,n_Q}\,g\,|\eta^\star|\,\big(|b(x_1)|+|b(x_2)|+\Delta x\max_p|\beta_p|\big)/\Delta x$:
it grows linearly in $|\eta^\star|$, scales with the unit round-off and grows
like $\Delta x^{-1}$ under refinement, while on a piecewise-constant bed the
residual itself vanishes exactly, by the last statement of the proposition. The
dependence on $|\eta^\star|$ alone does not discriminate, since an unbalanced
truncation error $g\eta^\star\Delta t\,[(b(x_2)-b(x_1))-G_Q]$ is also linear in
$\eta^\star$; but such an error does not depend on the arithmetic. A measured
residual that falls by $\mathbf u_{64}/\mathbf u_{32}\approx1.86\times10^{-9}$
when the precision is raised is therefore a round-off floor, not an imbalance.

\subsection{Consistency of the bed-slope correction}
\label{sec:shift}

Theorem~\ref{thm:wb} says nothing about the size of the perturbation introduced
by \eqref{eq:def_shift} away from equilibrium. The next two results show that it
is of the order of the quadrature error already committed. Both rest on the
classical Gauss--Legendre error formula \cite[Sec.~2.7]{DavisRabinowitz1984}:
for $f\in C^{2Q}([x_1,x_2])$,
\begin{equation}
\int_{x_1}^{x_2}f-\frac{\Delta x}{2}\sum_{p=1}^{Q}w_pf(x_p)
=C_Q\,\Delta x^{2Q+1}f^{(2Q)}(\zeta),
\qquad
C_Q=\frac{(Q!)^{4}}{(2Q+1)\big[(2Q)!\big]^{3}},
\label{eq:gl_error}
\end{equation}
for some $\zeta\in(x_1,x_2)$.

\begin{proposition}[Order of the shift]
\label{prop:shift_order}
Let $K\in\mathcal K_k$ and let
$\delta_K:=\big[(b(x_2)-b(x_1))-G_Q(K;\beta)\big]/\Delta x$ be the constant
shift of \eqref{eq:def_shift}.
\begin{enumerate}
\item[(i)] If $b\in C^{2Q+1}([x_1,x_2])$ and $\beta_p=b'(x_p)$, then
\begin{equation}
|\delta_K|\ \le\ C_Q\,\Delta x^{2Q}\,
\sup_{[x_1,x_2]}\big|b^{(2Q+1)}\big| ,
\label{eq:shift_bound}
\end{equation}
\item[(ii)] If $b$ is Lipschitz on $[x_1,x_2]$ with constant $L_b$ and
$|\beta_p|\le L_b$ for every $p$, then $|\delta_K|\le2L_b$.
\item[(iii)] If $b=b_{\mathrm{Lip}}+J\,\mathbf 1_{[x_\ast,x_2]}$ on $[x_1,x_2]$
for some $x_\ast\in(x_1,x_2)$ and $J\neq0$, with $b_{\mathrm{Lip}}$ Lipschitz
of constant $L_b$ and $|\beta_p|\le L_b$, then
\begin{equation}
\frac{|J|-2L_b\Delta x}{\Delta x}\ \le\ |\delta_K|\
\le\ \frac{|J|+2L_b\Delta x}{\Delta x} ,
\label{eq:shift_jump}
\end{equation}
so that $\delta_K=J/\Delta x+O(1)$.
\end{enumerate}
\end{proposition}

\begin{proof}
(i) Since $b(x_2)-b(x_1)=\int_{x_1}^{x_2}b'$, $\Delta x\,\delta_K$ is the error
of the $Q$-point Gauss--Legendre rule applied to $b'\in C^{2Q}$, and
\eqref{eq:gl_error} with $f=b'$, divided by $\Delta x$, gives
\eqref{eq:shift_bound}. (ii) $|b(x_2)-b(x_1)|\le L_b\Delta x$
and, by \eqref{eq:weights_sum}, $|G_Q(K;\beta)|\le L_b\Delta x$. (iii) Here
$b(x_2)-b(x_1)=b_{\mathrm{Lip}}(x_2)-b_{\mathrm{Lip}}(x_1)+J$, and the bounds of
(ii) on the remaining terms give
$|J|-2L_b\Delta x\le|\Delta x\,\delta_K|\le|J|+2L_b\Delta x$.
\end{proof}

\begin{proposition}[Consistency of the modified source]
\label{prop:source_consistency}
Let $b\in C^{2Q+1}(\overline K)$ with $\beta_p=b'(x_p)$, and let
$\eta\in C^{2Q+1}(\overline K)$. Then
\begin{equation}
\Sm_Q(K)=\int_K g\,b\,\partial_x\eta
+\mathcal E_{\mathrm{quad}}(K)
-g\,\delta_K\!\int_K\eta ,
\label{eq:source_error}
\end{equation}
where the error of the uncorrected quadrature satisfies
\begin{equation}
\big|\mathcal E_{\mathrm{quad}}(K)\big|
\ \le\ c_Q\,g\,|K|\,
\big(\Delta x^{2Q}\,M_x+\Delta t^{2Q}\,M_t\big) ,
\label{eq:quad_error}
\end{equation}
with $M_x,M_t$ the suprema over $\overline K$ of the $2Q$-th partial
derivatives in $x$ and $t$ of the integrands and $c_Q$ depending on $Q$ only,
and the correction term obeys
\begin{equation}
\Big|g\,\delta_K\!\int_K\eta\Big|
\ \le\ g\,C_Q\,\Delta x^{2Q}\,
\sup_{[x_1,x_2]}\big|b^{(2Q+1)}\big|\,\|\eta\|_{L^\infty(K)}\,|K| .
\label{eq:corr_error}
\end{equation}
The correction is therefore of the same order in $\Delta x$ as the quadrature
error, and the modified source is consistent with $\int_Kg\,b\,\partial_x\eta$
at that order.
\end{proposition}

\begin{proof}
Inserting $\tilde\beta_p=\beta_p+\delta_K$ in \eqref{eq:def_SQ} splits
$\Sm_Q(K)=\Sm_Q^{0}(K)-R_\delta$, where $\Sm_Q^{0}(K)$ is the uncorrected
Gauss--Legendre approximation of \eqref{eq:source_ibp} and
$R_\delta=\frac{g\Delta x\Delta t}{4}\delta_K\sum_{p,r}w_pw_r\eta(x_p,t_r)$. By
Lemma~\ref{lem:source_pairing} the exact value of \eqref{eq:source_ibp} is
$\int_Kg\,b\,\partial_x\eta$, and the classical Gauss--Legendre error bound
\eqref{eq:gl_error} \cite{DavisRabinowitz1984} applied to each face integral and, in tensor-product form, to the volume
integral gives
$\Sm_Q^{0}(K)=\int_Kg\,b\,\partial_x\eta+\mathcal E_{\mathrm{quad}}(K)$ with
\eqref{eq:quad_error}. The same bound applied to $\eta$ gives
$R_\delta=g\,\delta_K\int_K\eta$ up to a higher-order term absorbed into
$\mathcal E_{\mathrm{quad}}$, and \eqref{eq:corr_error} follows from
\eqref{eq:shift_bound}.
\end{proof}

Two qualifications apply. In the implementation, $\beta_p$ is a second-order
central difference of step $\varepsilon_b=10^{-6}|\Omega_k|$, which adds
$O(\varepsilon_b^{2}\sup|b'''|)$ to $\delta_K$, negligible beside
\eqref{eq:shift_bound} for smooth beds. For a bed obtained by piecewise-linear
interpolation of tabulated data, $b\in C^{2Q+1}$ fails and only part~(ii) of
Proposition~\ref{prop:shift_order} applies on cells containing a kink, so the
source loses its formal order there. Neither case affects Theorem~\ref{thm:wb},
which does not use $\beta$.

\subsection{Necessity of the sub-domain decomposition}
\label{sec:necessity}

The construction assigns an independent network to each sub-domain of smooth
bed, which multiplies the parameter count and requires the interface residuals
\eqref{eq:iface_res}. The following result, an elementary consequence of the
head condition, shows that this is forced. The underlying fact, that the free
surface must jump across a step carrying a non-zero discharge, is implicit in
the exact solutions of \cite{AlcrudoBenkhaldoun2001}; the contribution here is
the explicit lower bound \eqref{eq:head_defect} and its interpretation as an
obstruction for continuous parameterisations.

\begin{theorem}[Incompatibility of a continuous free surface]
\label{thm:no_continuous}
Let $x_0$ be a bed discontinuity with $\jump{b}\neq0$, and let $(\eta,q)$
possess one-sided traces at $x_0$ with $h_\pm=\eta_\pm-b_\pm>0$ and
$u_\pm=q_\pm/h_\pm$. Assume $\jump{q}=0$, write $q$ for the common value, and
assume that $\eta$ is continuous at $x_0$. Then
\begin{equation}
\big|\jump{H}\big|
=\frac{q^{2}}{2g}\,\frac{(h_-+h_+)\,\big|\jump{b}\big|}{h_-^{2}h_+^{2}}
\ \ge\ \frac{q^{2}\,\big|\jump{b}\big|}{2g\,h_{\max}^{3}} ,
\qquad h_{\max}:=\max\{h_-,h_+\} ,
\label{eq:head_defect}
\end{equation}
and consequently $\jump{H}=0$ forces $q=0$.
\end{theorem}

\begin{proof}
Continuity of $\eta$ gives $h_+-h_-=-\jump{b}\neq0$. Since $\jump{\eta}=0$ and
$u_\pm=q/h_\pm$,
\[
\jump{H}=\frac{q^{2}}{2g}\Big(\frac{1}{h_+^{2}}-\frac{1}{h_-^{2}}\Big)
=\frac{q^{2}}{2g}\,\frac{(h_-+h_+)(h_--h_+)}{h_-^{2}h_+^{2}} ,
\]
and $h_--h_+=\jump{b}$ gives the equality. The inequality follows from
$h_-+h_+\ge h_{\max}$ and $h_-^{2}h_+^{2}\le h_{\max}^{4}$. If $\jump{H}=0$,
all other factors being non-zero, $q=0$.
\end{proof}

Theorem~\ref{thm:no_continuous} applies to any parameterisation producing a
free surface continuous at $x_0$, in particular a single globally smooth
network. If $\jump{q}=0$ is enforced, the head condition cannot be satisfied,
and the objective retains a residual of size at least
$q^{2}|\jump{b}|/(2gh_{\max}^{3})$ at that abscissa. A penalty formulation then
trades the head condition against the continuity imposed by the
parameterisation, and its minimiser solves a modified problem in which
$\jump{\eta}=0$ has replaced $\jump{H}=0$: it may converge cleanly, but to the
solution of a \emph{different} Riemann problem. The bound grows quadratically
in the discharge and cubically in the inverse depth, so the obstruction is
most severe for fast, shallow flow over a large step, and it vanishes at
$q=0$, which is why the lake at rest imposes no such requirement. For instance,
for the exact states of configuration~6.3.3 of \cite{AlcrudoBenkhaldoun2001},
with $q=-3.468$~m$^{2}$\,s$^{-1}$, $\jump{b}=1$~m and $h_{\max}=2.634$~m,
\eqref{eq:head_defect} gives $|\jump{H}|\ge3.3\times10^{-2}$~m against the
exact value $\jump{H}=0$.

\subsection{Reduction, conservation and the momentum defect}
\label{sec:reduction}

\begin{proposition}[Strict extension]
\label{prop:reduction}
Let $b\equiv b_0$ be constant on $\Omega_x$. Then $M=0$, the interface terms
are absent, and for every $K\in\mathcal K_0$
\[
\Rm(K)=\int_{x_1}^{x_2}\!\big[h(x,t_2)-h(x,t_1)\big]\dd x
+\int_{t_1}^{t_2}\!\big[q(x_2,t)-q(x_1,t)\big]\dd t ,
\]
\[
\Rp(K)=\int_{x_1}^{x_2}\!\big[q(x,t_2)-q(x,t_1)\big]\dd x
+\int_{t_1}^{t_2}\!\big[\mathcal F(x_2,t)-\mathcal F(x_1,t)\big]\dd t ,
\qquad
\mathcal F:=\frac{q^{2}}{h}+\tfrac12 g h^{2} ,
\]
which is the weak residual \eqref{eq:we_rect} of the homogeneous shallow water
system in the variables $(h,q)$. The same identities hold for the discrete
residuals \eqref{eq:def_RmQ}--\eqref{eq:def_RpQ}.
\end{proposition}

\begin{proof}
Since $b$ is constant, $\eta(x,t_2)-\eta(x,t_1)=h(x,t_2)-h(x,t_1)$. The volume
terms of \eqref{eq:source_ibp} and \eqref{eq:def_SQ} vanish (in the discrete
case $\beta=0$ and $b(x_1)=b(x_2)$ give $\tilde\beta=0$), and the lateral
integrand of $\Rp$ becomes
$\Phi-g\,b_0\,\eta=q^{2}/h+\tfrac12g(\eta-b_0)^{2}-\tfrac12gb_0^{2}
=\mathcal F-\tfrac12gb_0^{2}$; the constant cancels between the two lateral
faces. The discrete statement follows node by node.
\end{proof}

If $b$ is smooth but not constant, the decomposition still collapses to a
single network and the method reduces to the inherited residual augmented only
by the source treatment of Sections~\ref{sec:source}--\ref{sec:correction}.

The global mass balance below is the usual telescoping property of the
integral conservation law \cite{LeVeque2002}, written for the decomposed
domain so as to expose the interface defect.

\begin{proposition}[Global mass balance]
\label{prop:mass}
Let $K_k^{\mathrm{full}}:=\overline{\Omega_k}\times[0,T]$ and assume
$\Rm(K_k^{\mathrm{full}})=0$ for every $k$. Then
\begin{equation}
\underbrace{\int_{\Omega_x}\!\big[\eta(\cdot,T)-\eta(\cdot,0)\big]\dd x}_{\text{storage}}
+\underbrace{\int_0^T\!\big[q(x_{\mathrm R},t)-q(x_{\mathrm L},t)\big]\dd t}_{\text{net outflux}}
+\underbrace{\sum_{m=1}^{M}\int_0^T\!\jump{q}(e_m,t)\,\dd t}_{\text{interface defect}}=0 ,
\label{eq:mass_balance}
\end{equation}
with $\jump{q}(e_m,t):=q(e_m^-,t)-q(e_m^+,t)$. In particular, global mass is
exactly conserved if \eqref{eq:sst_q} holds exactly at every interface.
\end{proposition}

\begin{proof}
Summing $\Rm(K_k^{\mathrm{full}})=0$ over $k$, the horizontal faces give the
storage term, the sub-domains being disjoint with union $\Omega_x$ up to a
finite set, and the vertical faces give
$\sum_k\int_0^T[q(e_{k+1}^-,t)-q(e_k^+,t)]\dd t$; grouping the two traces at
each interior node $e_m$ yields the net outflux plus the interface defect.
\end{proof}

The maximal cells $K_k^{\mathrm{full}}$ do not belong to the sampled family and
the trained residuals do not vanish exactly, so \eqref{eq:mass_balance} is not
enforced by the objective; it serves instead as an \emph{a posteriori}
diagnostic, its three terms being computed independently after training.

Momentum, by contrast, is not conserved across a bed step, and the defect is
fixed by the closure rather than by the equations. This is classical in
open-channel hydraulics, where the jump of the momentum function across an
obstacle equals the force it exerts \cite{Chow1959,Henderson1966}, and it is
discussed for the Riemann problem at a step in
\cite{Bernetti2008,ValianiCaleffi2017}. We state it in the form used below.

\begin{proposition}[Momentum defect at a step~\cite{Chow1959,Henderson1966}]
\label{prop:momentum_defect}
Let $(\eta,q)$ satisfy \eqref{eq:sst_q}--\eqref{eq:sst_H} at $x_0$ with $q\neq0$
and $\jump{b}\neq0$, and assume $0<h_{\min}\le h_\pm\le h_{\max}<\infty$. Then
$h_+\neq h_-$ and
\begin{equation}
\mathcal D:=\jump{\,q\,u+\tfrac12 g h^{2}\,}
=\big(h_+-h_-\big)\Big[\tfrac{g}{2}\big(h_++h_-\big)-\tfrac{q^{2}}{h_+h_-}\Big] ,
\label{eq:mom_defect}
\end{equation}
which vanishes if and only if $q^{2}=\tfrac{g}{2}h_+h_-(h_++h_-)$, the
classical conjugate-depth relation of the hydraulic jump \cite{Chow1959}, a
codimension-one set. Moreover, if the nonconservative product
$g\,h\,\partial_xb$ is given a meaning by a regularisation for which the
corresponding spatial integrals converge and are dominated, then at every
Lebesgue point in time of both sides
\begin{equation}
\mathcal D=-\lim_{\varepsilon\to0}\int_{x_0-\varepsilon}^{x_0+\varepsilon}
g\,h\,\partial_x b ,
\label{eq:mom_defect_force}
\end{equation}
that is, the net horizontal force per unit width exerted by the step on the
flow.
\end{proposition}

\begin{proof}
If $h_+=h_-$, \eqref{eq:sst_q}--\eqref{eq:sst_H} would give $b_+=b_-$; hence
$h_+\neq h_-$. With $u_\pm=q/h_\pm$,
$\mathcal D=q^{2}(1/h_+-1/h_-)+\tfrac12g(h_+^{2}-h_-^{2})$, and factoring out
$h_+-h_-$ gives \eqref{eq:mom_defect}. For \eqref{eq:mom_defect_force},
integrate the momentum equation of \eqref{eq:swe_conservative} over
$(x_0-\varepsilon,x_0+\varepsilon)\times(t,t+\sigma)$ and let $\varepsilon\to0$:
the storage term is bounded by $4\varepsilon\|q\|_{L^\infty}$, and the flux
term converges to $\int_t^{t+\sigma}\mathcal D\,\dd s$ by dominated
convergence, the flux being bounded by
$\|q\|_{L^\infty}^{2}/h_{\min}+\tfrac12gh_{\max}^{2}$. Dividing by $\sigma$ and
letting $\sigma\to0$ at a Lebesgue point gives \eqref{eq:mom_defect_force}.
\end{proof}

Proposition~\ref{prop:momentum_defect} localises the modelling content of the
formulation: the closure enters the objective only through $r_H$ in
\eqref{eq:iface_res}, and replacing \eqref{eq:sst_H} by a momentum-based
condition as in \cite{Bernetti2008} substitutes one term of the objective. Any
closure satisfied by $q\equiv0$, $\eta\equiv\eta^\star$ leaves
Theorem~\ref{thm:wb} and Corollary~\ref{cor:minimiser} unchanged.

\subsection{Limitations and future work}
\label{sec:open}

The analysis has several limitations, which we state explicitly because each
defines a direction for further work. \emph{No estimate links the trained objective to the error.}
Converting a small value of $\Lc$ into a bound on
$\|(\eta_\theta,q_\theta)-(\eta,q)\|_{L^1}$ would require a quantitative
sampling-consistency estimate and a stability estimate for the limit problem;
the conditional $L^{1}$ estimate available for the inherited framework rests on
the Bouchut--Perthame theory~\cite{bouchut1998} for scalar conservation laws,
which does not extend to a $2\times2$ system with a geometric source and
internal stationary interfaces. \emph{Uniqueness of the zero set of $\Lc$ is not established}: a negative
answer would upgrade Corollary~\ref{cor:minimiser} from a representability
statement to a characterisation. \emph{The closure itself carries
non-uniqueness}, and no mechanism in the objective selects among admissible
branches where the Riemann problem at the step is non-unique
\cite{IsaacsonTemple1995,GoatinLeFloch2004,LeFlochThanh2007,Andrianov2005}.
Finally, dry states are excluded by
Assumption~\ref{ass:states}: the barrier \eqref{eq:barrier} discourages but does
not prevent $h\to0$, and where the floor $h_{\min}$ becomes active the residual
no longer represents \eqref{eq:swe_eta}.

These limitations point to the following directions. The first is a stability framework
for the shallow water system with discontinuous topography, which, combined with
a sampling-consistency estimate for the random control-volume family, would turn
a small value of $\Lc$ into an error bound for the trained network. Finally, the
treatment of wet--dry fronts and the extension to two space dimensions, where the
control volumes become space--time prisms and the bed discontinuities curves, are
needed for applications.

\section{Implementation}\label{sec:implementation}

Each sub-domain carries an independent fully connected network with two inputs
$(x,t)$, affinely rescaled to $[-1,1]^2$ on $\overline{{\Omega_k^{T}}}$, four hidden
layers of $n_{\mathrm h}$ units and two outputs. The activation is
$\tanh(a\,z)$ with a trainable slope $a$ shared by all layers; its initial
value $a_0$ and the width $n_{\mathrm h}$ are the only architectural quantities
that vary between test cases. The output layer is initialised with its weight
matrix scaled by $0.1$ and its bias set to $(\overline{\eta_0}|_{\Omega_k},0)$,
the mean of $\eta_0$ over the sub-domain, so that the initial network represents
a nearly flat free surface at the correct level and a nearly vanishing
discharge, from which the equilibrium of Corollary~\ref{cor:minimiser} is
reachable without a large excursion.

Since the control-volume family is fixed before training, every quantity that
depends only on the geometry and the initial conditions is evaluated once in double
precision and stored as a constant tensor: the Gauss nodes on the four faces of
each cell and, on variable-bed sub-domains, in its interior; the bed values at
those nodes, read one-sidedly as in Section~\ref{sec:detection}; the nominal
slopes $\beta_p$, obtained by central differences of step $10^{-6}|\Omega_k|$,
and the corrected slopes \eqref{eq:bp_correction}; and the values of
$\eta_0,q_0$ at the same nodes. The training loop therefore never evaluates
$b$, $\eta_0$ or $u_0$, and each evaluation of the objective reduces to a
forward and a backward pass over a fixed point set. Its dominant cost is
\begin{equation}
N_{\mathrm{eval}}=\sum_{k=0}^{M}\Big[4\,N_kQ
+\mathbf 1_{\{b\ \mathrm{variable\ on}\ \Omega_k\}}\,N_kQ^{2}\Big]
\label{eq:cost}
\end{equation}
network evaluations, to which the interface and boundary nodes add
$2MN_\Gamma Q+2(4N_\Gamma+1)$ further evaluations. A variable-bed sub-domain
contributes $Q^{2}$ interior nodes per cell against $4Q$ face nodes, so
activating the volume source multiplies its cost by $1+Q/4$, a factor of four
at $Q=12$. The constant-bed flag of Section~\ref{sec:detection} is obtained by
sampling $b$ at $513$ interior points of each $\Omega_k$.

Training follows the two-stage schedule of \cite{Oubarka2026WEPINNs}. The first
stage performs $N_{\mathrm{Adam}}$ iterations of full-batch
Adam~\cite{kingma2014adam} with learning rate $2\times10^{-3}$, cosine-annealed
to $2\times10^{-5}$. The second performs at most $N_{\mathrm{rounds}}$
restarts of L-BFGS~\cite{liu1989limited,nocedal1980updating} with strong Wolfe
line search, history size $60$, gradient tolerance $10^{-9}$ and step tolerance
$10^{-12}$, each restart limited to $N_{\mathrm{LBFGS}}$ iterations; the
sequence stops as soon as a round fails to reduce the objective by a relative
$10^{-3}$. Restarting discards stale curvature pairs accumulated across the
sharp changes of the landscape that occur when the interface residuals become
active. Network evaluation and training are carried out in single precision,
whose only consequence for the well-balanced property is the round-off floor of
Proposition~\ref{prop:roundoff}. The procedure is summarised in
Algorithm~\ref{alg:wb_wepinn}.

After training, the solution at $(x,t)$ is evaluated by the network owning $x$,
sub-domains being taken left-closed and right-open, and the bed value combined
with the network output to form $h_\theta=\eta_\theta-b$ is read from the same
side, at an abscissa pulled $10^{-9}|\Omega_k|$ inside $\Omega_k$. Combining
$\eta_\theta$ from the network owning a jump abscissa with $b$ evaluated at the
raw abscissa would mix the two sides of the step and produce a spurious depth
of the size of the bed jump at that point, at every time including $t=0$. The
same convention must be applied when a reference solution is sampled for error
norms.

The parameters that are fixed independently of the test case are collected in
Table~\ref{tab:global_params}. The case-dependent inputs are the total number
of control volumes $N$, the number of interface and boundary time intervals
$N_\Gamma$, the quadrature order $Q$, the width $n_{\mathrm h}$, the initial
slope $a_0$, the boundary treatment and the optimiser budgets
$N_{\mathrm{Adam}}$, $N_{\mathrm{rounds}}$, $N_{\mathrm{LBFGS}}$; the blending
constant is tied to the interface resolution by $c_\tau=2N_\Gamma$. All
computations use PyTorch~2.4.1 with CUDA~12.1 on a single NVIDIA
A100-SXM4-80GB accelerator.

\begin{table}[H]
\centering
\caption{Parameters common to every test case.}
\label{tab:global_params}
\begin{tabular}{llc}
\toprule
Symbol & Meaning & Value \\
\midrule
$\gamma$                  & control-volume overlap factor, \eqref{eq:cell_sizes} & $2.5$ \\
$\alpha_{\mathrm c}$      & cell aspect ratio, \eqref{eq:cell_sizes} & $1.0$ \\
--                        & hidden layers per sub-domain network & $4$ \\
$\lambda_{\mathrm I},\lambda_{\mathrm B},\lambda_{\mathrm P}$ & interface, boundary, positivity weights, \eqref{eq:total_loss} & $10,\,1,\,1$ \\
$h_{\mathrm{safe}}$       & positivity barrier threshold, \eqref{eq:barrier} & $5\times10^{-2}$~m \\
$h_{\min}$                & depth floor in $q^{2}/h$ and $q/h$ & $10^{-3}$~m \\
$c_\tau$                  & blending constant, \eqref{eq:ansatz} & $2N_\Gamma$ \\
--                        & detector tolerance and probes, \eqref{eq:tol}--\eqref{eq:detect_test} & $d_1=10^{-7}L$, $d_2=10^{-9}L$ \\
--                        & one-sided offset of bed and initial values at sub-domain ends & $10^{-9}|\Omega_k|$ \\
--                        & Adam learning rate (cosine annealed) & $2\times10^{-3}\to2\times10^{-5}$ \\
--                        & L-BFGS history / line search / restart tolerance & $60$ / strong Wolfe / $10^{-3}$ \\
$g$                       & gravitational acceleration & $9.81$~m\,s$^{-2}$ \\
--                        & working precision, random seed & \texttt{float32}, $0$ \\
\bottomrule
\end{tabular}
\end{table}

\begin{algorithm}[H]
\caption{Well-balanced WE-PINN for the shallow water equations}
\label{alg:wb_wepinn}
\begin{algorithmic}[1]
\State \textbf{Input:} $\Omega_x$, $T$, bed $b$, initial conditions $\eta_0$, $u_0$; $N$, $N_\Gamma$, $Q$, $\gamma$, $\alpha_{\mathrm c}$, $c_\tau$, loss weights
\State Locate the bed discontinuities from $b$ alone (Sec.~\ref{sec:detection}) and build the sub-domains $\Omega_0,\dots,\Omega_M$
\State Compute the scales \eqref{eq:scales} from $\eta_0$, $u_0$
\For{$k=0$ to $M$}
  \State Allot $N_k$ cells by \eqref{eq:cells_share}; build the cell family on $\Omega_k$ from \eqref{eq:cell_centres}--\eqref{eq:cell_sizes}
  \State Pre-evaluate $b$, $\eta_0$, $q_0$ at all face nodes, one-sidedly inside $\Omega_k$
  \If{$b$ is non-constant on $\Omega_k$}
     \State Pre-evaluate the interior nodes and assemble the corrected slope $\tilde b'$ from \eqref{eq:bp_correction}
  \EndIf
  \State Initialise $\mathcal N_{\theta_k}$ with output bias $(\overline{\eta_0}|_{\Omega_k},0)$
\EndFor
\State Build the interface nodes at every $e_m$ and the boundary nodes at $x_{\mathrm L}$, $x_{\mathrm R}$
\For{$\mathrm{it}=1$ to $N_{\mathrm{Adam}}$}
  \State Evaluate $\Rm_Q,\Rp_Q$ on every cell, $r_q,r_H$ at every interface node, $r_{\mathrm B}$ at every boundary node
  \State Assemble $\Lc$ from \eqref{eq:total_loss} and update $\theta$ with Adam
\EndFor
\For{$r=1$ to $N_{\mathrm{rounds}}$}
  \State Run L-BFGS for at most $N_{\mathrm{LBFGS}}$ iterations on $\Lc$; stop if the round gains less than a relative $10^{-3}$
\EndFor
\State \textbf{Output:} trained networks and the evaluator $(x,t)\mapsto(\eta_\theta,q_\theta,h_\theta,u_\theta)$
\end{algorithmic}
\end{algorithm}

\section{Numerical experiments}\label{sec:results}

\subsection{Residual-operator and trained-network diagnostics}
\label{sec:protocol}

Two logically distinct statements are tested throughout, and they are reported
separately. The first concerns the \emph{residual operator}: the exact state is
substituted into the same code that assembles the training objective, with no
training involved, which tests Theorem~\ref{thm:wb} and
Proposition~\ref{prop:roundoff} directly. The second concerns the \emph{trained
network}, whose deviation from the exact solution is, by
Corollary~\ref{cor:minimiser}, an optimisation error and not a property of the
construction. Five diagnostics are used.
\begin{enumerate}
\item[D1] \emph{Residual operator.} The exact state is fed through the residual
code; reported are $\max_i|\Rm_Q(K_i)|/|K_i|$, $\max_i|\Rp_Q(K_i)|/|K_i|$, the
largest interface and boundary residuals, the objective $\Lc(\theta^\star)$, and
the size of each of the two source contributions of \eqref{eq:def_SQ} whose
difference forms the momentum residual.
\item[D2] \emph{Errors} of the trained network in $\eta$, $h$, $u$ and $q$, in
$L^1$, $L^2$ and $L^\infty$, on a uniform grid of $8000$ points at
$t\in\{T/2,T\}$; the error at $t=0$ vanishes identically by \eqref{eq:ansatz}.
\item[D3] \emph{Space--time extrema} of the deviation from the exact solution on
a $601\times101$ grid, normalised by the physical scales $\eta^\star$,
$\sqrt{gh_{\max}}$ and $h_{\max}\sqrt{gh_{\max}}$, and the minimum depth.
\item[D4] \emph{Residuals of the trained network}: the weak cell residuals on the
training control volumes.
\item[D5] \emph{Interfaces and conservation}: the traces and the jumps
$\jump{q}$, $\jump{H}$, $\jump{\eta}$ at every detected bed jump, and the three
terms of the global mass balance \eqref{eq:mass_balance}, computed
independently after training.
\end{enumerate}
Unless stated otherwise, the parameters are those of
Table~\ref{tab:global_params}.

\subsection{Preservation of the lake at rest}
\label{sec:results_wb}

\subsubsection{Configurations}
\label{sec:wb_setup}

\begin{figure}[!htbp]
  \centering
  \includegraphics[width=\textwidth]{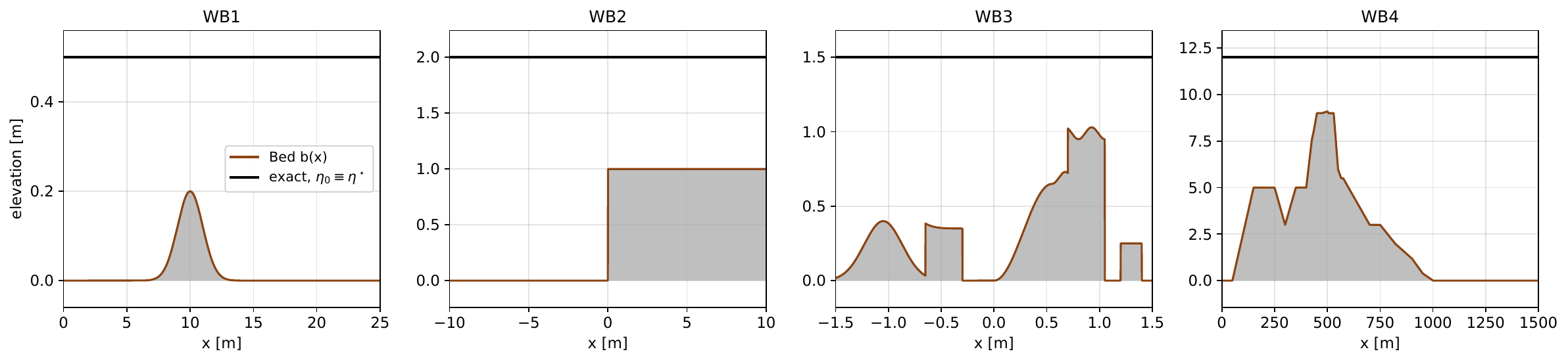}
  \caption{The four lake-at-rest bathymetries: bed $b(x)$ (grey) and exact free
  surface $\eta\equiv\eta^\star$ (black), which is also the initial condition. From
  left to right: WB1, Gaussian bump; WB2, bottom step of height $1$~m; WB3,
  compound bed with six discontinuities and two smooth features; WB4,
  piecewise-linear interpolant of $28$ tabulated points. Horizontal and vertical
  extents differ by up to three orders of magnitude between panels.}
  \label{fig:wb_beds}
\end{figure}

Four bathymetries are used (Fig.~\ref{fig:wb_beds},
Table~\ref{tab:wb_cases}). In every case the initial state is the lake at rest,
$\eta_0\equiv\eta^\star$, $u_0\equiv0$, so the exact solution on $[0,T]$ is the
initial state itself and every reported number is an absolute error. Nothing
about the solution is supplied: the bed discontinuities, the decomposition, the
constant-bed flags, the scales \eqref{eq:scales} and the boundary data are all
derived from $(b,\eta_0,u_0)$. The four cases isolate different parts of the
mechanism. WB1, a $C^\infty$ Gaussian bump, involves one network, no interface
and an active volume source, and therefore tests the bed-slope correction
alone. WB2, a single step, produces two constant-bed sub-domains on which the
volume source is not assembled, and tests the interface coupling and the
constant-bed case of Proposition~\ref{prop:roundoff}. WB3 combines smooth
features with six bed jumps of both signs, the largest ($0.95$~m) exceeding the
local depth ($0.55$~m), and contains a jump that the detector misses. WB4 is the
piecewise-linear interpolant of $28$ tabulated points (\ref{app:bed}),
continuous with $23$ kinks, over a domain and a time interval two to three
orders of magnitude larger than the others. In every case
$\min h_0/h_{\mathrm{safe}}\ge6$, so the domain is fully wet and the barrier
\eqref{eq:barrier} is inactive.

Each configuration is run twice with identical geometry generation, seed and
physics. Run~A uses a small network, two hidden layers of $3$ neurons per
sub-domain: the solution is simple, and by Lemma~\ref{lem:representable} the
equilibrium requires no approximation capacity, so a small network suffices and
needs few resources. Its budget is deliberately austere ($N=900$ requested cells,
$10^{3}$ Adam iterations, no L-BFGS, $\lambda_{\mathrm B}=100$). Run~B uses the
same network as the Riemann problems of Section~\ref{sec:results_ab}, four hidden
layers of $80$ neurons per sub-domain; it doubles the resolution in each direction
and adds up to four L-BFGS restarts ($N=3600$, $4\times10^{3}$ Adam iterations,
$\le4\times2\times10^{3}$ L-BFGS iterations, $\lambda_{\mathrm B}=1$). The residual
operator is structurally the same in both runs, whereas the network and the
optimisation budget differ, which separates the two statements of
Section~\ref{sec:protocol}. Apart from the network of Run~A, two settings depart from Table~\ref{tab:global_params}, each forced by the object
under test: the working precision is \texttt{float64}, since the quantities to be
measured are round-off floors; and both components of the state are pinned at
each end at their initial values,
\begin{equation}
r_{\mathrm B}^{\eta}=\frac{\eta(x_{\mathrm e},t)-\eta_0(x_{\mathrm e})}{S_H},
\qquad
r_{\mathrm B}^{q}=\frac{q(x_{\mathrm e},t)-q_0(x_{\mathrm e})}{S_q},
\qquad x_{\mathrm e}\in\{x_{\mathrm L},x_{\mathrm R}\},
\label{eq:bc_dirichlet}
\end{equation}
imposed as a soft penalty. This is the most constraining boundary term
available, and the proof of Corollary~\ref{cor:minimiser} applies to it
verbatim, since both residuals vanish at the equilibrium.

\begin{table}[!htbp]
\centering\footnotesize
\setlength{\tabcolsep}{4pt}
\caption{Lake-at-rest configurations and run settings. ``Jumps'' gives the
number of discontinuities detected from $b$ alone and, in parentheses, the true
number; ``sub-domains'' gives the total and, in parentheses, the number with a
variable bed, on which the volume source and the correction
\eqref{eq:bp_correction} are assembled. Networks are per sub-domain (hidden
layers $\times$ width); parameter counts are totals. For WB1, WB2 and WB4 the
cells are shared by \eqref{eq:cells_share}; for WB3, whose six sub-domains differ
widely in width, the counts were set per sub-domain as listed in
Table~\ref{tab:wb3}.}
\label{tab:wb_cases}
\begin{tabular}{lcccc}
\toprule
 & WB1 & WB2 & WB3 & WB4 \\
\midrule
$\Omega_x$ [m], $T$ [s]   & $(0,25)$, $2$ & $(-10,10)$, $1$ & $(-1.5,1.5)$, $0.15$ & $(0,1500)$, $100$ \\
$\eta^\star$ [m]; $b$ range [m] & $0.5$; $[0,0.2]$ & $2$; $\{0,1\}$ & $1.5$; $[0,1.03]$ & $12$; $[0,9.10]$ \\
$\min h_0/h_{\mathrm{safe}}$ & $6.0$ & $20.0$ & $9.4$ & $58.0$ \\
Bed regularity            & $C^\infty$ & piecewise const. & jumps $+$ smooth & Lipschitz, $23$ kinks \\
Jumps; sub-domains; interfaces & $0\,(0)$; $1\,(1)$; $0$ & $1\,(1)$; $2\,(0)$; $1$ & $5\,(6)$; $6\,(3)$; $5$ & $0\,(0)$; $1\,(1)$; $0$ \\
\midrule
Cells, Run A / B          & $900$ / $3600$ & $900$ / $3600$ & $1590$ / $4560$ & $900$ / $3600$ \\
Network, Run A / B        & $2\times3$ / $4\times80$ & $2\times3$ / $4\times80$ & $2\times3$ / $4\times80$ & $2\times3$ / $4\times80$ \\
Parameters, Run A / B     & $30$ / $19\,843$ & $59$ / $39\,685$ & $175$ / $119\,053$ & $30$ / $19\,843$ \\
Training time [s], Run A / B & $27$ / $41$ & $26$ / $61$ & $78$ / $181$ & $14$ / $36$ \\
\bottomrule
\end{tabular}
\end{table}

The pass criteria follow from the analysis. For the residual operator: (C1) the
mass, interface and boundary residuals are zero bitwise; (C2) the momentum
residual sits at the round-off floor and grows with $\eta^\star$; (C3) it is
zero bitwise on a piecewise-constant bed; (C4) it falls by
$\mathbf u_{64}/\mathbf u_{32}=2^{-53}/2^{-24}\approx1.86\times10^{-9}$ when
passing from \texttt{float32} to \texttt{float64}. Criteria C1--C4 are the
content of Theorem~\ref{thm:wb} and Proposition~\ref{prop:roundoff}. For the
trained network: (C5) the deviations are orders of magnitude below the physical
scales; (C6) they do not grow over $[0,T]$. Criteria C5--C6 measure the
optimiser and are, by Corollary~\ref{cor:minimiser}, not claims about the
construction.

\FloatBarrier
\subsubsection{The residual operator at the exact equilibrium}
\label{sec:wb_operator}

\begin{figure}[!htbp]
  \centering
  \begin{subfigure}{\textwidth}
    \centering\includegraphics[width=\textwidth]{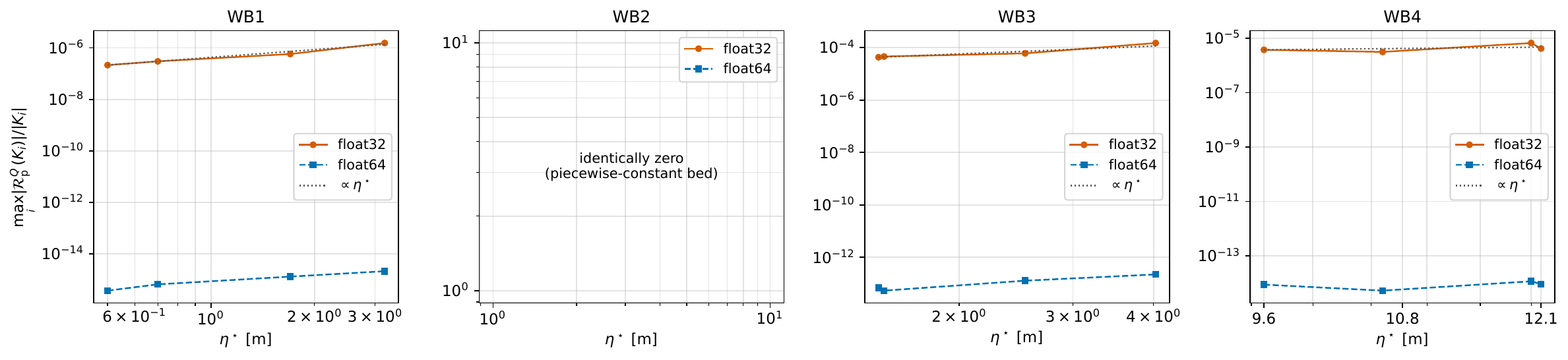}
    \caption{Run~A.}
  \end{subfigure}\\[6pt]
  \begin{subfigure}{\textwidth}
    \centering\includegraphics[width=\textwidth]{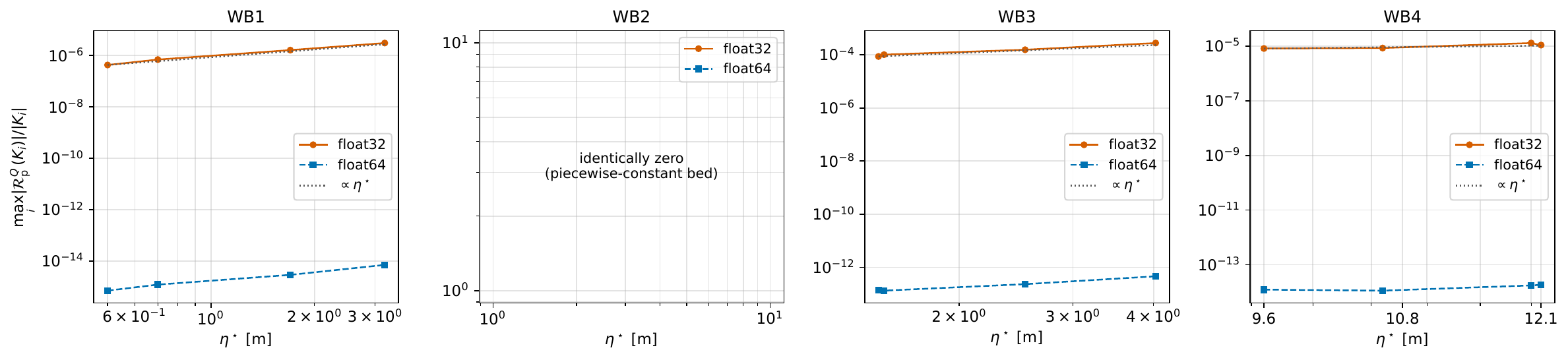}
    \caption{Run~B.}
  \end{subfigure}
  \caption{The residual-operator floor $\max_i|\Rp_Q(K_i)|/|K_i|$ at the exact
  lake at rest, as a function of the surface level $\eta^\star$, in
  \texttt{float32} (circles) and \texttt{float64} (squares), on logarithmic
  axes; the dotted line is the reference slope $\propto\eta^\star$. No training
  is involved. The two precisions are separated by about nine decades in every
  variable-bed configuration, as Proposition~\ref{prop:roundoff} requires of a
  round-off floor; the WB2 panel is empty because the residual is zero bitwise
  in both precisions on a piecewise-constant bed.}
  \label{fig:wb_precision}
\end{figure}

Diagnostic D1 is evaluated at four surface levels per configuration: the level
of the configuration itself and $\max b+0.5$, $\max b+1.5$ and
$\max b+3.0$~m, each with its own boundary datum. Every configuration is built
in both precisions with identical geometry. Table~\ref{tab:d1} summarises the
results over the four levels, and Fig.~\ref{fig:wb_precision} shows the
momentum floor level by level.

\begin{table}[!htbp]
\centering\footnotesize
\setlength{\tabcolsep}{3.5pt}
\caption{Diagnostic D1: the exact lake at rest fed through the residual code,
with no training. Each entry is the interval $[\min,\max]$ over the four surface levels. In all
$32$ evaluations the mass residual $\max_i|\Rm_Q|/|K_i|$, the interface
residuals $\max(|r_q|,|r_H|)$ and the boundary residual are exactly zero, in
both precisions; they are not tabulated. ``Source'' is the size
$\max_i|\Sm_Q^{\mathrm{face}}|/|K_i|$ of each of the two source contributions whose
difference is the momentum residual. ``Ratio'' is \texttt{f64}/\texttt{f32},
to be compared with $\mathbf u_{64}/\mathbf u_{32}\approx1.86\times10^{-9}$.}
\label{tab:d1}
\begin{tabular}{llccccc}
\toprule
Run & Case & $\max_i\frac{|\Rp_Q|}{|K_i|}$ (\texttt{f64}) &
$\max_i\frac{|\Rp_Q|}{|K_i|}$ (\texttt{f32}) & Ratio [$\times10^{-9}$] &
$\Lc(\theta^\star)$ (\texttt{f64}) & Source \\
\midrule
\multirow{4}{*}{A}
 & WB1 & $[3.6\times10^{-16},\,2.1\times10^{-15}]$ & $[2.2\times10^{-7},\,1.6\times10^{-6}]$ & $[1.35,\,2.21]$ & $[2.4\times10^{-31},\,1.3\times10^{-29}]$ & $[0.58,\,3.7]$ \\
 & WB2 & $0$ & $0$ & exact & $0$ & $0$ \\
 & WB3 & $[5.4\times10^{-14},\,2.2\times10^{-13}]$ & $[4.3\times10^{-5},\,1.5\times10^{-4}]$ & $[1.16,\,2.12]$ & $[1.6\times10^{-30},\,1.0\times10^{-29}]$ & $[78,\,210]$ \\
 & WB4 & $[5.2\times10^{-15},\,1.2\times10^{-14}]$ & $[3.2\times10^{-6},\,6.7\times10^{-6}]$ & $[1.66,\,2.33]$ & $[1.2\times10^{-30},\,2.6\times10^{-30}]$ & $[9.8,\,12]$ \\
\midrule
\multirow{4}{*}{B}
 & WB1 & $[7.0\times10^{-16},\,7.0\times10^{-15}]$ & $[4.3\times10^{-7},\,3.1\times10^{-6}]$ & $[1.63,\,2.29]$ & $[7.5\times10^{-31},\,4.8\times10^{-29}]$ & $[0.59,\,3.8]$ \\
 & WB2 & $0$ & $0$ & exact & $0$ & $0$ \\
 & WB3 & $[1.3\times10^{-13},\,4.6\times10^{-13}]$ & $[8.7\times10^{-5},\,2.8\times10^{-4}]$ & $[1.31,\,1.67]$ & $[5.1\times10^{-30},\,3.9\times10^{-29}]$ & $[154,\,413]$ \\
 & WB4 & $[1.1\times10^{-14},\,1.8\times10^{-14}]$ & $[8.2\times10^{-6},\,1.3\times10^{-5}]$ & $[1.32,\,1.69]$ & $[3.9\times10^{-30},\,8.5\times10^{-30}]$ & $[14,\,18]$ \\
\bottomrule
\end{tabular}
\end{table}

Criteria C1 and C3 are met exactly. The mass, interface and boundary residuals
are exactly zero, bitwise, on every
control volume, at every level, for all four bathymetries, both geometries and
both precisions. This is what Proposition~\ref{prop:roundoff} predicts for
$\Rm_Q$, every term of which is a difference of identical stored values, and what
Corollary~\ref{cor:minimiser} predicts for $r_q$, $r_H$ and the boundary
residuals. For WB2, the momentum residual and the objective are also zero
bitwise, at every level and in single as well as double precision: on a
constant-bed sub-domain the volume source is not assembled and the face term is
a difference of identical partial sums, as stated in the last part of
Proposition~\ref{prop:roundoff}.

Criterion C2 concerns the configurations in which the bed varies. There the
momentum residual is non-zero but confined to
values between $10^{-16}$ and $10^{-13}$, while each of the two source
contributions whose difference produces it is of size between $0.6$ and
$4\times10^{2}$. The cancellation made
exact by \eqref{eq:correction_property} is therefore observed over about
fifteen orders of magnitude. The floor increases with $\eta^\star$ in every
variable-bed configuration (Fig.~\ref{fig:wb_precision}), consistently with the
bound \eqref{eq:roundoff_bound}; the precision test C4 below is what separates
round-off from an imbalance. Over the
probe range the floor grows by factors of $5.7$ and $10.1$ for WB1 (level ratio
$6.4$) and $3.3$ for WB3 (level ratio $2.7$); exact proportionality is not
expected, \eqref{eq:roundoff_bound} being an upper bound attained on different
cells at different levels. For WB4 the probe levels span a ratio of only
$1.26$, and the level test is uninformative there. Halving the cell width
$\Delta x=2h_x$ from Run~A to Run~B, by quadrupling the number of randomly placed
cells (Section~\ref{sec:cells}; no grid is refined), increases the floor by factors of $1.9$ (WB1), $2.0$ (WB3) and
$1.5$ (WB4), against the factor $2$ implied by the per-volume form of
\eqref{eq:roundoff_bound}. Finally, $\Lc(\theta^\star)$ lies between $10^{-31}$
and $5\times10^{-29}$, the square of the floor, and is zero for WB2:
Corollary~\ref{cor:minimiser} is confirmed numerically.

The decisive test, criterion C4, is the change of precision: a round-off floor
must fall by
$\mathbf u_{64}/\mathbf u_{32}\approx1.86\times10^{-9}$, whereas a discretisation
imbalance would not move. The measured ratios lie between $1.16\times10^{-9}$ and
$2.33\times10^{-9}$ over all $24$ variable-bed evaluations, that is between $0.62$ and $1.25$
times the prediction, while the \texttt{float32} floors themselves span three
decades. Together with C1--C3 this identifies everything that remains in
Table~\ref{tab:d1} as floating-point arithmetic.

The compound bed WB3 provides a further test, of a situation that the analysis
covers but the construction was not designed for. The detector accepts five of
the six jumps of WB3, recovering their heights to
six decimal places, and rejects every smooth feature. It misses the jump of
$+0.30$~m at $x=0.70$, which lies immediately to the right of the steepest part
of the oscillation on the ramp: the bracketing stage converges to a smooth
point to its left, which the confirmation test correctly rejects. The missed
jump is retained in the test, because it realises the situation of
Remark~\ref{rem:wb_weak}(iii): some cells of the third sub-domain contain a bed
discontinuity, the nominal slope is meaningless there, and the shift must
absorb the whole increment. Table~\ref{tab:wb3} shows that it does. On the two
sub-domains without an undetected jump, $\max|\delta_K|\le2.4\times10^{-8}$,
of the order of the quadrature error, as Proposition~\ref{prop:shift_order}(i)
describes. On the sub-domain containing $x=0.70$, $|\delta_K|>1$ on $17$ of $390$
(Run~A) and $33$ of $1620$ (Run~B) cells, with maxima $5.27$ and $10.7$. These
reproduce Proposition~\ref{prop:shift_order}(iii), $\delta_K=J/\Delta x+O(1)$:
with $J=0.30$~m and $\Delta x=0.056967$ and $0.027951$~m, $J/\Delta x=5.266$ and
$10.733$, against the measured $5.266$ and $10.73$. The momentum residual on that
sub-domain is nevertheless at the same floor as on the others, so the
$C$-property holds on cells that violate the hypothesis $b\in C^{1}(K)$ under
which the residual was constructed. On the constant-bed sub-domains the residual
is $0$ on sub-domain~$4$, about $8\times10^{-28}$ on sub-domain~$3$ and
$2\times10^{-38}$ on sub-domain~$5$: the two
non-zero values are not round-off, but the face term
$g\eta^\star\Delta t\,(b(x_2)-b(x_1))$ produced by the exponentially small tail of
the Gaussian bump ($1.1\times10^{-30}$~m at $x=1.05$ and $2.4\times10^{-41}$~m
at $x=1.40$). The missed jump has no consequence for this test, because the flow
is at rest and, by Theorem~\ref{thm:no_continuous}, a continuous free surface is
compatible with the step transition when $q=0$; for a moving flow over the same
bed it would be an error.

\begin{table}[!htbp]
\centering\footnotesize
\setlength{\tabcolsep}{4pt}
\caption{WB3 at $\eta^\star=1.5$~m, by sub-domain: momentum residual of the exact
equilibrium and statistics of the bed-slope shift $\delta_K$. Sub-domains $0$--$2$
carry a variable bed; sub-domain $2$ contains the undetected jump at $x=0.70$.
\naSubNote{}Not applicable: sub-domains $3$--$5$ are flagged constant-bed
(Section~\ref{sec:detection}), so the volume source and the correction
\eqref{eq:bp_correction} are not assembled there (Section~\ref{sec:source}) and no
shift $\delta_K$ is computed; the momentum residual is still reported. Together
these three sub-domains carry $720$ (Run~A) and $1440$ (Run~B) cells.}
\label{tab:wb3}
\begin{tabular}{llcccccc}
\toprule
& & \multicolumn{3}{c}{Run~A} & \multicolumn{3}{c}{Run~B} \\
\cmidrule(lr){3-5}\cmidrule(lr){6-8}
Sub-domain & Extent [m] & $\max\frac{|\Rp_Q|}{|K|}$ & $\max|\delta_K|$ &
$\#\{|\delta_K|>1\}$ & $\max\frac{|\Rp_Q|}{|K|}$ & $\max|\delta_K|$ &
$\#\{|\delta_K|>1\}$ \\
\midrule
$0$ & $[-1.50,-0.65]$ & $2.6\times10^{-14}$ & $1.1\times10^{-8}$ & $0/240$ & $5.8\times10^{-14}$ & $2.4\times10^{-8}$ & $0/1020$ \\
$1$ & $[-0.65,-0.30]$ & $5.9\times10^{-14}$ & $7.9\times10^{-9}$ & $0/240$ & $8.5\times10^{-14}$ & $1.1\times10^{-8}$ & $0/480$ \\
$2$ & $[-0.30,+1.05]$ & $6.8\times10^{-14}$ & $5.27$ & $17/390$ & $1.4\times10^{-13}$ & $10.7$ & $33/1620$ \\
$3$ & $[+1.05,+1.20]$ & $7.4\times10^{-28}$ & \naSub & \naSub & $8.2\times10^{-28}$ & \naSub & \naSub \\
$4$ & $[+1.20,+1.40]$ & $0$ & \naSub & \naSub & $0$ & \naSub & \naSub \\
$5$ & $[+1.40,+1.50]$ & $1.9\times10^{-38}$ & \naSub & \naSub & $2.0\times10^{-38}$ & \naSub & \naSub \\
\bottomrule
\end{tabular}
\end{table}

\FloatBarrier
\subsubsection{Trained-network diagnostics across the suite}
\label{sec:wb_trained}

Because the exact solution is the initial state, every quantity reported for the
trained network is an absolute error. Table~\ref{tab:wb_trained} gives the
space--time extrema of the deviation and the residuals reached (D3, D4).

\begin{table}[!htbp]
\centering\footnotesize
\setlength{\tabcolsep}{3.5pt}
\caption{Diagnostics D3 and D4. Columns 3--5: space--time extrema of the
deviation from the exact equilibrium on a $601\times101$ grid, normalised by
$\eta^\star$, $\sqrt{gh_{\max}}$ and $h_{\max}\sqrt{gh_{\max}}$. Column 6: minimum
depth attained minus its exact value. Columns 7--8: largest weak cell residuals
of the trained network, to be compared with the residual-operator floor of
Table~\ref{tab:d1}. Column 9: objective reached.}
\label{tab:wb_trained}
\begin{tabular}{llccccccc}
\toprule
Run & Case & $\frac{\max|\eta_\theta-\eta^\star|}{\eta^\star}$ &
$\frac{\max|u_\theta|}{\sqrt{gh_{\max}}}$ &
$\frac{\max|q_\theta|}{h_{\max}\sqrt{gh_{\max}}}$ &
$\min h_\theta-\min h$ [m] &
$\max_i\frac{|\Rm_Q|}{|K_i|}$ & $\max_i\frac{|\Rp_Q|}{|K_i|}$ & $\Lc(\theta)$ \\
\midrule
\multirow{4}{*}{A}
 & WB1 & $8.4\times10^{-5}$ & $9.7\times10^{-5}$ & $5.9\times10^{-5}$ & $-3.9\times10^{-5}$ & $1.1\times10^{-4}$ & $1.3\times10^{-4}$ & $9.8\times10^{-7}$ \\
 & WB2 & $4.5\times10^{-5}$ & $1.4\times10^{-4}$ & $7.0\times10^{-5}$ & $-2.6\times10^{-5}$ & $2.1\times10^{-4}$ & $2.0\times10^{-3}$ & $8.7\times10^{-8}$ \\
 & WB3 & $1.3\times10^{-4}$ & $1.4\times10^{-4}$ & $1.1\times10^{-4}$ & $-5\times10^{-6}$ & $7.4\times10^{-3}$ & $4.0\times10^{-2}$ & $1.7\times10^{-6}$ \\
 & WB4 & $4.7\times10^{-6}$ & $1.0\times10^{-5}$ & $2.9\times10^{-6}$ & $-5.4\times10^{-5}$ & $1.8\times10^{-6}$ & $1.7\times10^{-5}$ & $8.8\times10^{-10}$ \\
\midrule
\multirow{4}{*}{B}
 & WB1 & $1.1\times10^{-16}$ & $1.9\times10^{-18}$ & $1.3\times10^{-18}$ & $0$ & $7.8\times10^{-17}$ & $6.9\times10^{-16}$ & $7.5\times10^{-31}$ \\
 & WB2 & $1.1\times10^{-16}$ & $1.7\times10^{-17}$ & $8.7\times10^{-18}$ & $0$ & $1.2\times10^{-15}$ & $1.2\times10^{-14}$ & $3.8\times10^{-32}$ \\
 & WB3 & $3.0\times10^{-16}$ & $2.6\times10^{-16}$ & $2.6\times10^{-16}$ & $0$ & $2.3\times10^{-14}$ & $2.9\times10^{-13}$ & $1.2\times10^{-29}$ \\
 & WB4 & $1.5\times10^{-16}$ & $7.6\times10^{-17}$ & $2.8\times10^{-17}$ & $0$ & $2.2\times10^{-16}$ & $1.7\times10^{-14}$ & $6.3\times10^{-30}$ \\
\bottomrule
\end{tabular}
\end{table}

\begin{figure}[!htbp]
  \centering
  \begin{subfigure}{\textwidth}
    \centering\includegraphics[width=0.8\textwidth]{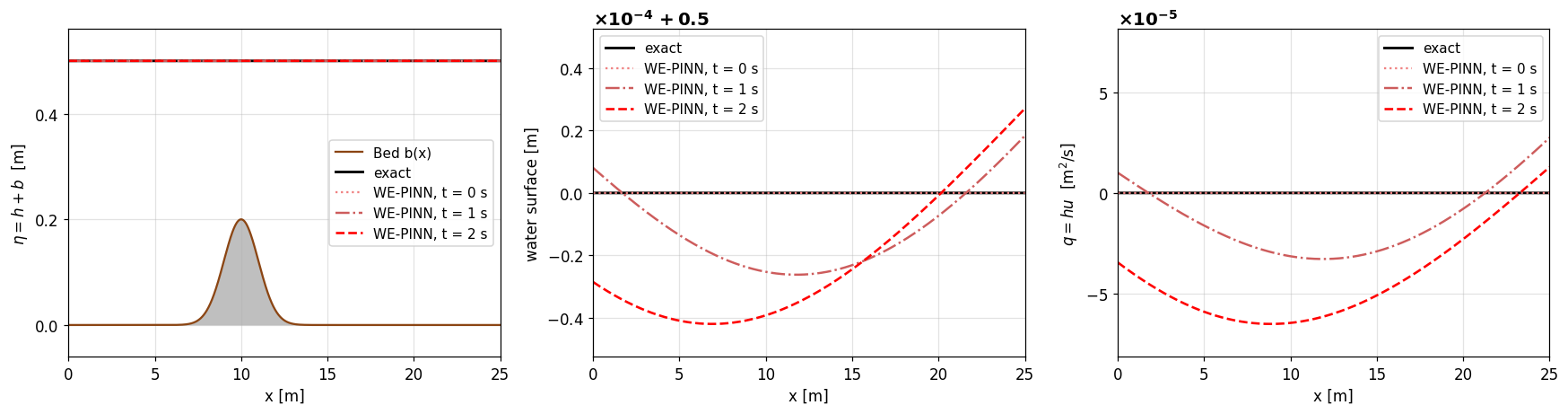}
    \caption{Run~A.}
  \end{subfigure}\\[2pt]
  \begin{subfigure}{\textwidth}
    \centering\includegraphics[width=0.8\textwidth]{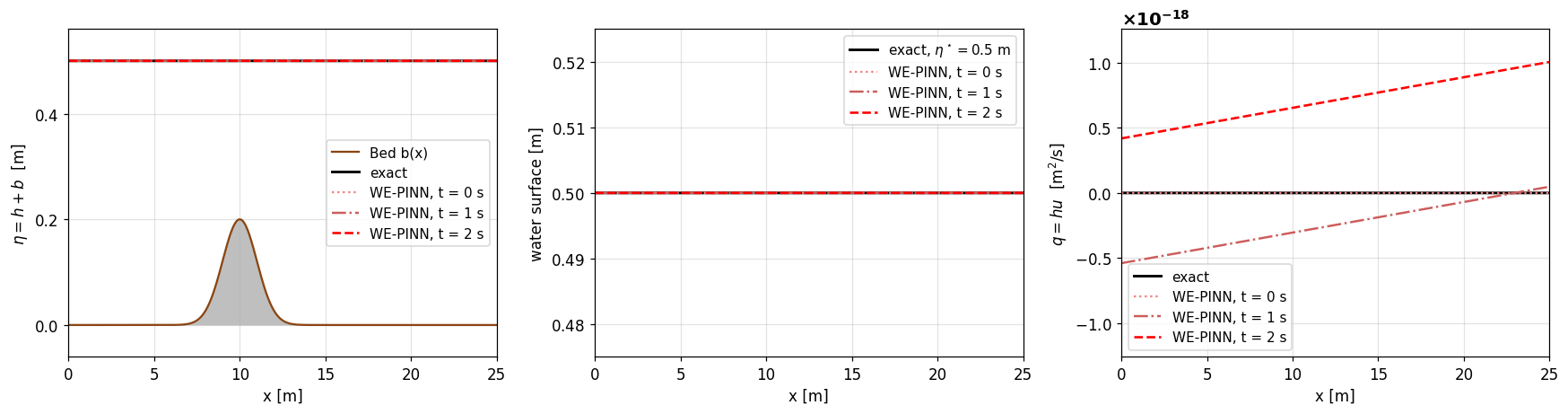}
    \caption{Run~B.}
  \end{subfigure}
  \caption{Lake at rest over a smooth Gaussian bump (WB1, $\eta^\star=0.5$~m, $T=2$~s). Left: free surface $\eta=h+b$ and bed $b$; centre: free surface magnified about the exact level $\eta^\star$; right: discharge $q$, whose exact value is zero. Black: exact equilibrium; red: WE-PINN prediction at $t=0$, $T/2$ and $T$ (three dash patterns). The factors in the corners of the centre and right panels give the magnification and offset of the vertical axis. (a) Run~A, austere budget; (b) Run~B, extended budget.}
  \label{fig:wb1}
\end{figure}

Criterion C5 is met by a wide margin: the normalised deviations never exceed $1.4\times10^{-4}$ in Run~A and
$3.0\times10^{-16}$ in Run~B. In Run~B the free surface of WB1 (Fig.~\ref{fig:wb1}) is flat to
$5.55\times10^{-17}$~m, which is $2^{-54}$, the spacing of double-precision
numbers just below $0.5$, and the minimum depth is reproduced exactly in every
case; in Run~A it falls below its exact value by at most $5.4\times10^{-5}$~m, so
the barrier \eqref{eq:barrier} is never active. Between the two runs the trained
deviation changes by ten to fourteen decades, whereas the residual-operator
floor of Table~\ref{tab:d1} changes by less than a factor of two. Comparing the
trained weak residual with that floor isolates the optimisation gap: for WB3 it
is a factor $2.1$ in Run~B ($2.9\times10^{-13}$ against $1.4\times10^{-13}$) and
a factor $5.9\times10^{11}$ in Run~A ($4.0\times10^{-2}$ against
$6.8\times10^{-14}$). The construction contributes the floor; everything above it
is the optimiser, as Corollary~\ref{cor:minimiser} states.

\begin{figure}[!htbp]
  \centering
  \begin{subfigure}{\textwidth}
    \centering\includegraphics[width=0.8\textwidth]{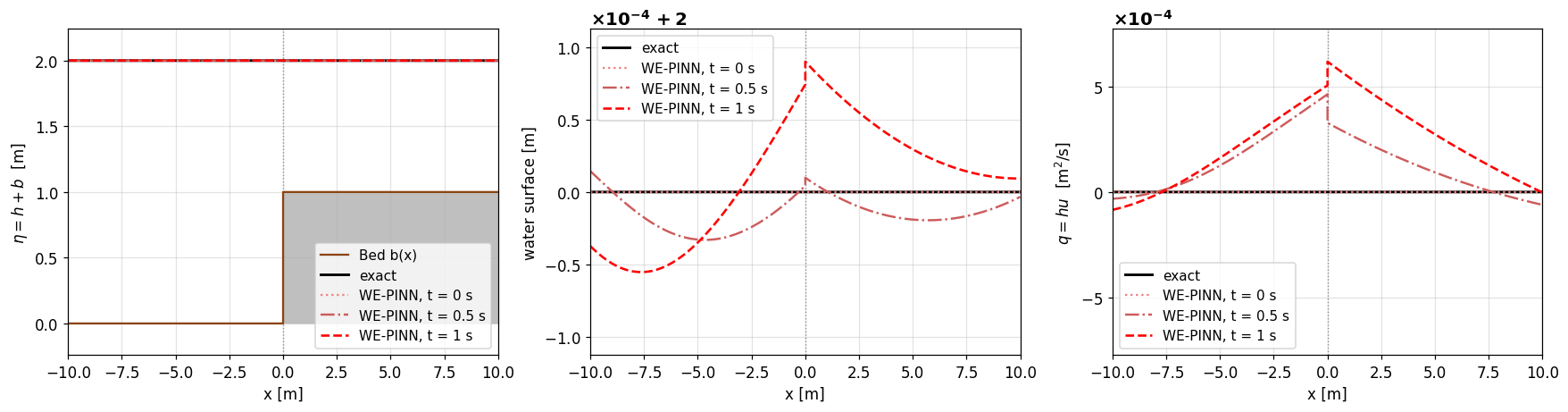}
    \caption{Run~A.}
  \end{subfigure}\\[2pt]
  \begin{subfigure}{\textwidth}
    \centering\includegraphics[width=0.8\textwidth]{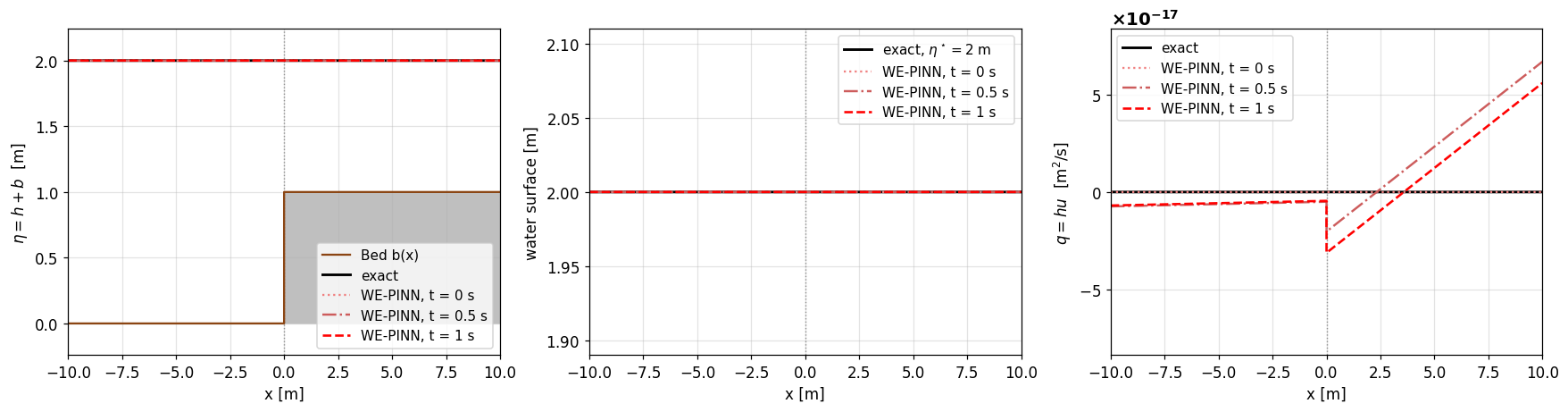}
    \caption{Run~B.}
  \end{subfigure}
  \caption{Lake at rest over a single bottom step of height $1$~m (WB2, $\eta^\star=2$~m, $T=1$~s). The dotted vertical line marks the bed jump at $x=0$, detected from $b$ alone; the depth is discontinuous there while the free surface is not. Left: free surface $\eta=h+b$ and bed $b$; centre: free surface magnified about the exact level $\eta^\star$; right: discharge $q$, whose exact value is zero. Black: exact equilibrium; red: WE-PINN prediction at $t=0$, $T/2$ and $T$ (three dash patterns). The factors in the corners of the centre and right panels give the magnification and offset of the vertical axis. (a) Run~A, austere budget; (b) Run~B, extended budget.}
  \label{fig:wb2}
\end{figure}

\begin{figure}[!htbp]
  \centering
  \begin{subfigure}{\textwidth}
    \centering\includegraphics[width=0.8\textwidth]{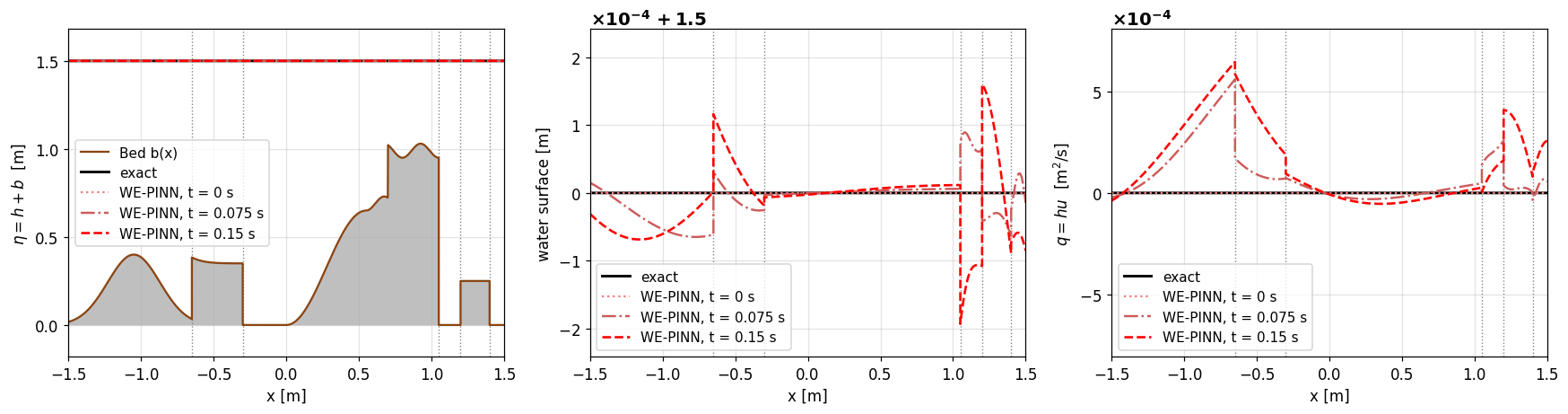}
    \caption{Run~A.}
  \end{subfigure}\\[2pt]
  \begin{subfigure}{\textwidth}
    \centering\includegraphics[width=0.8\textwidth]{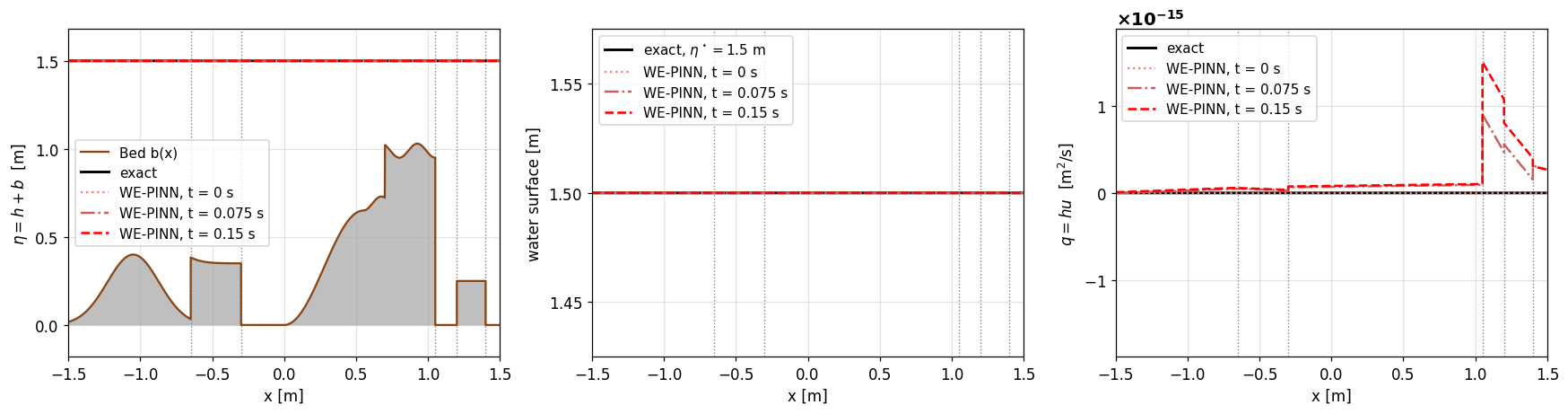}
    \caption{Run~B.}
  \end{subfigure}
  \caption{Lake at rest over the compound bathymetry (WB3, $\eta^\star=1.5$~m, $T=0.15$~s). The dotted vertical lines mark the five bed jumps detected from $b$ alone; the sixth jump, at $x=0.70$~m, is not detected and appears in the left panel as the step on the rising ramp. Left: free surface $\eta=h+b$ and bed $b$; centre: free surface magnified about the exact level $\eta^\star$; right: discharge $q$, whose exact value is zero. Black: exact equilibrium; red: WE-PINN prediction at $t=0$, $T/2$ and $T$ (three dash patterns). The factors in the corners of the centre and right panels give the magnification and offset of the vertical axis. (a) Run~A, austere budget; (b) Run~B, extended budget.}
  \label{fig:wb3}
\end{figure}

\begin{figure}[!htbp]
  \centering
  \begin{subfigure}{\textwidth}
    \centering\includegraphics[width=0.8\textwidth]{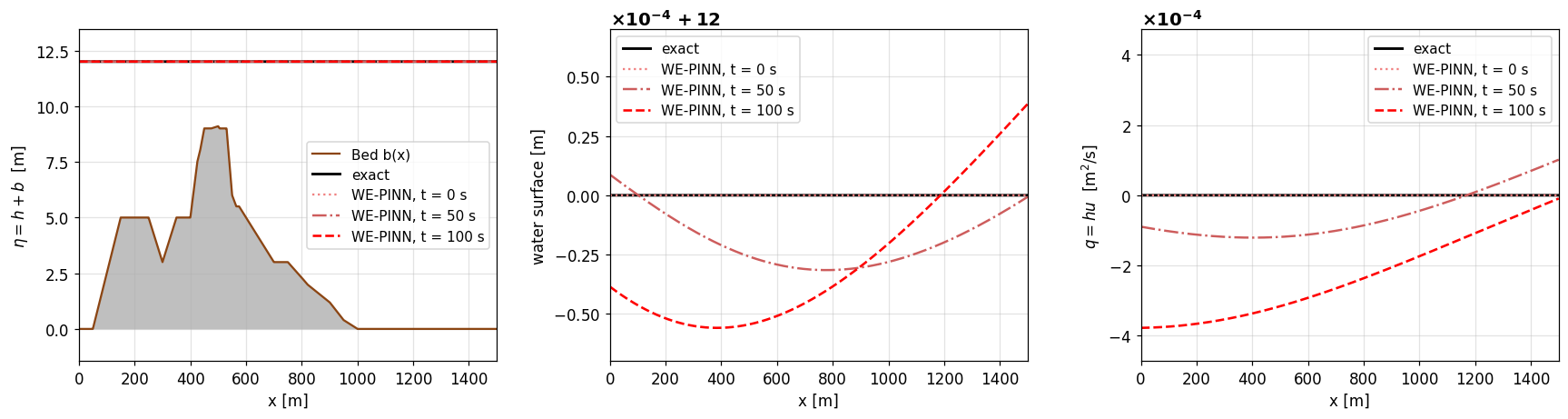}
    \caption{Run~A.}
  \end{subfigure}\\[2pt]
  \begin{subfigure}{\textwidth}
    \centering\includegraphics[width=0.8\textwidth]{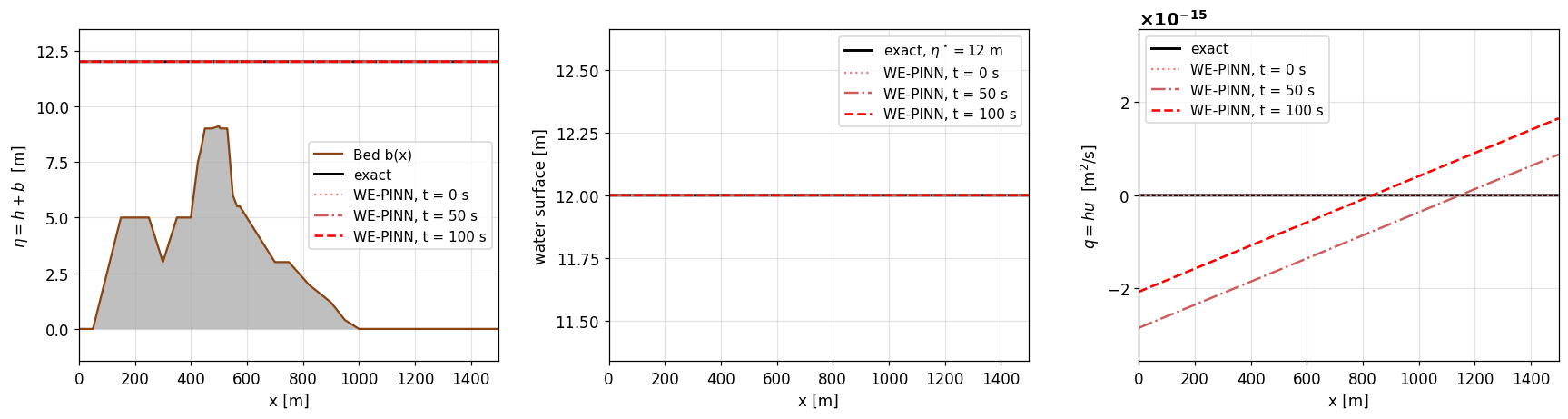}
    \caption{Run~B.}
  \end{subfigure}
  \caption{Lake at rest over the tabulated bed with $23$ kinks (WB4, $\eta^\star=12$~m, $T=100$~s); the prediction is shown at $t=0$, $50$ and $100$~s. The bed is continuous, so no jump is detected and no interface is marked. Left: free surface $\eta=h+b$ and bed $b$; centre: free surface magnified about the exact level $\eta^\star$; right: discharge $q$, whose exact value is zero. Black: exact equilibrium; red: WE-PINN prediction at $t=0$, $T/2$ and $T$ (three dash patterns). The factors in the corners of the centre and right panels give the magnification and offset of the vertical axis. (a) Run~A, austere budget; (b) Run~B, extended budget.}
  \label{fig:wb4}
\end{figure}

Criterion C6, the absence of growth in time, is examined through
Table~\ref{tab:wb_errors}, which reports the relative errors at the final time
$t=T$;
at $t=0$ they vanish identically by \eqref{eq:ansatz}. Since $b$ is exact and
steady, the absolute errors in $h$ and $\eta$ coincide, and $h$ is omitted. In
Run~A the errors at $T$ exceed those at $T/2$ by factors of $1.3$ to $2.9$, a
growth roughly linear in $t$; this reflects the blending $\tau(t)$ of \eqref{eq:ansatz}, which
hands progressively more of the state to the network, and not an accumulating
imbalance, which on a domain closed at both ends would compound. In Run~B the
errors at $T/2$ and $T$ are of the same order and decrease for $u$ and $q$ in
WB2 and WB4; the errors in $\eta$ of WB2 and WB4 are exactly zero at all $8000$
evaluation points (Figs.~\ref{fig:wb2} and~\ref{fig:wb4}). In Run~A the residual deviation is a smooth
structure organised by the bed, with extrema at the interfaces and, for WB3, at
the undetected jump, rather than a wave emitted from a discontinuity
(Figs.~\ref{fig:wb2} and~\ref{fig:wb3}).

\begin{table}[!htbp]
\centering\scriptsize
\setlength{\tabcolsep}{2.6pt}
\caption{Diagnostic D2: relative $L^1$, $L^2$ and $L^\infty$ errors at the final
time $t=T$, on $8000$ points. $\eta$ is normalised by the exact field; the exact
$u$ and $q$ vanish, so they are normalised by $\sqrt{gh_{\mathrm{ex}}}$ and
$h_{\mathrm{ex}}\sqrt{gh_{\mathrm{ex}}}$, $h_{\mathrm{ex}}=\eta^\star-b$.}
\label{tab:wb_errors}
\begin{tabular}{llccccccccc}
\toprule
& & \multicolumn{3}{c}{$\eta$} & \multicolumn{3}{c}{$u$} & \multicolumn{3}{c}{$q$} \\
\cmidrule(lr){3-5}\cmidrule(lr){6-8}\cmidrule(lr){9-11}
Run & Case & $L^1$ & $L^2$ & $L^\infty$ & $L^1$ & $L^2$ & $L^\infty$ & $L^1$ & $L^2$ & $L^\infty$ \\
\midrule
\multirow{4}{*}{A}
 & WB1 & $5.4{\times}10^{-5}$ & $6.0{\times}10^{-5}$ & $8.4{\times}10^{-5}$ & $4.4{\times}10^{-5}$ & $5.0{\times}10^{-5}$ & $9.7{\times}10^{-5}$ & $4.2{\times}10^{-5}$ & $4.5{\times}10^{-5}$ & $5.9{\times}10^{-5}$ \\
 & WB2 & $1.9{\times}10^{-5}$ & $2.2{\times}10^{-5}$ & $4.5{\times}10^{-5}$ & $5.1{\times}10^{-5}$ & $6.7{\times}10^{-5}$ & $1.4{\times}10^{-4}$ & $4.0{\times}10^{-5}$ & $4.5{\times}10^{-5}$ & $7.0{\times}10^{-5}$ \\
 & WB3 & $2.4{\times}10^{-5}$ & $3.5{\times}10^{-5}$ & $1.3{\times}10^{-4}$ & $4.1{\times}10^{-5}$ & $5.8{\times}10^{-5}$ & $1.4{\times}10^{-4}$ & $4.1{\times}10^{-5}$ & $5.7{\times}10^{-5}$ & $1.1{\times}10^{-4}$ \\
 & WB4 & $3.0{\times}10^{-6}$ & $3.3{\times}10^{-6}$ & $4.7{\times}10^{-6}$ & $3.3{\times}10^{-6}$ & $4.3{\times}10^{-6}$ & $1.0{\times}10^{-5}$ & $2.5{\times}10^{-6}$ & $2.6{\times}10^{-6}$ & $2.9{\times}10^{-6}$ \\
\midrule
\multirow{4}{*}{B}
 & WB1 & $2.2{\times}10^{-19}$ & $5.0{\times}10^{-18}$ & $1.1{\times}10^{-16}$ & $6.9{\times}10^{-19}$ & $7.1{\times}10^{-19}$ & $9.8{\times}10^{-19}$ & $6.8{\times}10^{-19}$ & $6.9{\times}10^{-19}$ & $9.1{\times}10^{-19}$ \\
 & WB2 & $0$ & $0$ & $0$ & $3.5{\times}10^{-18}$ & $5.2{\times}10^{-18}$ & $1.3{\times}10^{-17}$ & $2.4{\times}10^{-18}$ & $3.0{\times}10^{-18}$ & $6.3{\times}10^{-18}$ \\
 & WB3 & $5.5{\times}10^{-17}$ & $9.9{\times}10^{-17}$ & $3.0{\times}10^{-16}$ & $4.2{\times}10^{-17}$ & $7.3{\times}10^{-17}$ & $2.6{\times}10^{-16}$ & $4.1{\times}10^{-17}$ & $7.8{\times}10^{-17}$ & $2.6{\times}10^{-16}$ \\
 & WB4 & $0$ & $0$ & $0$ & $1.2{\times}10^{-17}$ & $1.5{\times}10^{-17}$ & $2.9{\times}10^{-17}$ & $1.0{\times}10^{-17}$ & $1.1{\times}10^{-17}$ & $1.6{\times}10^{-17}$ \\
\bottomrule
\end{tabular}
\end{table}

The interface coupling is examined in Table~\ref{tab:wb_iface}, which reports,
at the six detected interfaces, the largest jumps over $t\in\{T/2,T\}$. All three jumps $\jump{q}$, $\jump{H}$ and
$\jump{\eta}$ vanish at the lake at rest. In Run~B they are bounded by
$1.4\times10^{-15}$, $2.2\times10^{-16}$ and $2.2\times10^{-16}$, with
$\jump{H}=\jump{\eta}=0$ exactly at ten of the twelve sampled states, and the
traces reproduce the exact depths to all printed digits, including the contrast
$0.55$/$1.50$~m at $x=1.05$, where the bed jump exceeds the shallow-side depth.
In Run~A the jumps reach $3.9\times10^{-4}$ and $2.7\times10^{-4}$, and the
traces deviate from the exact depths in the fifth decimal place. $\jump{H}$ and
$\jump{\eta}$ coincide to the digits shown in all but one state, because the
velocity head $u^{2}/2g$ is below $2\times10^{-8}$~m in Run~A and $10^{-31}$~m in
Run~B; only $\jump{H}$ is tabulated.

\begin{table}[!htbp]
\centering\footnotesize
\setlength{\tabcolsep}{4pt}
\caption{Diagnostic D5a: largest interface jumps of the trained network over
$t\in\{T/2,T\}$, with the exact depths on either side. Both jumps vanish at the
lake at rest; $\jump{\eta}$ coincides with $\jump{H}$ to the digits shown.}
\label{tab:wb_iface}
\begin{tabular}{lccccccc}
\toprule
& & & & \multicolumn{2}{c}{Run~A} & \multicolumn{2}{c}{Run~B} \\
\cmidrule(lr){5-6}\cmidrule(lr){7-8}
Case & $x_0$ [m] & $\jump{b}$ [m] & exact $h^-/h^+$ [m] & $\max|\jump{q}|$ &
$\max|\jump{H}|$ & $\max|\jump{q}|$ & $\max|\jump{H}|$ \\
\midrule
WB2 & $0.00$ & $+1.00$ & $2.00/1.00$ & $1.4\times10^{-4}$ & $1.6\times10^{-5}$ & $2.7\times10^{-17}$ & $0$ \\
\addlinespace
\multirow{5}{*}{WB3}
 & $-0.65$ & $+0.35$ & $1.466/1.116$ & $3.9\times10^{-4}$ & $1.2\times10^{-4}$ & $2.3\times10^{-18}$ & $0$ \\
 & $-0.30$ & $-0.35$ & $1.150/1.500$ & $9.4\times10^{-5}$ & $2.2\times10^{-5}$ & $4.0\times10^{-17}$ & $0$ \\
 & $+1.05$ & $-0.95$ & $0.550/1.500$ & $9.2\times10^{-5}$ & $2.1\times10^{-4}$ & $1.4\times10^{-15}$ & $2.2\times10^{-16}$ \\
 & $+1.20$ & $+0.25$ & $1.500/1.250$ & $2.5\times10^{-4}$ & $2.7\times10^{-4}$ & $2.7\times10^{-16}$ & $0$ \\
 & $+1.40$ & $-0.25$ & $1.250/1.500$ & $3.9\times10^{-5}$ & $3.8\times10^{-5}$ & $1.6\times10^{-16}$ & $2.2\times10^{-16}$ \\
\bottomrule
\end{tabular}
\end{table}

Finally, Table~\ref{tab:wb_mass} evaluates the three terms of \eqref{eq:mass_balance}
independently on a fine grid, a posteriori. Their sum, which must vanish,
closes to between $1.7\times10^{-6}$ and $4.4\times10^{-5}$ of the initial
water volume in Run~A and to between $10^{-19}$ and $4\times10^{-17}$ in Run~B. Where an interface exists,
the interface defect is of the same order as the other two terms, consistent
with Proposition~\ref{prop:mass}, in which the mass error is carried by
$\jump{q}$. The boundary discharges, constrained directly by
\eqref{eq:bc_dirichlet}, stay below $3.8\times10^{-4}$ and
$3.6\times10^{-15}$~m$^{2}$\,s$^{-1}$; since the exact discharge is zero, they also
confirm that no spurious wave reaches the ends of the domain.

\begin{table}[!htbp]
\centering\footnotesize
\setlength{\tabcolsep}{3.5pt}
\caption{Diagnostic D5b: the three terms of the global mass balance
\eqref{eq:mass_balance} [m$^{2}$], their sum normalised by the initial volume
$V_0$, and the largest boundary discharge over $[0,T]$ [m$^{2}$\,s$^{-1}$]. All
vanish at the exact equilibrium. \naIfaceNote{}Not applicable: WB1 and WB4 have
no bed discontinuity and hence no interface (Table~\ref{tab:wb_cases}).}
\label{tab:wb_mass}
\begin{tabular}{llcccccc}
\toprule
Run & Case & storage & outflux & interface & $V_0$ [m$^{2}$] & sum$/V_0$ &
$\max_t|q(x_{\mathrm L,R})|$ \\
\midrule
\multirow{4}{*}{A}
 & WB1 & $-5.5\times10^{-4}$ & $2.4\times10^{-5}$ & \naIface & $12.00$ & $4.4\times10^{-5}$ & $3.5\times10^{-5}$ \\
 & WB2 & $1.9\times10^{-4}$ & $1.3\times10^{-5}$ & $5.4\times10^{-5}$ & $30.00$ & $8.7\times10^{-6}$ & $1.2\times10^{-4}$ \\
 & WB3 & $-3.6\times10^{-5}$ & $9.8\times10^{-6}$ & $7.1\times10^{-5}$ & $3.52$ & $1.3\times10^{-5}$ & $2.6\times10^{-4}$ \\
 & WB4 & $-4.2\times10^{-2}$ & $1.9\times10^{-2}$ & \naIface & $14191$ & $1.7\times10^{-6}$ & $3.8\times10^{-4}$ \\
\midrule
\multirow{4}{*}{B}
 & WB1 & $-2.8\times10^{-18}$ & $1.2\times10^{-18}$ & \naIface & $12.00$ & $1.3\times10^{-19}$ & $1.5\times10^{-18}$ \\
 & WB2 & $0$ & $7.4\times10^{-17}$ & $1.5\times10^{-17}$ & $30.00$ & $3.0\times10^{-18}$ & $7.7\times10^{-17}$ \\
 & WB3 & $2.5\times10^{-16}$ & $4.1\times10^{-17}$ & $-1.6\times10^{-16}$ & $3.52$ & $3.6\times10^{-17}$ & $2.7\times10^{-16}$ \\
 & WB4 & $0$ & $3.7\times10^{-13}$ & \naIface & $14191$ & $2.6\times10^{-17}$ & $3.6\times10^{-15}$ \\
\bottomrule
\end{tabular}
\end{table}

\FloatBarrier
\subsubsection{Summary and limitations}
\label{sec:wb_summary}

The suite separates the two statements of Section~\ref{sec:protocol}. The
residual operator annihilates the lake at rest in the sense of
Theorem~\ref{thm:wb}: every residual is zero bitwise or at a round-off floor with
the signatures predicted by Proposition~\ref{prop:roundoff}, including on cells
that violate the hypotheses under which the residual was constructed. The
trained network approaches the equilibrium to a level set by the optimisation
budget alone, which is the behaviour described by Corollary~\ref{cor:minimiser}.
Three limitations apply. (i) None of the four time windows covers a full
gravity-wave transit of its domain (WB4, the longest, covers between $0.36$ and $0.72$
of one), so the statements about the trained solution are bounded in time; those
about the residual operator involve no time integration. (ii) Exact preservation
is a property of the objective and its minimiser, obtained only to the extent
that the optimiser reaches it, not of a marching update rule. (iii) In
\texttt{float32}, the working precision of Table~\ref{tab:global_params}, the
round-off floor is nine decades higher than in the \texttt{float64} runs reported
here (Table~\ref{tab:d1}).

\FloatBarrier

\subsection{Riemann problems at a bottom step}
\label{sec:results_ab}

\subsubsection{Benchmark and set-up}
\label{sec:ab_config}

The seven problems of this subsection, RP1--RP7, are taken from the wet-bed
Riemann problems at a bottom step of Alcrudo and
Benkhaldoun~\cite{AlcrudoBenkhaldoun2001}, whose closed-form solutions
serve as reference. All share the bed $b=0$ for $x<0$ and $b=1$~m for $x\ge0$,
with the initial conditions $(h_l,u_l)$ for $x<0$ and $(h_r,u_r)$ for $x\ge0$, so that
the initial discontinuity and the bed discontinuity coincide. The exact solution
consists of simple waves travelling away from the step together with the
stationary step transition of Section~\ref{sec:sst} pinned to it, across which
$h$, $u$ and $\eta$ jump while $q$ and $H$ are conserved. The amplitude of that
jump is fixed by the cubic \eqref{eq:sst_cubic}, which has two admissible roots
(Lemma~\ref{lem:conjugate}). A correct solution therefore requires that the
head condition be enforced instead of the continuity of $\eta$ that a single
smooth network would impose (Theorem~\ref{thm:no_continuous}), and that the
physically relevant root be selected. Neither requirement is stated in the
objective, which imposes only $\jump{q}=0$ and $\jump{H}=0$ at $x=0$; both must
emerge from the minimisation.

\begin{figure}[!htbp]
  \centering
  \includegraphics[width=\textwidth]{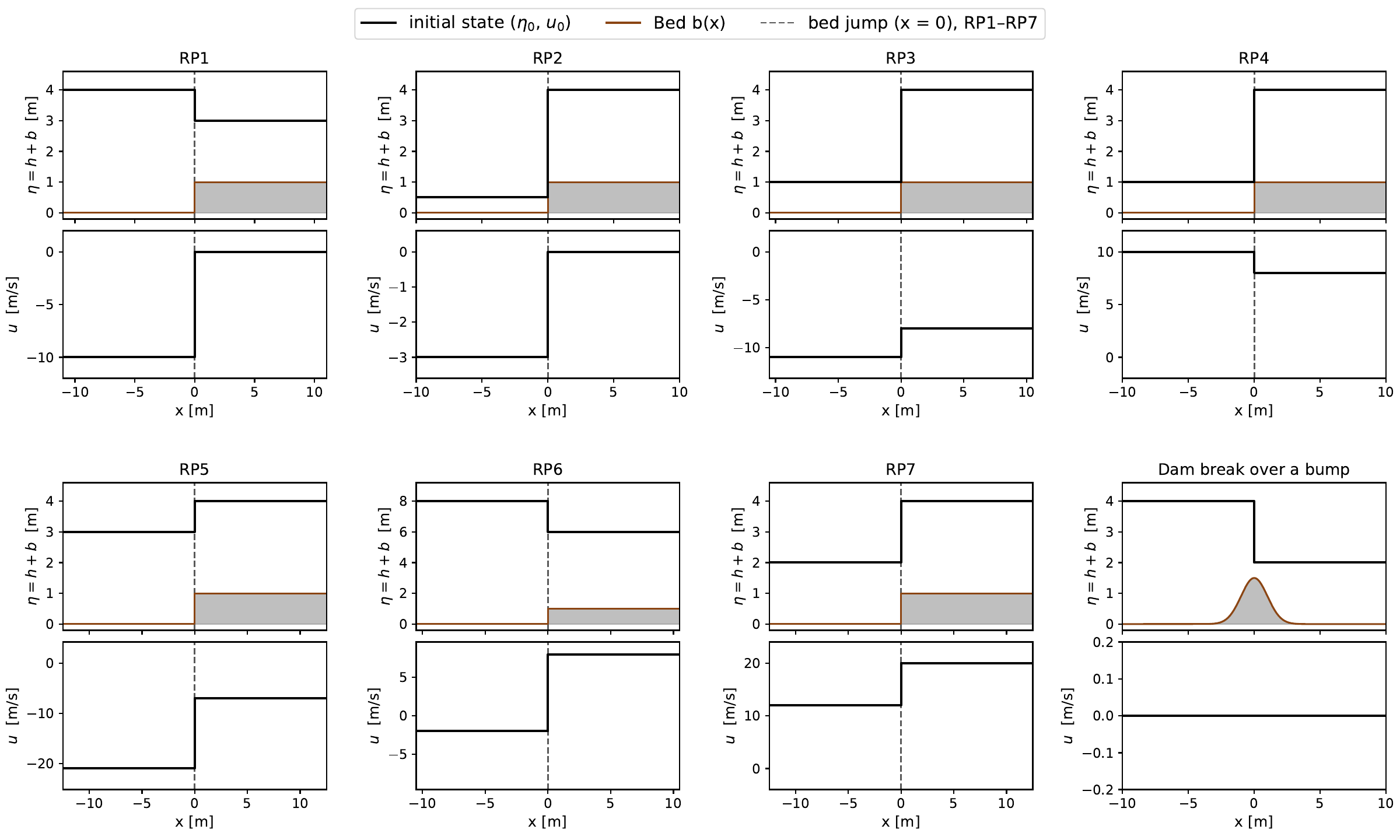}
  \caption{Initial conditions of the seven Riemann problems at a bottom step,
  RP1--RP7, and of the smooth-bed control of Section~\ref{sec:res_bump}. Each
  panel shows the bed, the initial free surface and the initial velocity.
  RP1--RP7 share $b=0$ for $x<0$ and $b=1$~m for $x\ge0$ and differ only in
  $(h_l,u_l,h_r,u_r)$ and in the final time (Table~\ref{tab:ab_config}); the
  dam break over a bump (last panel) has the smooth bed
  $b=\frac32e^{-x^{2}/2}$, no bed jump, and water initially at rest.}
  \label{fig:rp_overview}
\end{figure}

The seven problems are listed in Table~\ref{tab:ab_config} and shown in
Fig.~\ref{fig:rp_overview}. They cover the three characters of the transition
(subcritical in RP6, critical in RP1 and RP2, supercritical in RP3, RP4, RP5 and
RP7), wave patterns with and without shocks, a step twice the upstream depth
(RP2), Froude numbers up to $3.2$ (RP4), velocities up to $20$~m\,s$^{-1}$ (RP7)
and both signs of the free-surface jump. The ratio $s_{\max}h_t/h_x$ measures how
many cell widths the fastest wave crosses per cell height; it ranges from $5.4$
to $26.4$ and, by the quadrature argument following
Proposition~\ref{prop:consistency}, controls the resolution of moving shocks.

\begin{table}[!htbp]
\centering\footnotesize
\setlength{\tabcolsep}{4pt}
\caption{The seven Riemann problems, with the step $b=0\to1$~m at $x=0$ and
domain $[-L,L]$, taken from \cite{AlcrudoBenkhaldoun2001}. Wave pattern: $R^k$ and $S^k$ denote a
$k$-rarefaction and a $k$-shock, and $|\,\mathrm{SST}\,|$ the stationary step
transition, with its character (subcritical, critical, supercritical). $U_2$ and
$U_3$ are the exact traces $(h,u)$ at $x=0^-$ and $x=0^+$, $F_2$, $F_3$ the
corresponding Froude numbers and $\jump{\eta}$ the exact free-surface jump. All
values are properties of the benchmark, not results of the method.}
\label{tab:ab_config}
\begin{tabular}{llccccll}
\toprule
Case & Name & $(h_l,u_l)$ & $(h_r,u_r)$ & $T$ [s] & $L$ [m] & Wave pattern & $s_{\max}h_t/h_x$ \\
\midrule
RP1 & critical transition & $(4,-10)$ & $(2,0)$ & $0.5$ & $11$ & $R^1S^2\,|\,\mathrm{SST_{crit}}\,|\,R^2$ & $16.3$ \\
RP2 & high step, shallow inflow & $(0.5,-3)$ & $(3,0)$ & $1.0$ & $10$ & $S^1S^2\,|\,\mathrm{SST_{crit}}\,|\,R^2$ & $5.4$ \\
RP3 & undisturbed far side & $(1,-11)$ & $(3,-8)$ & $0.5$ & $10.5$ & $S^1R^2\,|\,\mathrm{SST_{sup}}\,|$ & $14.1$ \\
RP4 & strongly supercritical & $(1,10)$ & $(3,8)$ & $0.5$ & $10$ & $|\,\mathrm{SST_{sup}}\,|\,S^1R^2$ & $13.4$ \\
RP5 & fast rarefactions & $(3,-21)$ & $(3,-7)$ & $0.35$ & $12.5$ & $R^1R^2\,|\,\mathrm{SST_{sup}}\,|$ & $26.4$ \\
RP6 & subcritical, deep flow & $(8,-2)$ & $(5,8)$ & $0.5$ & $10.5$ & $R^1\,|\,\mathrm{SST_{sub}}\,|\,R^2$ & $15.0$ \\
RP7 & high-velocity inflow & $(2,12)$ & $(3,20)$ & $0.35$ & $12.5$ & $|\,\mathrm{SST_{sup}}\,|\,R^1R^2$ & $25.4$ \\
\bottomrule
\end{tabular}

\vspace{6pt}

\begin{tabular}{lcccc}
\toprule
Case & $U_2=(h,u)(0^-)$ & $U_3=(h,u)(0^+)$ & $F_2$ / $F_3$ & $\jump{\eta}$ [m] \\
\midrule
RP1 & $(0.4295,-6.1118)$ & $(0.8889,-2.9530)$ & $-2.98$ / $-1.00$ & $+1.4594$ \\
RP2 & $(0.7212,-6.6866)$ & $(1.3333,-3.6166)$ & $-2.51$ / $-1.00$ & $+1.6122$ \\
RP3 & $(2.4770,-9.6893)$ & $(3.0000,-8.0000)$ & $-1.97$ / $-1.47$ & $+1.5230$ \\
RP4 & $(1.0000,10.0000)$ & $(1.1341,8.8175)$ & $+3.19$ / $+2.64$ & $+1.1341$ \\
RP5 & $(2.3195,-9.0539)$ & $(3.0000,-7.0000)$ & $-1.90$ / $-1.29$ & $+1.6805$ \\
RP6 & $(4.2906,2.7423)$ & $(2.7205,4.3250)$ & $+0.42$ / $+0.84$ & $-0.5701$ \\
RP7 & $(2.0000,12.0000)$ & $(2.1839,10.9897)$ & $+2.71$ / $+2.37$ & $+1.1839$ \\
\bottomrule
\end{tabular}
\end{table}

Every run uses two sub-domains, obtained from the bathymetry alone by the
detector of Section~\ref{sec:detection}, which locates the jump at $x=0$ to
machine precision and returns $\jump{b}=+1.000000$ in all seven cases; the
network is never told that a step is present. The working precision is
\texttt{float32} and the parameters are those of
Table~\ref{tab:global_params}, with the transmissive boundary condition
\eqref{eq:bc_transmissive}; the case-dependent settings are listed in
Table~\ref{tab:ab_settings}. The domain is sized from the fastest exact wave with
a $35\%$ margin, so that no wave reaches a boundary before $T$ and the
transmissive condition is exact; this is part of the test definition, and the
boundary diagnostic of Section~\ref{sec:ab_cross} verifies it a posteriori. The
reference is the closed-form solution of \cite{AlcrudoBenkhaldoun2001}, computed as described
in~\ref{app:exact} and never
used in the loss. Errors are relative, the norm of
$\phi_\theta-\phi_{\mathrm{ex}}$ divided by the same norm of $\phi_{\mathrm{ex}}$,
on a uniform grid of $8000$ points, for $\phi\in\{h,u,\eta,q,F\}$. The $L^1$ norm
measures the accuracy of the solution away from discontinuities, whereas the
$L^\infty$ norm measures the thickness of the smeared shock. None of the norms
tests the interface, which occupies a single abscissa; the interface is judged by
the traces and jumps of Table~\ref{tab:ab_d2}.

\begin{table}[!htbp]
\centering\footnotesize
\setlength{\tabcolsep}{4pt}
\caption{Case-dependent settings of the Riemann problems. $N$ is the requested
number of control volumes, $N_\Gamma$ the number of interface and boundary time
intervals ($c_\tau=2N_\Gamma$), $Q$ the Gauss--Legendre order, $a_0$ the initial
activation slope; every sub-domain network has four hidden layers of $80$ units.
L-BFGS is run as a single round ($N_{\mathrm{rounds}}=1$) of at most
$N_{\mathrm{LBFGS}}$ iterations.}
\label{tab:ab_settings}
\begin{tabular}{lccccccc}
\toprule
 & RP1 & RP2 & RP3 & RP4 & RP5 & RP6 & RP7 \\
\midrule
$N$; $N_\Gamma$       & $10^4$; $50$ & $10^4$; $100$ & $10^4$; $100$ & $10^4$; $100$ & $10^4$; $100$ & $10^4$; $100$ & $10^4$; $65$ \\
$Q$; $a_0$            & $12$; $0.90$ & $12$; $0.90$ & $12$; $0.90$ & $12$; $0.90$ & $12$; $0.90$ & $12$; $0.90$ & $12$; $0.90$ \\
$N_{\mathrm{Adam}}$   & $2000$ & $2000$ & $2000$ & $2000$ & $2000$ & $2000$ & $2000$ \\
$N_{\mathrm{LBFGS}}$  & $5000$ & $5000$ & $5000$ & $5000$ & $5000$ & $5000$ & $5000$ \\
\bottomrule
\end{tabular}
\end{table}

\FloatBarrier
\subsubsection{Accuracy of the wave pattern}
\label{sec:ab_accuracy}

\begin{figure}[!htbp]
  \centering
  \includegraphics[width=0.85\textwidth]{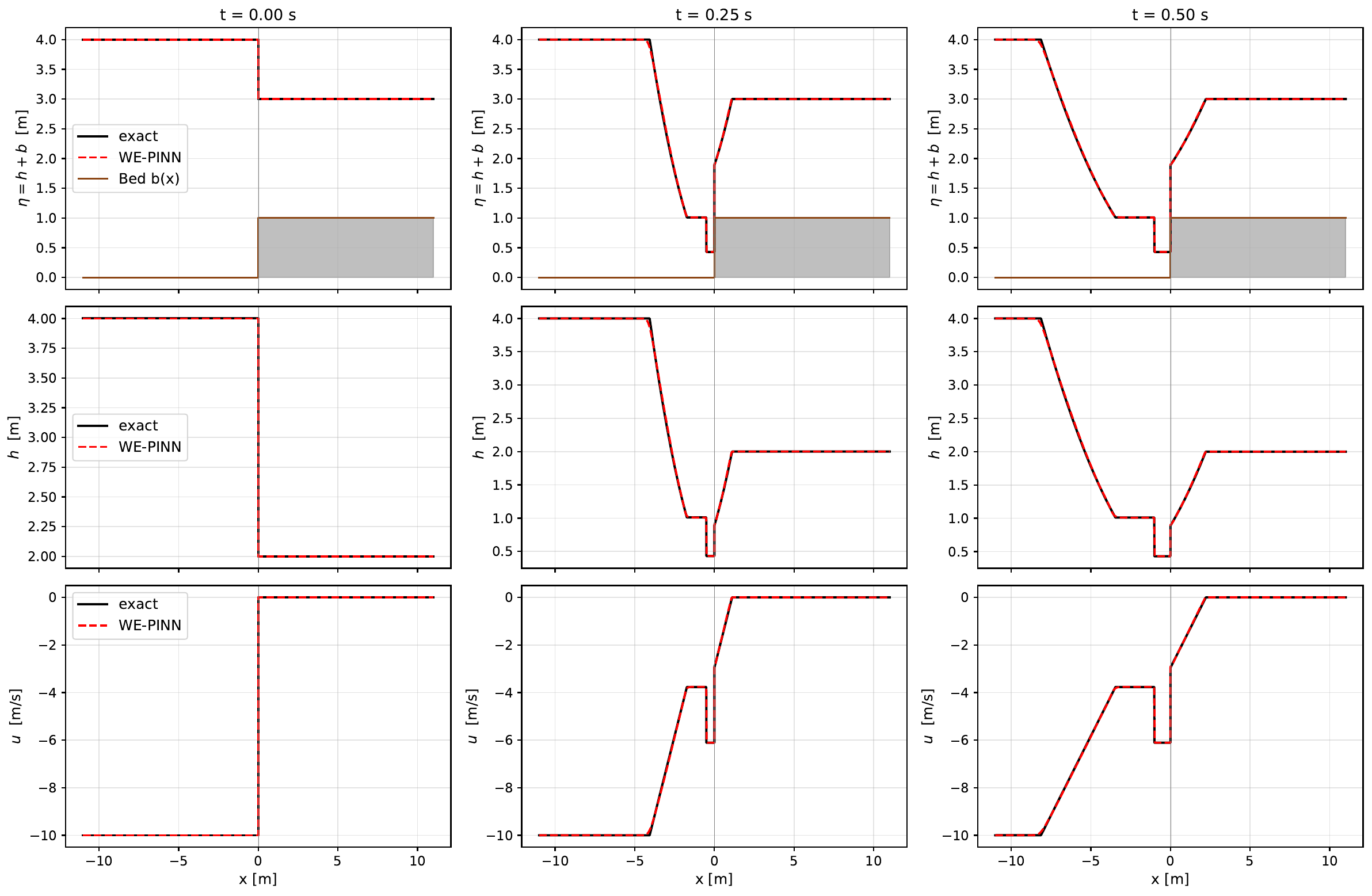}
  \caption{RP1, critical transition ($h=4/2$~m, $u=-10/0$~m\,s$^{-1}$): a $1$-rarefaction and a $2$-shock on the low bed, the critical transition and a $2$-rarefaction on the high bed. Free surface $\eta=h+b$ (top), depth (middle) and velocity (bottom) at $t=0$, $T/2$ and $T$. Black: exact solution \cite{AlcrudoBenkhaldoun2001}; red dashed: WE-PINN; grey: bed.}
  \label{fig:rp1}
\end{figure}

\begin{figure}[!htbp]
  \centering
  \includegraphics[width=0.85\textwidth]{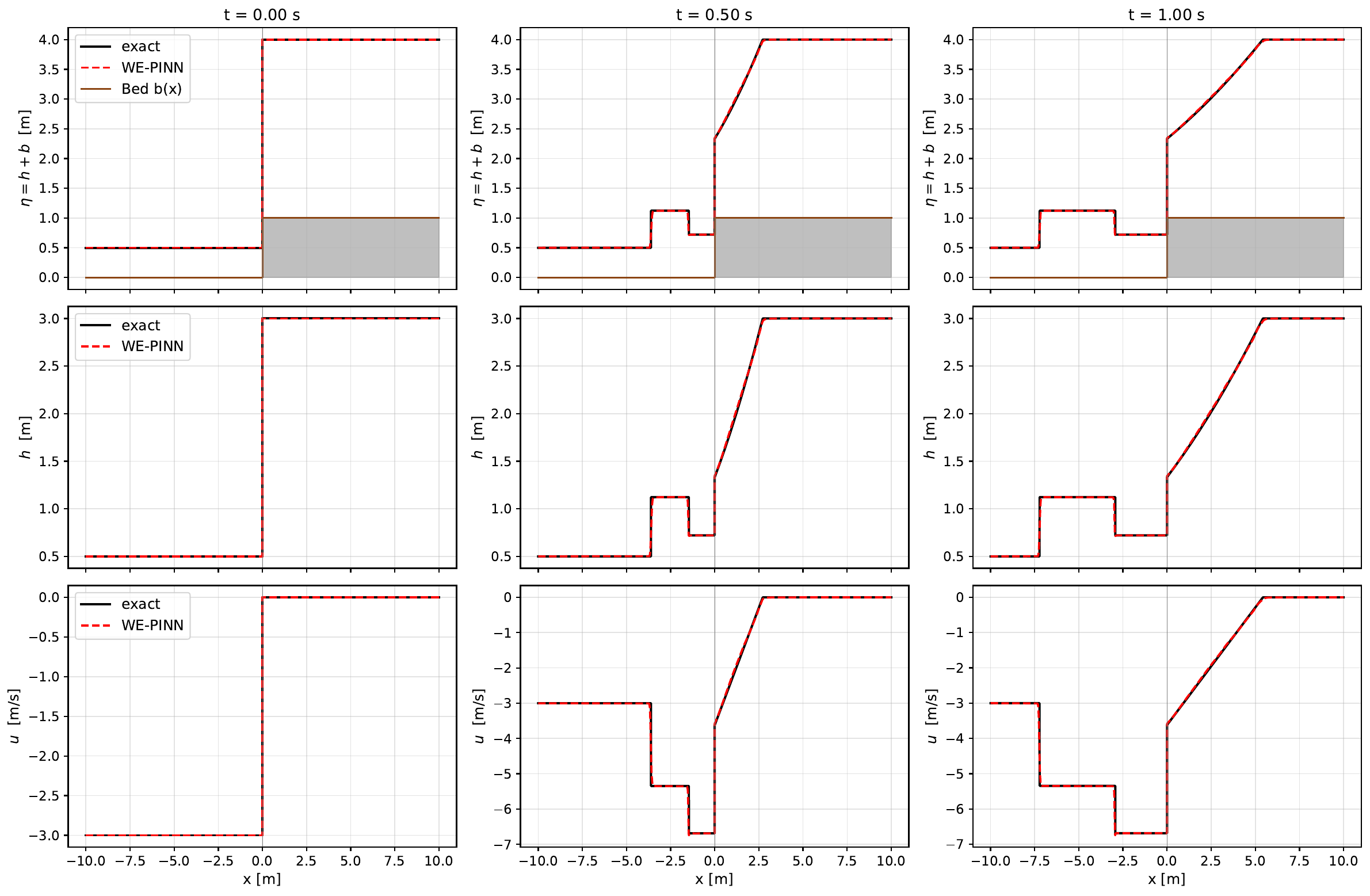}
  \caption{RP2, step twice the upstream depth ($h=0.5/3$~m, $u=-3/0$~m\,s$^{-1}$): two shocks on the low bed, a critical transition and a $2$-rarefaction on the high bed. Free surface $\eta=h+b$ (top), depth (middle) and velocity (bottom) at $t=0$, $T/2$ and $T$. Black: exact solution \cite{AlcrudoBenkhaldoun2001}; red dashed: WE-PINN; grey: bed.}
  \label{fig:rp2}
\end{figure}

\begin{figure}[!htbp]
  \centering
  \includegraphics[width=0.85\textwidth]{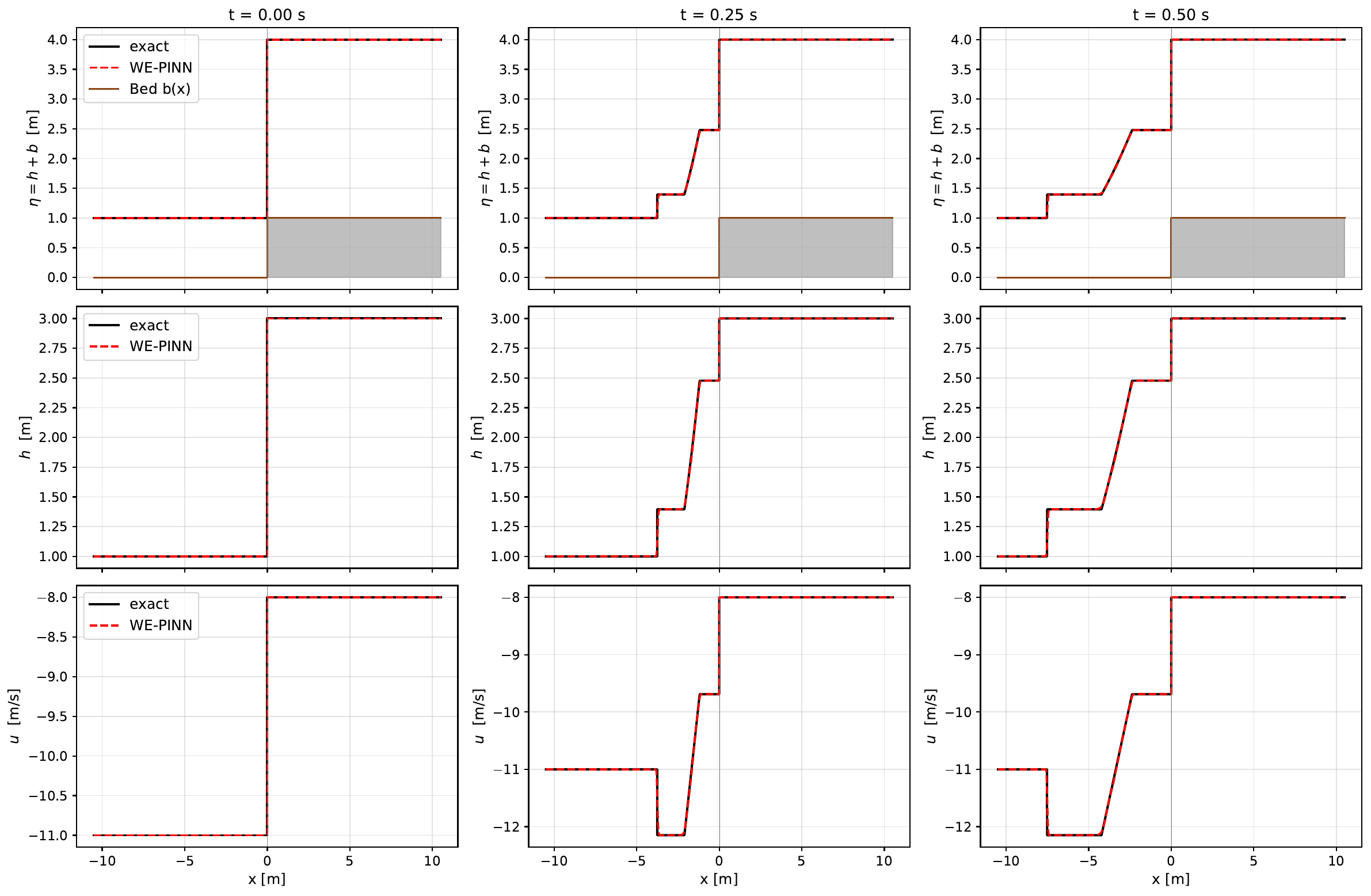}
  \caption{RP3, undisturbed far side ($h=1/3$~m, $u=-11/-8$~m\,s$^{-1}$): a $1$-shock and a $2$-rarefaction on the low bed and a supercritical transition, with the right state untouched. Free surface $\eta=h+b$ (top), depth (middle) and velocity (bottom) at $t=0$, $T/2$ and $T$. Black: exact solution \cite{AlcrudoBenkhaldoun2001}; red dashed: WE-PINN; grey: bed.}
  \label{fig:rp3}
\end{figure}

\begin{figure}[!htbp]
  \centering
  \includegraphics[width=0.85\textwidth]{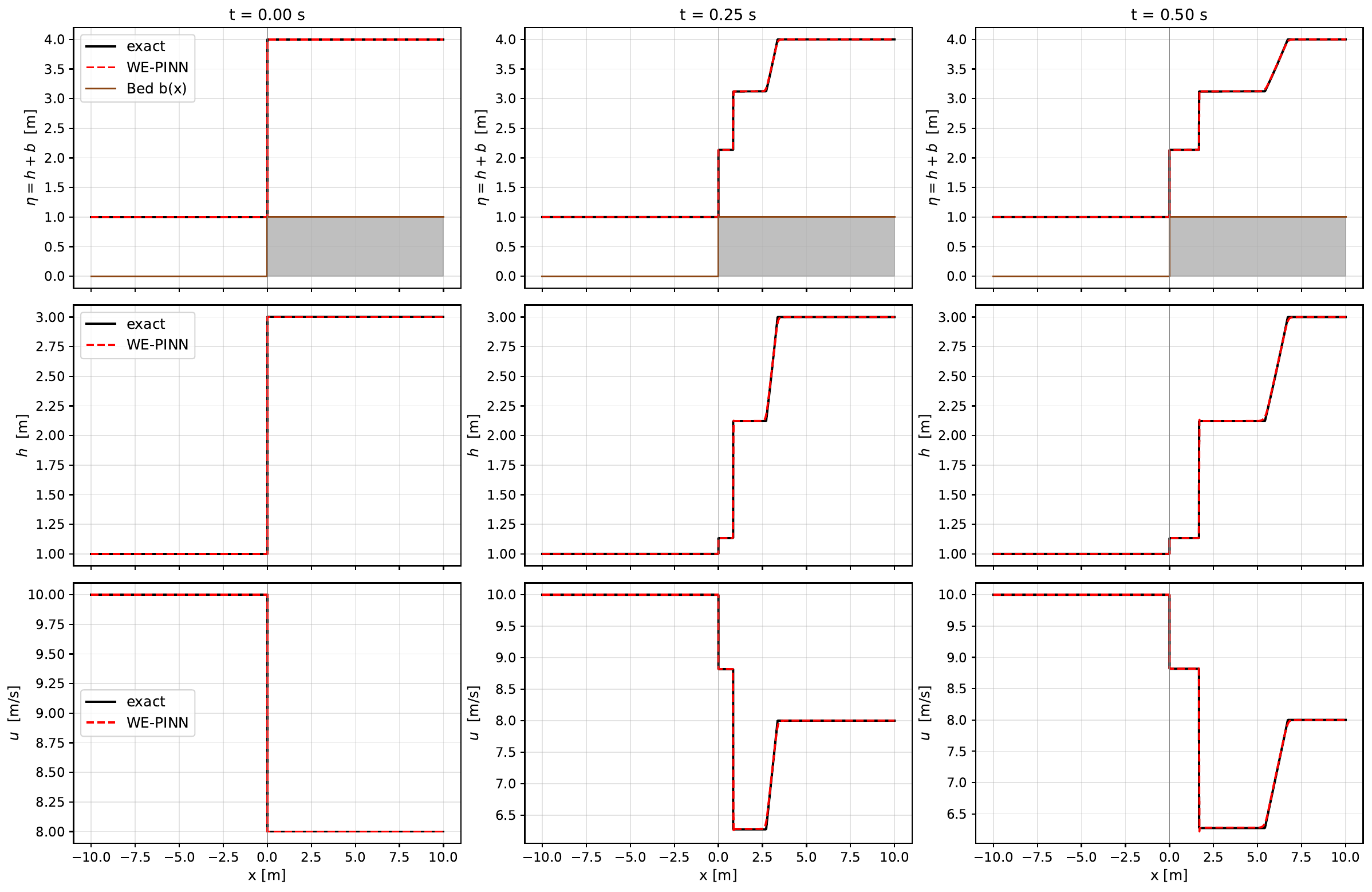}
  \caption{RP4, strongly supercritical flow ($h=1/3$~m, $u=10/8$~m\,s$^{-1}$, $F_2=+3.19$): supercritical transition followed by a $1$-shock and a $2$-rarefaction on the high bed. Free surface $\eta=h+b$ (top), depth (middle) and velocity (bottom) at $t=0$, $T/2$ and $T$. Black: exact solution \cite{AlcrudoBenkhaldoun2001}; red dashed: WE-PINN; grey: bed.}
  \label{fig:rp4}
\end{figure}

\begin{figure}[!htbp]
  \centering
  \includegraphics[width=0.85\textwidth]{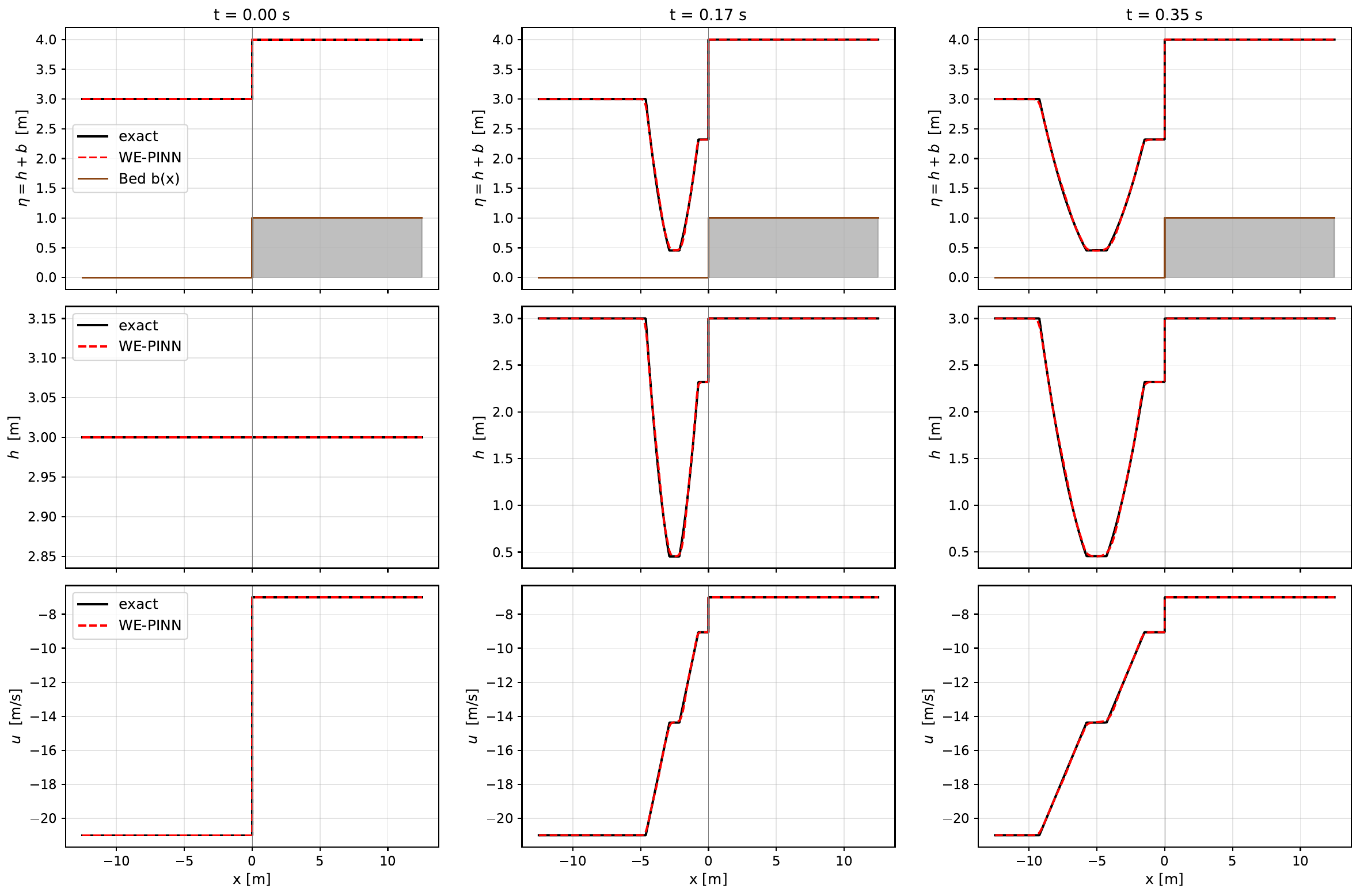}
  \caption{RP5, fast rarefactions ($h=3/3$~m, $u=-21/-7$~m\,s$^{-1}$): two rarefactions on the low bed and a supercritical transition, the fastest case of the subset. Free surface $\eta=h+b$ (top), depth (middle) and velocity (bottom) at $t=0$, $T/2$ and $T$. Black: exact solution \cite{AlcrudoBenkhaldoun2001}; red dashed: WE-PINN; grey: bed.}
  \label{fig:rp5}
\end{figure}

\begin{figure}[!htbp]
  \centering
  \includegraphics[width=0.85\textwidth]{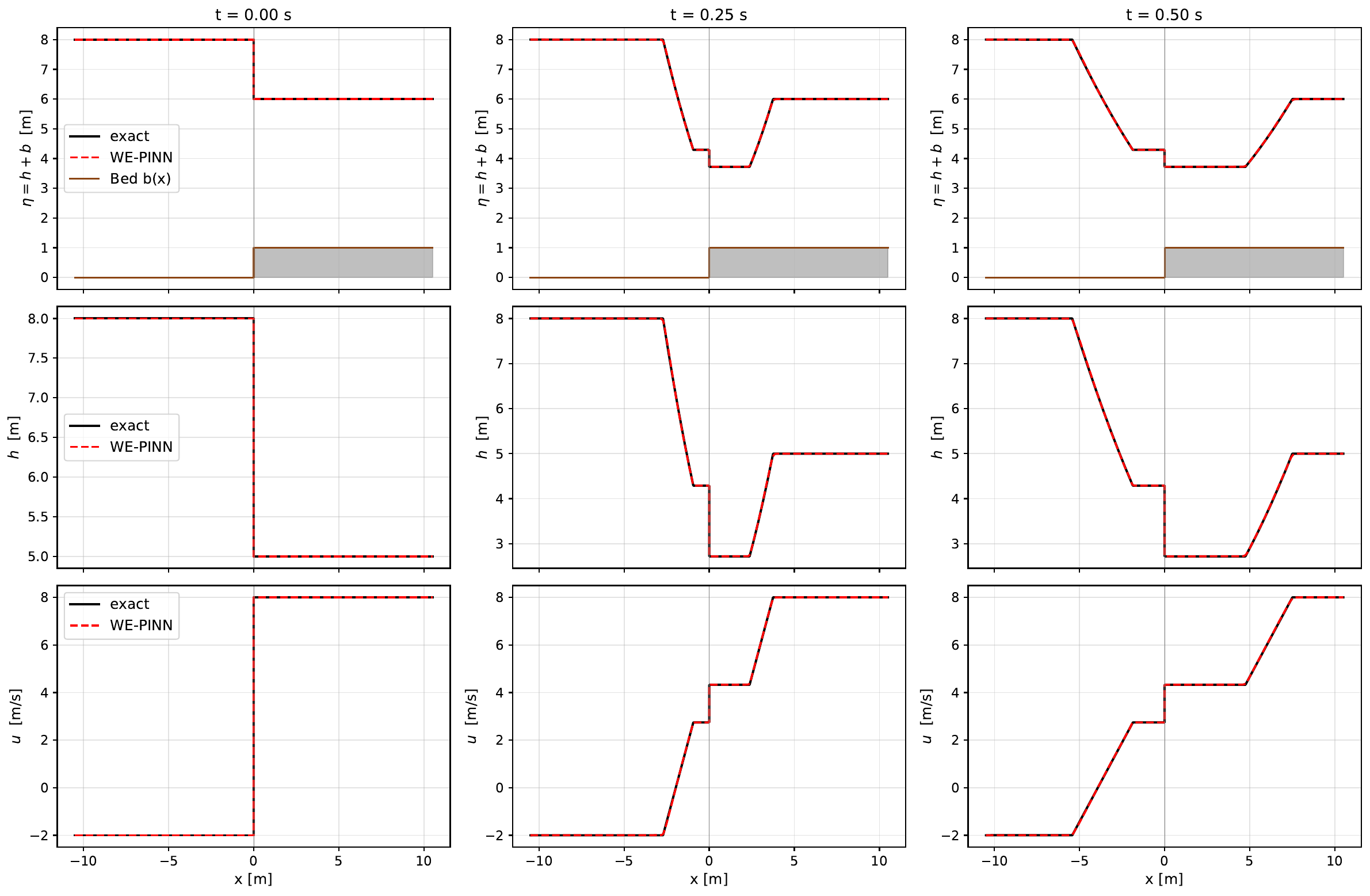}
  \caption{RP6, subcritical transition in deep flow ($h=8/5$~m, $u=-2/8$~m\,s$^{-1}$): two rarefactions separated by a subcritical transition with a negative free-surface jump. Free surface $\eta=h+b$ (top), depth (middle) and velocity (bottom) at $t=0$, $T/2$ and $T$. Black: exact solution \cite{AlcrudoBenkhaldoun2001}; red dashed: WE-PINN; grey: bed.}
  \label{fig:rp6}
\end{figure}

\begin{figure}[!htbp]
  \centering
  \includegraphics[width=0.85\textwidth]{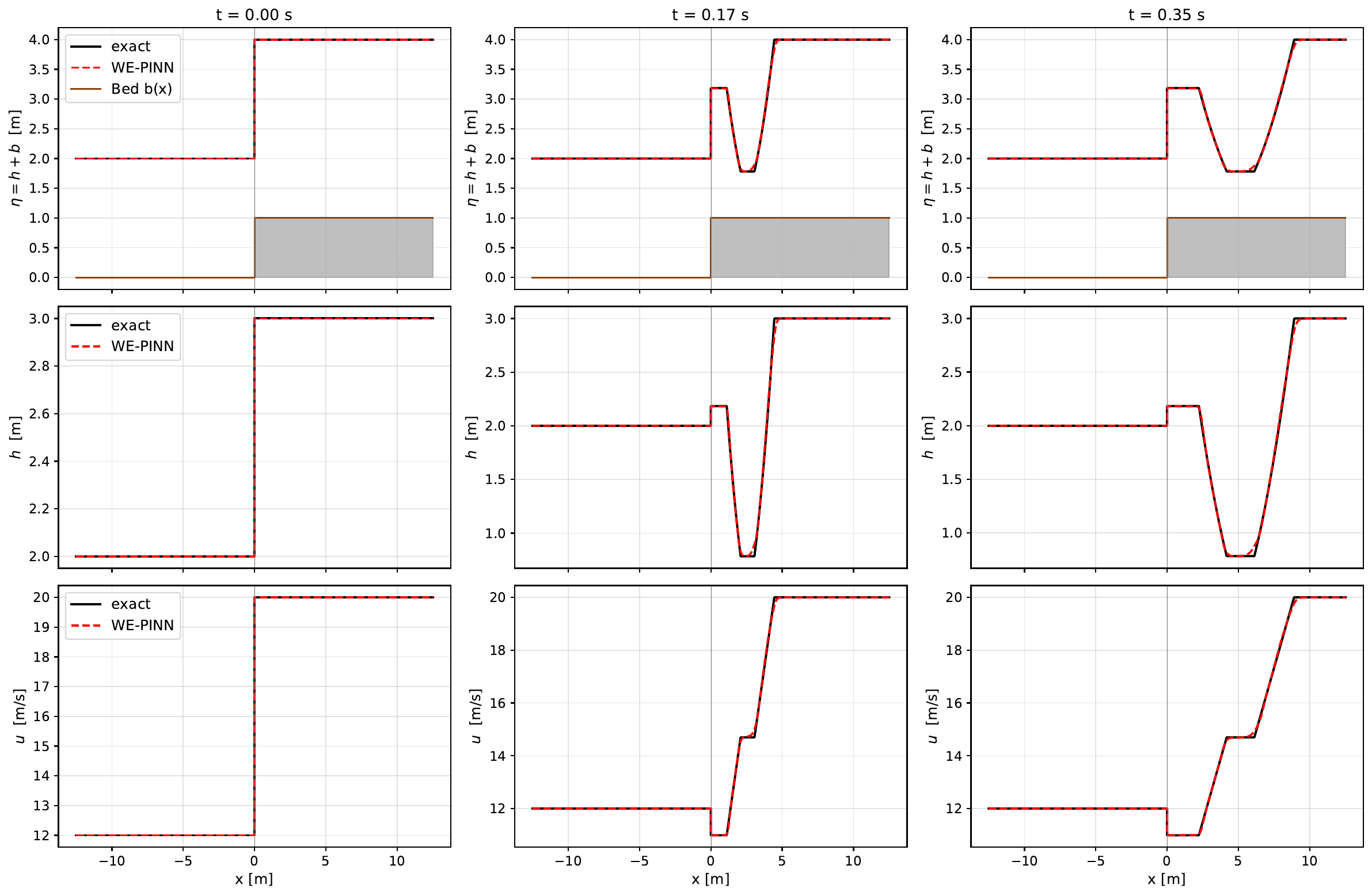}
  \caption{RP7, high-velocity inflow ($h=2/3$~m, $u=12/20$~m\,s$^{-1}$): supercritical transition with both rarefactions on the high bed. Free surface $\eta=h+b$ (top), depth (middle) and velocity (bottom) at $t=0$, $T/2$ and $T$. Black: exact solution \cite{AlcrudoBenkhaldoun2001}; red dashed: WE-PINN; grey: bed.}
  \label{fig:rp7}
\end{figure}

The trained solutions are compared with the exact ones in
Figs.~\ref{fig:rp1}--\ref{fig:rp7}, and the errors at the final time are given in
Table~\ref{tab:ab_err}. In every case the number, type and position of the waves
and the intermediate plateaux are recovered, and the stationary transition stays
pinned at $x=0$ throughout. On the free surface the relative $L^1$ error at
$t=T$ lies between $6.9\times10^{-4}$ (RP6) and $3.5\times10^{-3}$ (RP1), one to
four parts in a thousand, with the depth slightly larger and the discharge larger
again. Between $T/2$ and $T$ this error changes by factors between $0.96$ and
$1.33$: it does not accumulate in time, and in RP3 it decreases while the
solution develops.

\begin{table}[!htbp]
\centering\scriptsize
\setlength{\tabcolsep}{3.5pt}
\caption{Relative errors of the trained network at the final time $t=T$ against
the exact solution of \cite{AlcrudoBenkhaldoun2001}, on $8000$ points. Each entry
is the norm of the error divided by the same norm of the exact field; at $t=0$
all errors vanish identically by \eqref{eq:ansatz}.}
\label{tab:ab_err}
\begin{tabular}{llccccc}
\toprule
Case & Norm & $h$ & $u$ & $\eta$ & $q$ & $F$ \\
\midrule
\multirow{3}{*}{RP1} & $L^1$ & $4.4\times10^{-3}$ & $6.2\times10^{-3}$ & $3.5\times10^{-3}$ & $9.4\times10^{-3}$ & $7.2\times10^{-3}$ \\
 & $L^2$ & $1.4\times10^{-2}$ & $2.1\times10^{-2}$ & $1.1\times10^{-2}$ & $1.3\times10^{-2}$ & $6.8\times10^{-2}$ \\
 & $L^\infty$ & $1.5\times10^{-1}$ & $2.5\times10^{-1}$ & $1.5\times10^{-1}$ & $4.3\times10^{-2}$ & $6.3\times10^{-1}$ \\
\cmidrule(l){2-7}
\multirow{3}{*}{RP2} & $L^1$ & $2.6\times10^{-3}$ & $4.7\times10^{-3}$ & $2.0\times10^{-3}$ & $6.0\times10^{-3}$ & $5.2\times10^{-3}$ \\
 & $L^2$ & $9.6\times10^{-3}$ & $1.9\times10^{-2}$ & $7.1\times10^{-3}$ & $2.6\times10^{-2}$ & $2.2\times10^{-2}$ \\
 & $L^\infty$ & $1.4\times10^{-1}$ & $2.8\times10^{-1}$ & $1.0\times10^{-1}$ & $4.8\times10^{-1}$ & $3.8\times10^{-1}$ \\
\cmidrule(l){2-7}
\multirow{3}{*}{RP3} & $L^1$ & $1.1\times10^{-3}$ & $6.5\times10^{-4}$ & $8.7\times10^{-4}$ & $1.3\times10^{-3}$ & $1.0\times10^{-3}$ \\
 & $L^2$ & $5.2\times10^{-3}$ & $3.6\times10^{-3}$ & $4.2\times10^{-3}$ & $8.7\times10^{-3}$ & $3.9\times10^{-3}$ \\
 & $L^\infty$ & $7.8\times10^{-2}$ & $4.9\times10^{-2}$ & $5.8\times10^{-2}$ & $1.5\times10^{-1}$ & $3.9\times10^{-2}$ \\
\cmidrule(l){2-7}
\multirow{3}{*}{RP4} & $L^1$ & $1.1\times10^{-3}$ & $8.4\times10^{-4}$ & $8.6\times10^{-4}$ & $1.2\times10^{-3}$ & $1.0\times10^{-3}$ \\
 & $L^2$ & $1.0\times10^{-2}$ & $5.8\times10^{-3}$ & $7.3\times10^{-3}$ & $5.6\times10^{-3}$ & $1.0\times10^{-2}$ \\
 & $L^\infty$ & $3.0\times10^{-1}$ & $2.5\times10^{-1}$ & $2.3\times10^{-1}$ & $1.2\times10^{-1}$ & $3.8\times10^{-1}$ \\
\cmidrule(l){2-7}
\multirow{3}{*}{RP5} & $L^1$ & $3.3\times10^{-3}$ & $1.8\times10^{-3}$ & $2.8\times10^{-3}$ & $5.2\times10^{-3}$ & $6.6\times10^{-3}$ \\
 & $L^2$ & $6.5\times10^{-3}$ & $3.3\times10^{-3}$ & $5.3\times10^{-3}$ & $1.1\times10^{-2}$ & $1.3\times10^{-2}$ \\
 & $L^\infty$ & $4.1\times10^{-2}$ & $1.1\times10^{-2}$ & $3.1\times10^{-2}$ & $5.1\times10^{-2}$ & $3.9\times10^{-2}$ \\
\cmidrule(l){2-7}
\multirow{3}{*}{RP6} & $L^1$ & $7.6\times10^{-4}$ & $1.5\times10^{-3}$ & $6.9\times10^{-4}$ & $2.4\times10^{-3}$ & $1.0\times10^{-3}$ \\
 & $L^2$ & $1.5\times10^{-3}$ & $2.7\times10^{-3}$ & $1.4\times10^{-3}$ & $4.7\times10^{-3}$ & $1.7\times10^{-3}$ \\
 & $L^\infty$ & $7.2\times10^{-3}$ & $1.0\times10^{-2}$ & $7.2\times10^{-3}$ & $2.2\times10^{-2}$ & $7.1\times10^{-3}$ \\
\cmidrule(l){2-7}
\multirow{3}{*}{RP7} & $L^1$ & $4.3\times10^{-3}$ & $1.6\times10^{-3}$ & $3.4\times10^{-3}$ & $5.3\times10^{-3}$ & $3.8\times10^{-3}$ \\
 & $L^2$ & $9.9\times10^{-3}$ & $3.7\times10^{-3}$ & $7.9\times10^{-3}$ & $1.3\times10^{-2}$ & $9.6\times10^{-3}$ \\
 & $L^\infty$ & $4.7\times10^{-2}$ & $1.3\times10^{-2}$ & $3.5\times10^{-2}$ & $5.9\times10^{-2}$ & $4.6\times10^{-2}$ \\
\bottomrule
\end{tabular}
\end{table}

The two critical cases test a property that appears nowhere in the objective:
the flow must reach $F_3=-1$ exactly at $x=0^+$, the value at which the cubic
\eqref{eq:sst_cubic} has a double root. In RP1 (Fig.~\ref{fig:rp1}) the exact
solution has four waves: a $1$-rarefaction and a $2$-shock on the low bed, the
critical transition, and a transcritical $2$-rarefaction on the high bed. All
four are recovered with the correct plateau values. The trained traces
$(0.4295,-6.1102)$ and $(0.8845,-2.9655)$ match the exact $(0.4295,-6.1118)$ and
$(0.8889,-2.9530)$ to within $0.5\%$ in depth, and give $F_3=-1.007$ at $x=0^+$
against the exact $-1$: the critical state is attained without ever being
prescribed.

RP2 (Fig.~\ref{fig:rp2}) is the most severe geometric test, with a step twice
the upstream depth: two shocks on the low bed force a transcritical
$2$-rarefaction on the high bed. The four-wave structure is reproduced, and the
trained trace gives $F_3=-0.980$ at $x=0^+$ against the exact $-1$. The free-surface jump, which exceeds the bed jump by
$60\%$, is obtained as $+1.6299$~m against $+1.6122$~m, a $1.1\%$ error that is
the largest of the subset. The interface residuals of RP2 are nevertheless among
the smallest ($|\jump{q}|/|q|=2.8\times10^{-5}$), so the discrepancy lies in the
states that the closure connects, not in the closure itself, and the minimum
depth, $0.45$~m, stays an order of magnitude above $h_{\mathrm{safe}}$.

RP3 and RP4 (Figs.~\ref{fig:rp3} and~\ref{fig:rp4}) probe the robustness of the coupling and of
the parameterisation. In RP3 the transition is supercritical from right to left,
so the exact solution leaves the right state untouched, $U_3=U_r$. This targets a
failure mode specific to a global space--time minimisation, in which spurious
influence could reach regions that no wave has reached. The trained trace at
$x=0^+$ is $(3.0000,-8.0000)$, equal to $U_r$ to four decimal places, and is
held over the whole right half-domain, although the interface residuals carry
ten times the weight of the cell residuals. RP4 stresses the nonlinearity of the
momentum flux: with $F_2=+3.19$ the convective part $q^{2}/h$ dominates the
hydrostatic part by an order of magnitude at a depth of $1$~m, which is where a
parameterisation by $\eta$ rather than $h$ is most exposed. The traces
$(1.0004,9.9948)$ and $(1.1353,8.8098)$ match the exact ones, no oscillation
develops near the step, and the relative $L^1$ error in $\eta$,
$8.6\times10^{-4}$, is among the smallest of the subset.

The $L^\infty$ column of Table~\ref{tab:ab_err} separates the subset into two
groups that do not overlap. The three rarefaction-only cases RP5, RP6 and RP7
(Figs.~\ref{fig:rp5}, \ref{fig:rp6} and~\ref{fig:rp7}) reach $4.1\times10^{-2}$,
$7.2\times10^{-3}$ and $4.7\times10^{-2}$ in $h$, while the four cases containing
a shock lie between $7.8\times10^{-2}$ and $3.0\times10^{-1}$. The separation does
not follow the cell aspect: RP5 has the largest ratio $s_{\max}h_t/h_x=26.4$ and
one of the smallest $L^\infty$ errors, while RP2 has the smallest ratio, $5.4$,
and an error more than three times larger. The two groups differ in the
regularity of the lateral integrand, which is the mechanism identified after
Proposition~\ref{prop:consistency}: the Gauss--Legendre rule on the lateral faces
converges at its nominal rate on the Lipschitz integrand of a rarefaction however
fast the wave crosses the cell, and degrades on the discontinuous integrand of a
shock. Consistently, the free-surface jump of RP5, the largest of the subset, is
reproduced to four decimal places ($+1.6806$ against $+1.6805$~m), and that of
RP7 to the four decimals of the exact value ($+1.1839$~m). The shock smearing seen
in RP1--RP4 is thus a resolution effect of the lateral quadrature, not a defect of
the step treatment.

Finally, RP6 (Fig.~\ref{fig:rp6}) tests the selection of the subcritical root,
which nothing in the objective performs. The trained traces $(4.2907,2.7423)$ and
$(2.7203,4.3257)$ give $\jump{\eta}=-0.5704$~m against the exact $-0.5701$~m, and
both Froude numbers stay below unity across the step, whereas the supercritical
root for the same discharge and head would have produced a shallower, faster
state at $x=0^+$ and a larger negative jump, $\jump{\eta}\approx-1.13$~m. The selection is not forced by the
far field, since the right state itself is supercritical ($F_r=+1.14$). RP6 is
also the most accurate case, with a relative $L^1$ error of $6.9\times10^{-4}$ in
$\eta$.

\FloatBarrier
\subsubsection{Interface fidelity, conservation and cost}
\label{sec:ab_cross}

Table~\ref{tab:ab_d2} reports the traces of the trained network on both sides of
the step. The closure is enforced uniformly: the discharge defect never exceeds
$4.6\times10^{-4}$ of the local discharge and the head defect $1.6\times10^{-4}$
of the local head, with no systematic dependence on the character of the
transition, the depth ratio or the presence of a shock. At the same time the
free-surface jump, which must not vanish, is reproduced to within $0.3\%$ in six
of the seven cases and to $1.1\%$ in RP2, over a range from $-0.5704$ to
$+1.6806$~m that includes a negative jump. The two observations together
distinguish the correct closure from a plausible substitute: a formulation that
replaced the head condition by continuity of the free surface could also return
small $\jump{q}$ and $\jump{H}$, but would necessarily give $\jump{\eta}=0$ and
violate the head condition by the amount bounded in
Theorem~\ref{thm:no_continuous}.

\begin{table}[!htbp]
\centering\footnotesize
\setlength{\tabcolsep}{4pt}
\caption{Interface diagnostic at $t=T$: traces of the trained network at
$x=0^\mp$, to be compared with $U_2$, $U_3$ of Table~\ref{tab:ab_config}; the
discharge and head defects, which must vanish, normalised by the local $|q|$ and
$H$; and the free-surface jump, which must not vanish, against its exact value.}
\label{tab:ab_d2}
\begin{tabular}{lcccccccc}
\toprule
Case & $h(0^-)$ & $u(0^-)$ & $h(0^+)$ & $u(0^+)$ & $|\jump{q}|/|q|$ &
$|\jump{H}|/H$ & $\jump{\eta}_\theta$ [m] & $\jump{\eta}_{\mathrm{exact}}$ [m] \\
\midrule
RP1 & $0.4295$ & $-6.1102$  & $0.8845$ & $-2.9655$ & $4.6\times10^{-4}$ & $1.6\times10^{-4}$ & $+1.4550$ & $+1.4594$ \\
RP2 & $0.7210$ & $-6.6867$  & $1.3508$ & $-3.5687$ & $2.8\times10^{-5}$ & $2.7\times10^{-5}$ & $+1.6299$ & $+1.6122$ \\
RP3 & $2.4781$ & $-9.6882$  & $3.0000$ & $-8.0000$ & $3.3\times10^{-4}$ & $4.1\times10^{-6}$ & $+1.5219$ & $+1.5230$ \\
RP4 & $1.0004$ & $9.9948$   & $1.1353$ & $8.8098$  & $3.2\times10^{-4}$ & $1.5\times10^{-4}$ & $+1.1349$ & $+1.1341$ \\
RP5 & $2.3194$ & $-9.0541$  & $3.0000$ & $-7.0001$ & $1.0\times10^{-5}$ & $2.7\times10^{-5}$ & $+1.6806$ & $+1.6805$ \\
RP6 & $4.2907$ & $2.7423$   & $2.7203$ & $4.3257$  & $8.0\times10^{-5}$ & $3.7\times10^{-6}$ & $-0.5704$ & $-0.5701$ \\
RP7 & $2.0000$ & $12.0000$  & $2.1838$ & $10.9896$ & $1.3\times10^{-5}$ & $7.2\times10^{-6}$ & $+1.1839$ & $+1.1839$ \\
\bottomrule
\end{tabular}
\end{table}

Table~\ref{tab:ab_cons} collects the conservation, boundary and cost
diagnostics. The residual of the mass balance \eqref{eq:mass_balance}, with its
three terms computed independently, stays between $3.5\times10^{-6}$ and
$2.7\times10^{-4}$ of the initial water volume. The free-surface excursion at the
boundaries never exceeds $6.2\times10^{-4}$~m, which confirms that no wave
reached a boundary and that the transmissive condition was never stressed. The
minimum depth, from $0.42$~m (RP1) to $2.72$~m (RP6), remains at least eight
times $h_{\mathrm{safe}}=5\times10^{-2}$~m, so the positivity barrier is inactive
throughout.

\begin{table}[!htbp]
\centering\footnotesize
\setlength{\tabcolsep}{4pt}
\caption{Conservation, boundary and cost diagnostics. Mass: residual of
\eqref{eq:mass_balance} normalised by the initial volume $V_0$. Boundary:
$\max_t|\eta_\theta-\eta_0|$ at $x=-L$ and $x=+L$, which detects a wave reaching
a boundary. $\Lc(\theta)$: objective reached; $N_{\mathrm{eval}}$: number of loss
evaluations; $t_{\mathrm{train}}$: training time on one A100-SXM4-80GB.}
\label{tab:ab_cons}
\begin{tabular}{lcccccccc}
\toprule
Case & $V_0$ [m$^2$] & mass$/V_0$ & boundary $x=-L$ [m] & boundary $x=+L$ [m] &
$\min h_\theta$ [m] & $\Lc(\theta)$ & $N_{\mathrm{eval}}$ & $t_{\mathrm{train}}$ [s] \\
\midrule
RP1 & $66.0$  & $2.2\times10^{-4}$ & $5.8\times10^{-5}$ & $1.5\times10^{-4}$ & $0.420$ & $1.3\times10^{-3}$ & $6682$  & $224$ \\
RP2 & $35.0$  & $6.4\times10^{-5}$ & $6.1\times10^{-4}$ & $9.1\times10^{-6}$ & $0.453$ & $1.5\times10^{-4}$ & $11\,029$ & $370$ \\
RP3 & $42.0$  & $1.9\times10^{-5}$ & $9.8\times10^{-6}$ & $9.5\times10^{-7}$ & $1.000$ & $7.8\times10^{-4}$ & $5139$  & $183$ \\
RP4 & $40.0$  & $2.7\times10^{-4}$ & $3.7\times10^{-4}$ & $7.6\times10^{-5}$ & $1.000$ & $5.2\times10^{-4}$ & $3194$  & $110$ \\
RP5 & $75.0$  & $3.5\times10^{-6}$ & $2.5\times10^{-5}$ & $4.3\times10^{-6}$ & $0.462$ & $1.4\times10^{-3}$ & $3273$  & $114$ \\
RP6 & $136.5$ & $1.4\times10^{-5}$ & $3.2\times10^{-5}$ & $3.2\times10^{-5}$ & $2.720$ & $1.1\times10^{-4}$ & $7454$  & $267$ \\
RP7 & $62.5$  & $2.8\times10^{-5}$ & $3.6\times10^{-6}$ & $4.8\times10^{-5}$ & $0.781$ & $7.0\times10^{-4}$ & $14\,211$ & $443$ \\
\bottomrule
\end{tabular}
\end{table}

The cost is the main limitation. Training takes between $110$ and $443$~s per
case on a single A100, against a fraction of a second for a well-balanced
finite-volume scheme on the same problem, and, although the optimiser settings
are the same in all seven cases (Table~\ref{tab:ab_settings}), the effort is not
predictable from the data: the number of loss evaluations needed varies by a
factor of $4.4$ and the final objective by an order of magnitude. This sensitivity to the optimiser, not
the treatment of the step, is the least satisfactory aspect of the present
implementation. What the subset establishes is that the construction of
Section~\ref{sec:method} solves the right problem: over a discontinuous bed, at
non-zero discharge, in sub-, super- and critical regimes and with either sign of
the free-surface jump, the trained network converges to the correct weak
solution rather than to the continuous one that a single smooth network would be
forced to produce.

\subsection{A smooth-bed control: dam break over a bump}
\label{sec:res_bump}

The last test removes the bed discontinuity altogether: $b=\frac32e^{-x^{2}/2}$
on $[-10,10]$~m with $\eta_0=4$~m for $x<0$ and $2$~m for $x\ge0$, water
initially at rest, transmissive boundaries and $T=1$~s. The detector finds no
jump, the decomposition collapses to a single network ($4\times40$, $5123$
parameters), and no interface condition is imposed: the interface term of the
objective is identically zero at every iteration, not merely small. This is the
empirical counterpart of Proposition~\ref{prop:reduction}: without a bed jump the
additional machinery removes itself and only the source treatment of
Sections~\ref{sec:source}--\ref{sec:correction} remains. Since the bed is
smooth, non-constant and the flow unsteady, the test also checks what the
lake-at-rest suite cannot: that the bed-slope correction, which makes the
equilibrium exact, does not degrade the solution far from equilibrium, a question
that Proposition~\ref{prop:source_consistency} answers only asymptotically. The
settings are $N=8000$, $N_\Gamma=80$, $Q=12$, $a_0=0.9$, $5000$ Adam iterations and
one L-BFGS round of $10^4$ iterations. No closed-form solution exists, and the
reference is the second-order well-balanced MUSCL--HLL scheme of~\ref{app:hll}
\cite{HLL1983,Audusse2004} on $4000$ cells at $\mathrm{CFL}=0.4$, with minmod
reconstruction of $(\eta,u,b)$, hydrostatic reconstruction and two-stage SSP
Runge--Kutta time stepping.

\begin{figure}[!htbp]
  \centering
  \includegraphics[width=0.85\textwidth]{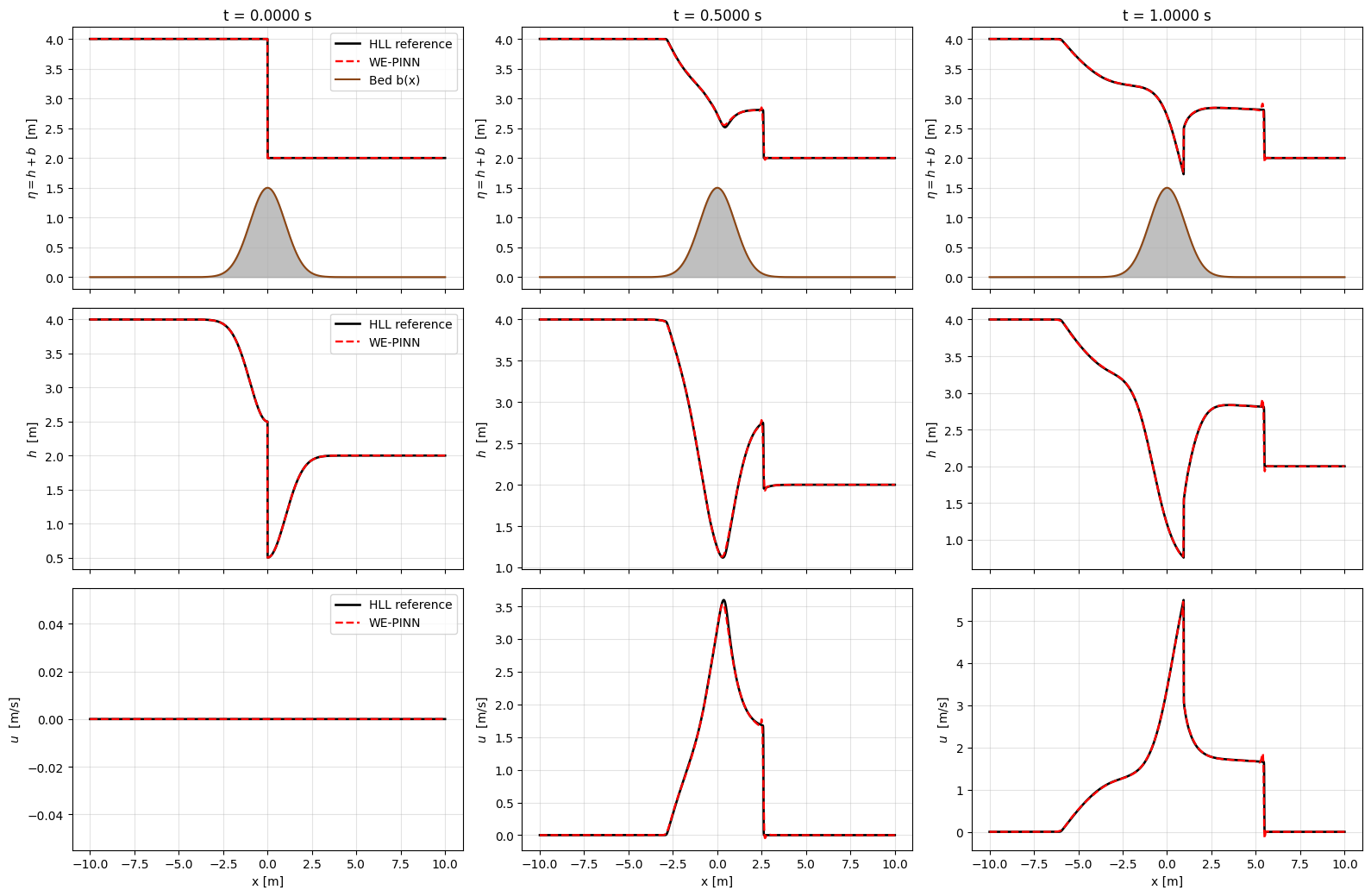}
  \caption{Dam break over a smooth Gaussian bump with no bed jump: free surface
  (top), depth (middle) and velocity (bottom) at $t=0$, $T/2$ and $T$. Black:
  well-balanced MUSCL--HLL reference; red dashed: WE-PINN with a single network
  and no interface condition; grey: bed.}
  \label{fig:bump}
\end{figure}

The solution (Fig.~\ref{fig:bump}) follows the reference through the
rarefaction, over the bump, across the strong acceleration on the lee slope, where
the velocity reaches $5.5$~m\,s$^{-1}$, and across the right-running shock, with
the correct plateau values on both sides. The relative errors at $t=T$
(Table~\ref{tab:bump_err}) are $1.4\times10^{-3}$ in $\eta$ and
$1.5\times10^{-3}$ in $h$ in $L^1$, of the same order as in the cases with an
interface. The $L^\infty$ error, $1.7\times10^{-1}$ in $h$, sits on the shock and is
the lateral-quadrature effect identified in Section~\ref{sec:ab_accuracy}
($s_{\max}h_t/h_x=6.3$). The global mass balance closes to $3.9\times10^{-5}$ of
the initial volume, and the final objective, $6.9\times10^{-5}$, shows that
removing the interface does not make the problem easier to optimise: the
difficulty lies in the shock, not in the coupling.

\begin{table}[!htbp]
\centering\footnotesize
\setlength{\tabcolsep}{4pt}
\caption{Dam break over a smooth bump: relative errors at $t=T$ against the
MUSCL--HLL reference, on $8000$ points.}
\label{tab:bump_err}
\begin{tabular}{lccccc}
\toprule
Norm & $h$ & $u$ & $\eta$ & $q$ & $F$ \\
\midrule
$L^1$      & $1.5\times10^{-3}$ & $9.4\times10^{-3}$ & $1.4\times10^{-3}$ & $9.4\times10^{-3}$ & $1.2\times10^{-2}$ \\
$L^2$      & $9.3\times10^{-3}$ & $4.5\times10^{-2}$ & $8.9\times10^{-3}$ & $4.4\times10^{-2}$ & $6.9\times10^{-2}$ \\
$L^\infty$ & $1.7\times10^{-1}$ & $3.9\times10^{-1}$ & $1.7\times10^{-1}$ & $7.4\times10^{-1}$ & $5.4\times10^{-1}$ \\
\bottomrule
\end{tabular}
\end{table}

\FloatBarrier

\section{Conclusions}\label{sec:conclusion}

We have extended the weak-form space--time control-volume formulation of
WE-PINNs from homogeneous conservation laws to the shallow water balance law
over smooth and discontinuous bathymetry, and made the extension well balanced
by construction. The construction combines three elements. The free-surface
variables, which play the role they play in surface-gradient and pre-balanced
finite-volume schemes, make the momentum flux constant and the geometric source
vanish pointwise at the lake at rest. The source is evaluated over each
control volume by an integration by parts that transfers the derivative from the
network to the bed, so that the residual never differentiates the
approximation. A cell-local shift of the bed slope, chosen so that the quadrature
reproduces the exact bed increment, makes the equilibrium cancellation exact for
any control-volume family, any quadrature order and any nominal bed slope; it
decouples well-balancedness from the mesh and from the accuracy, and even the
existence, of a bed derivative. Bed discontinuities are located from the
bathymetry alone, and the resulting sub-domains, each carrying its own network,
are coupled through conservation of discharge and total head, of which the
first is forced by mass conservation and the second is a closure isolated in a
single term of the objective.

The analysis establishes that the common zeros of the continuous residuals are
the weak solutions of the balance law, that the fully discrete residual vanishes
identically on the lake at rest, so that the equilibrium is an exactly
representable global minimiser of the training objective, that the bed-slope
correction is of the order of the quadrature error it cancels, and that an
approximation with a continuous free surface cannot satisfy the head condition
at a bed step at non-zero discharge, with an explicit lower bound on the defect.
The decomposition at bed jumps is therefore necessary, not a convenience.

The experiments confirm these statements and separate what the construction
guarantees from what the optimiser achieves. On four lake-at-rest
configurations, over smooth, discontinuous, compound and tabulated beds, the
residual operator annihilates the exact equilibrium bitwise or to a round-off
floor identified by its dependence on the working precision, including on cells
that contain an undetected bed jump; the trained network approaches the
equilibrium to $O(10^{-4})$ of the physical scales under an austere budget and to
round-off under an extended one, and the difference between the two is
optimisation error. On seven Riemann problems at a bottom step with closed-form
solutions, covering sub-, super- and critical transitions, the trained network
recovers the stationary transition and the surrounding wave pattern, with
relative $L^1$ errors in the free surface between $6.9\times10^{-4}$ and
$3.5\times10^{-3}$, free-surface jumps between $-0.57$ and $+1.68$~m reproduced,
critical states attained
without being prescribed, and the subcritical root selected without any explicit
mechanism. The remaining $L^\infty$ error is confined to shocks and is explained
by the lateral quadrature of a discontinuous integrand. A dam break over a smooth
bump confirms that, when the bed has no discontinuity, the decomposition and
interface terms vanish and only the source treatment remains.

Several limitations remain. First, no estimate links the trained objective to
the $L^1$ error. Second, the well-balanced guarantee concerns the objective and its
minimiser, not the optimiser. Finally, dry states are not treated. These points
set the directions of future work: a stability framework supporting a
convergence estimate, and the extension of the method to wet--dry fronts.

\section*{Acknowledgments}
The authors gratefully acknowledge the support and computing resources provided
by the Toubkal Supercomputer~\cite{Kissami2025TOUBKAL}
(\url{https://toubkal.um6p.ma}) at UM6P, Morocco, where all computational
experiments were carried out on NVIDIA A100 GPUs.

\section*{Code Availability}
The source code is publicly available at
\href{https://github.com/oubarka-ismail/Well-Balanced_WE-PINNs}{https://github.com/oubarka-ismail/Well-Balanced\_WE-PINNs}.

\appendix

\section{Tabulated bathymetry of configuration WB4}
\label{app:bed}

The bed of WB4 is the piecewise-linear interpolant of the following $28$ points
on $[0,1500]$~m; the initial free surface is $\eta_0\equiv12$~m and the water is
at rest.

\begin{table}[H]
\centering\footnotesize
\setlength{\tabcolsep}{3.5pt}
\begin{tabular}{lcccccccccccccc}
\toprule
$x_j$ [m] & 0 & 50 & 100 & 150 & 250 & 300 & 350 & 400 & 425 & 435 & 450 & 475 & 500 & 505\\
$b_j$ [m] & 0.0 & 0.0 & 2.5 & 5.0 & 5.0 & 3.0 & 5.0 & 5.0 & 7.5 & 8.0 & 9.0 & 9.0 & 9.1 & 9.0\\
\midrule
$x_j$ [m] & 530 & 550 & 565 & 575 & 600 & 650 & 700 & 750 & 800 & 820 & 900 & 950 & 1000 & 1500\\
$b_j$ [m] & 9.0 & 6.0 & 5.5 & 5.5 & 5.0 & 4.0 & 3.0 & 3.0 & 2.3 & 2.0 & 1.2 & 0.4 & 0.0 & 0.0\\
\bottomrule
\end{tabular}
\end{table}

\section{The exact Riemann solver at a bottom step}
\label{app:exact}

The reference solutions of Section~\ref{sec:results_ab} follow the closed-form
construction of \cite{AlcrudoBenkhaldoun2001}, reimplemented here. They are
computed after training, used for comparison only, and never enter the loss. The
construction applies when, and only when, the initial conditions are two constant
states separated by a single bed step at the same abscissa; the solution is then
self-similar in $\xi=x/t$, and the step is the stationary wave at $\xi=0$.

\subsection{Building blocks}
\label{app:exact_blocks}

Away from the step the bed is constant and the system is the homogeneous
shallow water system, so the two states connected by a $k$-wave lie on its
classical wave curves. Writing $c(h)=\sqrt{gh}$, a state $(h,u)$ is joined to
$(h_l,u_l)$ by a $1$-wave and to $(h_r,u_r)$ by a $2$-wave when, respectively,
\begin{equation}
u=\Psi^1(h;h_l,u_l):=
\begin{cases}
u_l+2c(h_l)-2c(h), & h\le h_l\quad\text{(rarefaction)},\\[3pt]
u_l-\sqrt{\tfrac12 g(h-h_l)^{2}\big(\tfrac1h+\tfrac1{h_l}\big)},
 & h>h_l\quad\text{(shock)},
\end{cases}
\label{eq:wave_curve_1}
\end{equation}
\begin{equation}
u=\Psi^2(h;h_r,u_r):=
\begin{cases}
u_r-2c(h_r)+2c(h), & h\le h_r\quad\text{(rarefaction)},\\[3pt]
u_r+\sqrt{\tfrac12 g(h-h_r)^{2}\big(\tfrac1h+\tfrac1{h_r}\big)},
 & h>h_r\quad\text{(shock)} .
\end{cases}
\label{eq:wave_curve_2}
\end{equation}
The associated speeds are, for the $1$-family, the fan
$[\,u_l-c(h_l),\,u-c(h)\,]$ in the rarefaction case and the single ray
$u_l-\sqrt{\tfrac12 g(h+h_l)h/h_l}$ in the shock case, and symmetrically for the
$2$-family. Two states admitting a classical solution (the case
$u_l+2c(h_l)\ge u_r-2c(h_r)$, outside which a vacuum forms) are connected
through the intermediate state $U_m=(h_m,u_m)$ obtained from the scalar equation
$\Psi^1(h_m;h_l,u_l)=\Psi^2(h_m;h_r,u_r)$; we write
$U_m=\mathcal{R}(U_l,U_r)$ for this classical Riemann solution.

Across the step the connection is not a wave curve but the stationary step
transition \eqref{eq:sst_q}--\eqref{eq:sst_H}. Given a trace $(h,u)$ on one side
and a bed increment $\Delta z$ towards the other, the conjugate state
$\mathcal{T}_{\Delta z}(h,u)=(\tilde h,\tilde u)$ is defined by $\tilde u=q/\tilde h$
with $q=hu$ and $\tilde h$ a positive root of the cubic \eqref{eq:sst_cubic}
with $A=H_\star:=h+u^{2}/2g-\Delta z$. By Lemma~\ref{lem:conjugate} such a root
exists if and only if $\alpha:=q^{2}/(2gH_\star^{3})\le4/27$, and the two roots
are given in closed form by
\begin{equation}
\tilde h=H_\star\Big[\tfrac13+\tfrac23\cos\big(\tfrac{\vartheta}{3}-\kappa\big)\Big],
\qquad
\vartheta=\arccos\!\big(1-\tfrac{27}{2}\alpha\big),
\qquad
\kappa=\begin{cases}0, &\text{subcritical branch},\\[2pt]
\tfrac{2\pi}{3}, &\text{supercritical branch},\end{cases}
\label{eq:cubic_roots}
\end{equation}
which coalesce at the critical depth $\tilde h=\tfrac23H_\star=(q^{2}/g)^{1/3}$
when $\alpha=4/27$. The implementation returns the branch matching the Froude
regime of the incoming state, and refuses the transition when $H_\star\le0$ or
$\alpha>4/27$: the head carried by the incoming state is then insufficient to
cross the step, and no solution of the assumed type exists.

\subsection{The five admissible orderings}
\label{app:exact_patterns}

A solution is a sequence of the $1$-wave $W^1$, the $2$-wave $W^2$ and the
transition $\mathrm{SST}$, arranged in order of increasing $\xi$ with the
transition fixed at $\xi=0$. Genuine nonlinearity forbids a wave of the same
family on both sides of the step, except in the resonant patterns (D) and (E)
below, where the second wave of that family is a transcritical rarefaction
attached to the step, and the transition can be crossed only in the
direction of the flow, so only five orderings are admissible. Following
\cite{AlcrudoBenkhaldoun2001} they are, with $U_2$ and $U_3$ the traces at
$\xi=0^-$ and $\xi=0^+$ and $U_m$ the intermediate state of a classical pair,
\begin{equation}
\begin{aligned}
\text{(A)}\quad & U_l \xrightarrow{\;W^1\;} U_2
  \xrightarrow{\;\mathrm{SST_{sub}}\;} U_3 \xrightarrow{\;W^2\;} U_r ,\\[2pt]
\text{(B)}\quad & U_l \xrightarrow{\;W^1\;} U_m \xrightarrow{\;W^2\;} U_2
  \xrightarrow{\;\mathrm{SST_{sup}}\;} U_3 = U_r ,\\[2pt]
\text{(C)}\quad & U_l = U_2 \xrightarrow{\;\mathrm{SST_{sup}}\;} U_3
  \xrightarrow{\;W^1\;} U_m \xrightarrow{\;W^2\;} U_r ,\\[2pt]
\text{(D)}\quad & U_l \xrightarrow{\;W^1\;} U_m \xrightarrow{\;W^2\;} U_2
  \xrightarrow{\;\mathrm{SST_{crit}}\;} U_3 \xrightarrow{\;R^2\;} U_r ,\\[2pt]
\text{(E)}\quad & U_l \xrightarrow{\;W^1\;} U_2
  \xrightarrow{\;\mathrm{SST_{crit}}\;} U_3 \xrightarrow{\;R^1\;} U_m
  \xrightarrow{\;W^2\;} U_r .
\end{aligned}
\label{eq:patterns}
\end{equation}
The five differ in where the transition sits relative to the two waves and in
its character, and each carries its own solvability conditions.

In (A) the transition is subcritical and separates the two waves, so both $U_2$
and $U_3$ must be subcritical and the $1$-wave must lie entirely in $\xi<0$ and
the $2$-wave entirely in $\xi>0$. In (B) both waves lie on the low bed to the
left of a supercritical transition, which leaves the right state untouched,
$U_3=U_r$; this requires $U_r$ itself to be supercritical and directed towards
the step. Configuration (C) is its mirror image, with $U_l$ supercritical and
$U_2=U_l$, and both waves on the far side of the step. In (D) and (E) the flow
is choked: the transition is critical, $|F_3|=1$, and the state $U_3$ is joined
to the far state not by a wave of arbitrary sign but by a \emph{transcritical}
rarefaction attached to the step, whose tail sits at $\xi=0$ and which
accelerates the flow away from the critical point. (D) and (E) are mirror images
of each other under reversal of the flow direction: in both the transcritical
fan lies on the right of the step, of the $2$-family in (D) and of the
$1$-family in (E).

\subsection{Reduction to a scalar equation, and selection}
\label{app:exact_solve}

Each ordering reduces to one scalar unknown, so all five are solved by the same
machinery.

Patterns (B), (C) and (D) are explicit up to a classical Riemann solve. In (B)
the right state is known, $U_3=U_r$, so $U_2=\mathcal{T}_{-\Delta z}(U_r)$ on the
supercritical branch and $U_m=\mathcal{R}(U_l,U_2)$; in (C), symmetrically,
$U_2=U_l$, $U_3=\mathcal{T}_{+\Delta z}(U_l)$ on the supercritical branch and
$U_m=\mathcal{R}(U_3,U_r)$. In (D) criticality on the right of the step fixes
$U_3$ outright: the transcritical fan carries the $2$-invariant of $U_r$, so
$c_3=\tfrac13\big(2c(h_r)-u_r\big)$, $h_3=c_3^{2}/g$ and $u_3=-c_3$, after which
$U_2=\mathcal{T}_{-\Delta z}(U_3)$ and $U_m=\mathcal{R}(U_l,U_2)$.

Patterns (A) and (E) require a root. In (A) the unknown is the depth $h_2$ at
$\xi=0^-$, and the residual to be annihilated states that the conjugate of the
$1$-wave state matches the $2$-wave curve,
\begin{equation}
\Phi_{\mathrm A}(h_2):=
\tilde u\big(h_2\big)-\Psi^2\big(\tilde h(h_2);h_r,u_r\big)=0,
\qquad
(\tilde h,\tilde u)=\mathcal{T}_{+\Delta z}\big(h_2,\Psi^1(h_2;h_l,u_l)\big),
\label{eq:resid_A}
\end{equation}
on the subcritical branch. In (E) the unknown is again $h_2$, and criticality at
$\xi=0^+$ with $h_3=(q^{2}/g)^{1/3}$ turns the head condition into
\begin{equation}
\Phi_{\mathrm E}(h_2):=h_2+\frac{u_2^{2}}{2g}-\frac32
\Big(\frac{(h_2u_2)^{2}}{g}\Big)^{1/3}-\Delta z=0,
\qquad u_2=\Psi^1(h_2;h_l,u_l).
\label{eq:resid_E}
\end{equation}
Both residuals are undefined wherever the conjugate state fails to exist, so
they are not simply bracketed between two guesses. Instead every sign change is
located: $\Phi$ is sampled on a geometric grid of $1500$ points spanning
$[10^{-4},\,20\max(h_l,h_r,1)]$, each bracketing pair of consecutive finite
values of opposite sign is refined by Brent's method to a tolerance of
$10^{-14}$, and a pair straddling the boundary of the domain of definition is
first bisected to locate that boundary. Every root so obtained is a candidate.

A candidate becomes the solution only if it passes the admissibility filters of
its ordering: the Froude conditions at $\xi=0^\pm$ ($|F_2|,|F_3|<1$ in~(A);
$|F|>1$ on the appropriate side in~(B), (C); $|F_3|=1$ in~(D), (E)) and the
ordering of the wave speeds, which must place the $1$-wave, the transition and
the $2$-wave at increasing $\xi$ without overlap. All five patterns are
attempted for every initial condition. If none succeeds no reference is produced and the
network is reported alone; if more than one succeeds the solver reports the
ambiguity explicitly, since \cite{AlcrudoBenkhaldoun2001} establishes no
uniqueness result. No such ambiguity arose in any of the seven configurations of
Section~\ref{sec:results_ab}.

The accepted solution is stored as an ordered list of self-similar segments in
$\xi$ (constant states, and rarefaction fans identified by their Riemann
invariant) and is evaluated at any $(x,t)$ by locating $\xi=x/t$ in that
list, so the reference is available pointwise and needs no interpolation. The
dry-bed limiting solutions of \cite[\S7]{AlcrudoBenkhaldoun2001} are not used;
the only vacuum handled is the degenerate right state of pattern~(E), which does
not occur in the wet-bed suite reported here.

\section{The well-balanced finite-volume reference}
\label{app:hll}

For the dam break over a smooth bump of Section~\ref{sec:res_bump}, which has
no closed-form solution, the reference is a second-order well-balanced
MUSCL--HLL scheme \cite{HLL1983,Audusse2004}, described here in full because
the accuracy of the comparison depends on it. It is computed after training, on the host CPU,
and never enters the loss.

\subsection{Discretisation and boundary treatment}
\label{app:hll_grid}

The domain $[x_{\mathrm L},x_{\mathrm R}]$ is divided into $n_x=4000$ uniform
cells of width $\Delta x$, with centres $x_i=x_{\mathrm L}+(i+\tfrac12)\Delta x$
and cell values $(h_i,q_i)$ and $b_i=b(x_i)$. The bed is sampled once from the
same function the network is given, so the two solvers see the same bathymetry. Two layers of ghost cells
implement the boundary condition, identical to the one imposed on the network:
constant extrapolation of $(h,q,b)$ for the transmissive condition, and
mirroring of $h$ and $b$ with sign reversal of $q$ for a reflecting wall.

\subsection{Reconstruction, flux and source}
\label{app:hll_flux}

The reconstruction is performed on $(\eta,u,b)$ rather than on $(h,q)$, which is
what makes the scheme well balanced: at the lake at rest $\eta$ and $u$ are
constant, every limited slope vanishes, and the reconstruction is exact. With
$\delta\phi_i=\mathrm{minmod}(\phi_i-\phi_{i-1},\,\phi_{i+1}-\phi_i)$ the
interface $i+\tfrac12$ receives
\begin{equation}
\phi_L=\phi_i+\tfrac12\delta\phi_i,
\qquad
\phi_R=\phi_{i+1}-\tfrac12\delta\phi_{i+1},
\qquad \phi\in\{\eta,u,b\} .
\label{eq:muscl}
\end{equation}
Two depths are then formed on each side: the \emph{consistent} depth on the
cell's own bed and the \emph{reconstructed} depth of \cite{Audusse2004} on the
common bed,
\begin{equation}
h^{c}_{L,R}=\big(\eta_{L,R}-b_{L,R}\big)_+ ,
\qquad
b^\star=\max(b_L,b_R),
\qquad
h^\star_{L,R}=\big(\eta_{L,R}-b^\star\big)_+ .
\label{eq:hydro_rec}
\end{equation}
The numerical flux is the HLL flux evaluated on the reconstructed states
$(h^\star_L,h^\star_Lu_L)$ and $(h^\star_R,h^\star_Ru_R)$,
\begin{equation}
\mathcal{F}^{\mathrm{HLL}}=
\begin{cases}
\mathcal{F}_L, & s_L\ge0,\\[2pt]
\dfrac{s_R\mathcal{F}_L-s_L\mathcal{F}_R+s_Ls_R(V_R-V_L)}{s_R-s_L},
 & s_L<0<s_R,\\[6pt]
\mathcal{F}_R, & s_R\le0,
\end{cases}
\label{eq:hll}
\end{equation}
with Einfeldt wave speeds built from the physical and Roe-averaged states,
\begin{equation}
s_L=\min\big(u_L-c(h^\star_L),\,\bar u-\bar c\big),
\qquad
s_R=\max\big(u_R+c(h^\star_R),\,\bar u+\bar c\big),
\label{eq:einfeldt}
\end{equation}
\[
\bar u=\frac{\sqrt{h^\star_L}\,u_L+\sqrt{h^\star_R}\,u_R}
            {\sqrt{h^\star_L}+\sqrt{h^\star_R}},
\qquad
\bar c=\sqrt{\tfrac12 g\big(h^\star_L+h^\star_R\big)} ,
\]
replaced by the dry-front speeds $(u_R-2c_R,\,u_R+c_R)$ or
$(u_L-c_L,\,u_L+2c_L)$ when one side is dry. The mass flux is common to the two
cells, but the momentum flux is not: each side receives the interface value
corrected by the difference between the hydrostatic pressures on its own bed and
on the common bed,
\begin{equation}
\mathcal{F}^{q}_{L}=\mathcal{F}^{q,\mathrm{HLL}}
 +\tfrac12 g\big(h^{c}_L\big)^{2}-\tfrac12 g\big(h^\star_L\big)^{2},
\qquad
\mathcal{F}^{q}_{R}=\mathcal{F}^{q,\mathrm{HLL}}
 +\tfrac12 g\big(h^{c}_R\big)^{2}-\tfrac12 g\big(h^\star_R\big)^{2} .
\label{eq:pressure_corr}
\end{equation}
This asymmetry is the whole of the well-balancing at second order in space: it
is exactly the pressure difference the bed exerts between the two faces, and at
the lake at rest it cancels the flux difference identically. The remaining
bed-slope contribution, which arises only from the variation of $b$ inside a
cell and therefore only at second order, is added as the centred source
\begin{equation}
S_i=-\,g\,\frac{h^{c}_{R,i-1/2}+h^{c}_{L,i+1/2}}{2}\;
\frac{b_{L,i+1/2}-b_{R,i-1/2}}{\Delta x} ,
\label{eq:hll_source}
\end{equation}
evaluated between the two reconstructed faces of cell $i$, so that it too
vanishes when $\eta$ and $u$ are constant.

\subsection{Time integration and cost}
\label{app:hll_time}

Time stepping is the two-stage SSP Runge--Kutta scheme of Heun, with the depth
clipped at zero after each stage and the step size taken as
$\Delta t=\min\big(\mathrm{CFL}\,\Delta x/\max|s|,\;t_{\mathrm{out}}-t\big)$
at $\mathrm{CFL}=0.4$, where $\max|s|$ is the largest Einfeldt speed of the
current sweep, so that the solution is advanced exactly onto each reporting
time. A run of the smooth-bump configuration takes about $3500$ steps and a few
seconds on one CPU core.

\bibliographystyle{elsarticle-num}
\bibliography{references}

\end{document}